\documentclass[a4paper,twoside]{article}

\newcommand\CustomTitle{%
    Review and computational
    comparison of
    adaptive least-squares finite element schemes%
}
\newcommand\CustomShortTitle{%
    Review and comparison of
    adaptive LSFEM
}
\newcommand\CustomAbstract{%
    The convergence analysis for least-squares finite element methods
    led to various adaptive mesh-refinement strategies:
    Collective marking algorithms driven by the built-in a~posteriori error
    estimator or an alternative explicit residual-based error estimator
    as well as a separate marking strategy based on the alternative error
    estimator and an optimal data approximation algorithm.
    This paper reviews and discusses available convergence results.
    In addition, all three strategies are investigated empirically
    for a set of benchmarks examples
    of second-order elliptic partial differential equations
    in two spatial dimensions.
    Particular interest is on
    the choice of the marking and refinement parameters and
    the approximation of the given data.
    The numerical experiments are reproducible using
    the author's software package octAFEM
    available on the platform Code Ocean.
}
\newcommand\CustomKeywords{%
    least-squares finite element method,
    adaptive mesh-refinement,
    alternative a~posteriori error estimation,
    elliptic PDEs,
    separate marking,
    data approximation,
    numerical experiments
}
\newcommand\CustomClassification{%
    65N12, 
    65N15, 
    65N30, 
    65N50, 
    65Y20  
}

\newcommand\AuthorA{Philipp Bringmann}
\newcommand\ShortAuthorA{P.~Bringmann}
\newcommand\AffiliationA{%
    TU Wien,
    Institute of Analysis and Scientific Computing,
    Wiedner Hauptstr.\ 8--10,
    1040 Vienna,
    Austria
}

\usepackage[%
    includehead,%
    inner  = 3cm,%
    outer  = 3cm,%
    top    = 2.4cm,%
    bottom = 3.4cm%
]{geometry}
\usepackage{lmodern}
\usepackage{microtype}
\usepackage{titlesec}
\titleformat{\section}[runin]{\bfseries}{\thesection.}{0.5em}{}[.]
\titleformat{\subsection}[runin]{\bfseries}{\thesubsection.}{0.5em}{}[.]
\titleformat{\subsubsection}[runin]{\bfseries}{\thesubsection.}{0.5em}{}[.]
\renewcommand\maketitle{
    \begin{center}
        {\bfseries\MakeUppercase{\CustomTitle}} \par
        \begin{center}
            \footnotesize
            Dedicated to Professor Leszek F.~Demkowicz on the
            occasion of his 70th birthday
        \end{center}
        {\footnotesize\MakeUppercase{\AuthorA}} \par
        {\footnotesize\AffiliationA}
    \end{center}
    \thispagestyle{plain}
}
\newcommand\CustomAbstractFormat{%
    {
        \footnotesize
        \hspace\parindent
        \textbf{Abstract.}
        \CustomAbstract
        \medskip

        \textbf{Keywords.}
        \CustomKeywords
        \medskip

        \textbf{AMS subject classification.}
        \CustomClassification
    }
}
\usepackage[british]{babel}
\renewenvironment{thebibliography}[1]{%
    \titleformat{\section}{\centering}{}{}{\MakeUppercase}%
    \begin{oldthebibliography}{#1}%
    \footnotesize
    \setlength{\parskip}{0ex}%
    \setlength{\itemsep}{0ex plus 0.25ex}%
}{%
    \end{oldthebibliography}%
}

\usepackage{amsmath}
\usepackage{amsthm}
\theoremstyle{plain}
\newtheorem{theorem}{Theorem}[section]

\theoremstyle{definition}

\theoremstyle{remark}

\makeatletter
\thm@headfont{\scshape}
\makeatother
\makeatletter
\let\thm@indent\indent

\makeatother

\usepackage{hyperref}
\usepackage{pdftexcmds}
\makeatletter
\newcommand{\strequal}[2]{\pdf@strcmp{#1}{#2}==0}
\makeatother
\usepackage{amssymb}
\usepackage{mathtools}
\usepackage[font=footnotesize,labelfont=bf]{caption}
\usepackage[caption=false]{subfig}
\usepackage{enumitem}
\usepackage{algorithm}
\usepackage{algpseudocode}
\algrenewcommand\algorithmicrequire{\textbf{Input:}}
\algrenewcommand\algorithmicensure{\textbf{Output:}}

\usepackage[caption=false]{subfig}
\usepackage{tikz} 
\usetikzlibrary{arrows,positioning,fit,calc,matrix,shapes,spy,patterns,decorations.pathreplacing}
\definecolor{TUblue}{rgb}{0,0.4,0.6}
\definecolor{TUgray}{rgb}{0.3922,0.3882,0.3882}
\definecolor{TUgreen}{rgb}{0,0.4941,0.4431}
\definecolor{TUmagenta}{rgb}{0.7294,0.2745,0.5098}
\definecolor{TUyellow}{rgb}{0.8824,0.5373,0.1333}
\colorlet{HUblue}{TUblue}
\colorlet{HUred}{TUmagenta}
\colorlet{HUsand}{TUyellow}
\colorlet{HUgreen}{TUgreen}
\usepackage{pgfplots}
\pgfplotsset{%
  compat=newest,%
  every axis/.style={scale only axis},%
  grid style={densely dotted, semithick},%
}
\usepackage{pgfplotstable}
\newcommand\drawslopetriangle[4][ST]{
  \pgfplotsextra
  {
    \pgfkeys{/pgf/fpu=true}
    \pgfmathsetmacro\leftcoord{#3};
    \pgfmathsetmacro\rightcoord{10*#3};
    \pgfmathsetmacro\bottomcoord{#4};
    \pgfmathsetmacro\topcoord{10^(#2)*#4};
    \pgfkeys{/pgf/fpu=false}

    \coordinate (#1-BL) at (axis cs:\leftcoord,\bottomcoord);
    \coordinate (#1-BR) at (axis cs:\rightcoord,\bottomcoord);
    \coordinate (#1-TL) at (axis cs:\leftcoord,\topcoord);

    \shadedraw[%
      bottom color = black!20,%
      middle color = black!5,%
      top color    = white,%
      draw         = black%
    ]
      (#1-TL) -- (#1-BL) node[midway, left] {\(#2\)}
      -- (#1-BR) node[midway, below] {\(1\)} -- (#1-TL);
  }
}
\newcommand\drawswappedslopetriangle[4][SST]{
  \pgfplotsextra
  {
    \pgfkeys{/pgf/fpu=true}
    \pgfmathsetmacro\leftcoord{#3/10};
    \pgfmathsetmacro\rightcoord{#3};
    \pgfmathsetmacro\topcoord{#4};
    \pgfmathsetmacro\bottomcoord{10^(-#2)*#4};
    \pgfkeys{/pgf/fpu=false}

    \coordinate (#1-TR) at (axis cs:\rightcoord,\topcoord);
    \coordinate (#1-BR) at (axis cs:\rightcoord,\bottomcoord);
    \coordinate (#1-TL) at (axis cs:\leftcoord,\topcoord);

    \shadedraw[%
      bottom color = black!20,%
      middle color = black!5,%
      top color    = white,%
      draw         = black%
    ]
      (#1-BR) -- (#1-TR) node[midway, right] {\(#2\)}
      -- (#1-TL) node[midway, above] {\(1\)} -- (#1-BR);
  }
}
\newcommand\drawslopetriangleup[4][STU]{
  \pgfplotsextra
  {
    \pgfkeys{/pgf/fpu=true}
    \pgfmathsetmacro\leftcoord{#3};
    \pgfmathsetmacro\rightcoord{10*#3};
    \pgfmathsetmacro\bottomcoord{#4};
    \pgfmathsetmacro\topcoord{10^(#2)*#4};
    \pgfkeys{/pgf/fpu=false}

    \coordinate (#1-BL) at (axis cs:\leftcoord,\bottomcoord);
    \coordinate (#1-BR) at (axis cs:\rightcoord,\bottomcoord);
    \coordinate (#1-TR) at (axis cs:\rightcoord,\topcoord);

    \shadedraw[%
      bottom color    = black!20,%
      middle color = black!5,%
      top color = white,%
      draw         = black%
    ]
      (#1-BL) -- (#1-BR) node[midway, below] {\(1\)}
      -- (#1-TR) node[midway, right] {\(#2\)} -- (#1-BL);
  }
}
\pgfplotsset{
    colormap={parula}{
        rgb=(0.2081,0.1663,0.5292)
        rgb=(0.2116,0.1898,0.5777)
        rgb=(0.2123,0.2138,0.627)
        rgb=(0.2081,0.2386,0.6771)
        rgb=(0.1959,0.2645,0.7279)
        rgb=(0.1707,0.2919,0.7792)
        rgb=(0.1253,0.3242,0.8303)
        rgb=(0.0591,0.3598,0.8683)
        rgb=(0.0117,0.3875,0.882)
        rgb=(0.006,0.4086,0.8828)
        rgb=(0.0165,0.4266,0.8786)
        rgb=(0.0329,0.443,0.872)
        rgb=(0.0498,0.4586,0.8641)
        rgb=(0.0629,0.4737,0.8554)
        rgb=(0.0723,0.4887,0.8467)
        rgb=(0.0779,0.504,0.8384)
        rgb=(0.0793,0.52,0.8312)
        rgb=(0.0749,0.5375,0.8263)
        rgb=(0.0641,0.557,0.824)
        rgb=(0.0488,0.5772,0.8228)
        rgb=(0.0343,0.5966,0.8199)
        rgb=(0.0265,0.6137,0.8135)
        rgb=(0.0239,0.6287,0.8038)
        rgb=(0.0231,0.6418,0.7913)
        rgb=(0.0228,0.6535,0.7768)
        rgb=(0.0267,0.6642,0.7607)
        rgb=(0.0384,0.6743,0.7436)
        rgb=(0.059,0.6838,0.7254)
        rgb=(0.0843,0.6928,0.7062)
        rgb=(0.1133,0.7015,0.6859)
        rgb=(0.1453,0.7098,0.6646)
        rgb=(0.1801,0.7177,0.6424)
        rgb=(0.2178,0.725,0.6193)
        rgb=(0.2586,0.7317,0.5954)
        rgb=(0.3022,0.7376,0.5712)
        rgb=(0.3482,0.7424,0.5473)
        rgb=(0.3953,0.7459,0.5244)
        rgb=(0.442,0.7481,0.5033)
        rgb=(0.4871,0.7491,0.484)
        rgb=(0.53,0.7491,0.4661)
        rgb=(0.5709,0.7485,0.4494)
        rgb=(0.6099,0.7473,0.4337)
        rgb=(0.6473,0.7456,0.4188)
        rgb=(0.6834,0.7435,0.4044)
        rgb=(0.7184,0.7411,0.3905)
        rgb=(0.7525,0.7384,0.3768)
        rgb=(0.7858,0.7356,0.3633)
        rgb=(0.8185,0.7327,0.3498)
        rgb=(0.8507,0.7299,0.336)
        rgb=(0.8824,0.7274,0.3217)
        rgb=(0.9139,0.7258,0.3063)
        rgb=(0.945,0.7261,0.2886)
        rgb=(0.9739,0.7314,0.2666)
        rgb=(0.9938,0.7455,0.2403)
        rgb=(0.999,0.7653,0.2164)
        rgb=(0.9955,0.7861,0.1967)
        rgb=(0.988,0.8066,0.1794)
        rgb=(0.9789,0.8271,0.1633)
        rgb=(0.9697,0.8481,0.1475)
        rgb=(0.9626,0.8705,0.1309)
        rgb=(0.9589,0.8949,0.1132)
        rgb=(0.9598,0.9218,0.0948)
        rgb=(0.9661,0.9514,0.0755)
        rgb=(0.9763,0.9831,0.0538)
    }
}

\newcommand\vvvert{|\mkern-1.5mu|\mkern-1.5mu|}
\DeclareMathOperator{\conv}{conv}
\DeclareMathOperator{\ddiv}{div}
\DeclareMathOperator{\ccurl}{curl}
\DeclareMathOperator{\osc}{osc}
\DeclareMathOperator{\interior}{int}
\newcommand\R{\ensuremath{\mathbb R}}
\newcommand\N{\ensuremath{\mathbb N}}
\newcommand\T{\ensuremath{\mathcal T}}
\newcommand\E{\ensuremath{\mathcal E}}
\newcommand\LS{\textup{LS}}
\newcommand\NAT{\textup{N}}

\newcommand\COL{\textup{C}}
\newcommand\SEP{\textup{S}}

\colorlet{RED}{red}

\begin{document}
    \maketitle
    %
    %
    \CustomAbstractFormat

\section{Introduction}

Least-squares finite element methods (LSFEMs)
are highly popular discretisation schemes
for partial differential equations.
One key feature is their
built-in a~posteriori error estimation
which renders this class of methods well-suited
for adaptive mesh-refining algorithms.
One of the first adaptive algorithms for LSFEMs
has been proposed by Jiang and Carey
\cite{10.1002/nme.1620240308}.
The theoretical basis relies on the equivalence
of the least-squares functional with the error
in the standard Sobolev norm \cite{MR1615154}
resp.\
the equality with the error
in the norm induced by the least-squares functional
\cite{MR1793582}.
This property transfers
to LSFEMs for regularised \(H^{-1}\) loads
up to an oscillation term \cite{MR4425909}.
A particular scaling of the residuals enables
the estimation of the contributions to the underlying norm
separately \cite{MR2671052}.
Further algorithmic contributions
deal with the iterative solution by algebraic multigrid
\cite{MR2765484} and parallelisation \cite{MR2896812}.

The built-in a~posteriori error estimation
and adaptive mesh-refinement for standard LSFEMs
have been established and investigated for a multitude
of problems.
The following non-exhaustive list
illustrates the variety of applications.
Adaptive LSFEMs in computational fluid mechanics
deal with
the shallow water equations
\cite{MR2139398,DanischDissertation},
coupled Stokes-Darcy flow
\cite{MR2783231,MR3343602},
viscoelastic fluids
\cite{10.1016/j.jnnfm.2009.02.004},
interface problems
\cite{MR3817763},
and
fluid-structure interaction
\cite{KayserHeroldDissertation}.
In computational solid mechanics,
adaptive LSFEMs have been investigated
for linear elasticity
\cite{MR2100304,MR2084237},
elasto-plasticity
\cite{MR2285860},
and
the Signorini contact problem
\cite{MR2551156,MR3580409},
Further applications include
convection-diffusion problems
\cite{MR1638080},
parabolic problems
\cite{MR1885710,MR3103833,MR4242919,MR4216839},
hyperbolic problems
\cite{OlsonDissertation},
the transport equation
\cite{MR4087177,MR4127415},
the Poisson-Boltzmann equation
\cite{MR2870014},
Maxwell and Helmholtz equation
\cite{MR3820383},
convex energy minimisation
\cite{MR4433562},
elliptic equations in nondivergence form
\cite{MR4173220},
and
the obstacle problem \cite{MR4050087}.

However, the advances in the convergence analysis with rates
for adaptive FEMs in the past 15 years
seem not to be applicable to this class of methods.
This is because the built-in estimator lacks prefactors
in terms of the mesh-size inhibiting all known arguments
for the proof of a local reduction of this estimator.
In order to overcome this,
an alternative explicit residual-based error estimator
for an adaptive mesh-refining algorithm with optimal convergence rates
is developed for the Poisson model problem in~\cite{MR3296614}
and for further linear model problems
in~\cite{MR3599566,MR3715170,MR3757107}.
The known convergence results for \(h\)-adaptive LSFEMs
are summarised and extended
in the thesis~\cite{BringmannDissertation}
for the Poisson model problem, the Stokes equations,
and the linear elasticity equations
with discretisation of arbitrary polynomial degree and
mixed boundary conditions in three spatial dimensions.
All these algorithms employ a separate marking strategy
with a quasi-optimal data approximation algorithm
\cite{MR3719030}.

The negligence of the divergence contribution to the flux error
allows for a collective marking strategy driven
by the (slightly) modified alternative error estimator
\cite{MR4011536,MR4271577}.
This guarantees optimal convergence rates in terms
of the energy error plus the
\(L^2\) error of the flux variable.

The analysis of the alternative error estimator in \cite{BringmannDissertation}
and in \cite{MR4011536,MR4271577}
requires the exact solution of the linear system of the FEM.
However, a modified adaptive algorithm with collective marking
in \cite{MR4280291} allows for an iterative solver
leading to optimal convergence rates with respect to
the overall computational costs
for a standard adaptive FEM.

While the plain convergence of adaptive LSFEM driven by the
built-in estimator has recently been shown in
\cite{MR4138307,MR4216839} for any bulk parameter
\(0 < \theta \leq 1\),
the only Q-linear convergence result in \cite{MR3671598}
requires a \emph{sufficiently large} bulk parameter
\(0 \ll \theta < 1 \).
This contrasts the established convergence analysis in
\cite{MR3170325,MR3719030} asserting optimal rates
for \emph{sufficiently small} bulk parameters.
The investigation of the bulk parameter is one goal of this paper.

Most but not all of the convergence results in
Section~\ref{sec:algorithms} below
hold for discretisations with
arbitrary fixed polynomial degree.
For the sake of concise statements,
the presentation in this paper restricts
to the lowest-order case.

Besides the theoretical review,
this paper provides an experimental investigation
of the performance of three adaptive LSFEM
based on different error estimates
applied to multiple benchmark problems.
The influence of the chosen marking and refinement parameters
is examined.
A benchmark problem with a scalable microstructure in the right-hand side
with exact integration
allows the investigation of the data approximation employed
in the adaptive algorithm with separate marking.
Another focus is put on the performance of the implementation.

Beforehand,
Section~\ref{sec:triangulation} introduces
the notation for the triangulations
and their adaptive refinement
and
Section~\ref{sec:LSFEM}
presents the least-squares discretisation
of the Poisson model problem.
The subsequent Section~\ref{sec:algorithms}
presents the three investigated adaptive LSFEM
algorithms and recalls the theoretical convergence results.
The first subsection of Section~\ref{sec:experiments}
discusses some aspects of the implementation such
as the employed numerical quadrature.
The Subsections~\ref{subsec:Lshape}--\ref{subsec:interface}
present the results of the experiments.
This paper ends with a conclusion in Section~\ref{sec:conclusion}.

\section{Triangulations and refinement}
\label{sec:triangulation}

Given a  bounded polygonal Lipschitz domain
$\Omega\subset\mathbb{R}^2$,
finite element discretisations typically
base on shape-regular triangulations
of \(\Omega\) into closed triangles \cite{MR2353951}.
Let $\T_0$ be an initial triangulation
of $\Omega$.
Note that the initial condition on \(\T_0\) from
\cite[Sect.~4]{MR2353951} is not required in 2D
\cite{MR3097045}.
Given a set \(\mathcal{M}_0 \subseteq \T_0\) of marked triangles,
the refinement algorithm from \cite[Sect.~6]{MR2353951}
creates the smallest regular refinement $\T_1$ of $\T_0$
such that all triangles in
$\mathcal{M}_0 \subseteq \T_0 \setminus \T_1$ are refined.
The algorithm employs the newest-vertex bisection (NVB) from
\cite{MR1311687,MR1329875,MR1475530}.
This defines the concept of a one-level refinement
\cite[Sect.~2]{MR2353951}
leading to the set of admissible triangulations
\begin{align*}
    \mathbb{T}
    &\coloneqq
    \{ \T_\ell \text{ regular triangulation of } \Omega
        \text{  into closed triangles} \;:\\
    &\qquad
    \exists \ell \in \mathbb{N}_0 \exists \T_0,\T_1,\dots,\T_\ell
    \text{ successive one-level refinements in the sense}\\
    &\qquad
    \text{that } \T_{j+1} \text{ is a one-level refinement of }
    \T_j \text{ for } j=0,1,\dots,\ell-1\}.
\end{align*}
It contains the finite subsets of triangulations
with at most \(N \in \N\) additional triangles
\[
    \mathbb{T}(N)
    \coloneqq
    \{
        \T \in \mathbb{T}
        \;:\;
        \vert \T \vert - \vert \T_0 \vert \leq N
    \}.
\]

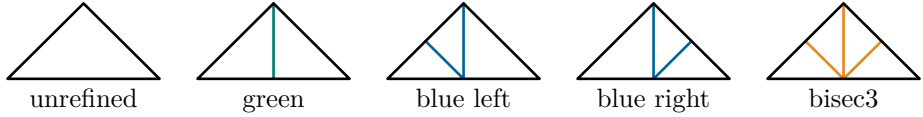
\begin{figure}
    \centering
    \begin{tikzpicture}[%
    line cap = round,%
    line join = round,%
    line width = 1pt]

    \draw (0, 0) -- node[midway, below] {unrefined\vphantom{fg}} (2, 0)
    -- (1, 1) -- cycle;

    \begin{scope}[xshift = 2.5cm]
        \draw[HUgreen] (1, 1) -- (1, 0);
        \draw (0, 0) -- node[midway, below] {green\vphantom{fg}} (2, 0)
        -- (1, 1) -- cycle;
    \end{scope}

    \begin{scope}[xshift = 5cm]
        \draw[HUblue] (1, 1) -- (1, 0) -- (0.5, 0.5);
        \draw (0, 0) -- node[midway, below] {blue left\vphantom{fg}} (2, 0)
        -- (1, 1) -- cycle;
    \end{scope}

    \begin{scope}[xshift = 7.5cm]
        \draw[HUblue] (1, 1) -- (1, 0) -- (1.5, 0.5);
        \draw (0, 0) -- node[midway, below] {blue right\vphantom{fg}} (2, 0)
        -- (1, 1) -- cycle;
    \end{scope}

    \begin{scope}[xshift = 10cm]
        \draw[TUyellow] (1, 1) -- (1, 0) -- (1.5, 0.5) (1, 0) -- (0.5, 0.5);
        \draw (0, 0) -- node[midway, below] {bisec3\vphantom{fg}} (2, 0)
        -- (1, 1) -- cycle;
    \end{scope}

\end{tikzpicture}
    \caption{Possible refinement of a triangle in 2D newest-vertex bisection.}
    \label{fig:refinement}
\end{figure}

In 2D, the one-level refinements result
in one of the five possible refinements
as displayed in Figure~\ref{fig:refinement}
for each triangle in \(\T\).
For a sequence of successively refined meshes
$\T_\ell$ for $\ell \in \mathbb{N}_0$,
the mesh-closure estimate bounds
the number of newly created triangles
\cite{MR1311687,MR1475530}
\begin{align*}
    |\T_\ell| - |\T_0|
    \leq
    C_\textup{NVB}
    \sum_{j=0}^{\ell-1} |\mathcal{M}_{j}|.
\end{align*}

The process of marking triangles for bisection employs
an error indicator \(\eta(\T, T) \in \R\) for each \(T \in \T\).
For any subset \(\mathcal{M} \subseteq \T\),
abbreviate the corresponding contributions
\(\eta^2(\T, \mathcal{M}) \coloneqq \sum_{T \in
\mathcal{M}} \eta^2(\T, T)\)
and \(\eta^2(\T) \coloneqq \eta^2(\T, \T)\).
Given a bulk parameter \(0 < \theta \leq 1\),
the D\"orfler marking from \cite{MR1393904}
selects a subset \(\mathcal{M} \subseteq \T\)
with minimal cardinality \cite{MR4136545}
according to the criterion
\begin{equation}
    \label{eq:bulk_criterion}
    \theta\,
    \eta^2(\T)
    \leq
    \eta^2(\T, \mathcal{M}).
\end{equation}

For a triangulation $\T \in \mathbb{T}$,
let  $\E(T)$ denote the set of all edges of a triangle $T\in\T$.
Let $\E$ denote the set of edges and
$\E(\partial\Omega)$ the edges on the boundary $\partial\Omega$.
Each triangle $T$ has an outer unit normal vector $\nu_T$
and the orientation of each edge $E \in \E$ is fixed
by one of the two possible orientations
of the unit normal vector $\nu_E$
(and unit tangential $\tau_E$).
For an interior edge
$E = \partial T_+ \cap \partial T_- \in \E(\Omega)$
shared by two triangles \(T_+, T_- \in \T\)
ensuring that $\nu_E = \nu_{T_+} = -\nu_{T_-}$.
Accordingly $[w_h]_E \coloneqq (w_h|_{T_+} - w_h|_{T_-})_E$
defines the jump of any piecewise $H^1$ function $w_h$.
Let $\omega_E \coloneqq \interior(T_+ \cup T_-)$ denote the
patch of the interior edge \(E\).
Along the boundary $E \subset \partial\Omega$,
the jump $[w_h]_E \coloneqq w_h|_E $ is the trace of
$w_h$ on the unique triangle \(T_+ \in \T\) with \(E \in \E(T_+)\)
and set $\omega_E \coloneqq \interior(T_+)$.

Throughout the paper,
\(A \lesssim B\) abbreviates the relation \(A \leq C\,B\)
with a positive generic constant \(0 < C\)
which solely depends on the initial triangulation \(\T_0\),
but is independent of the underlying
piecewise constant mesh-size function
\(h_\T \in P_0(\T)\) with
\(h_\T \vert_T \coloneqq h_T \coloneqq \vert T \vert^{1/2}\)
for \(T \in \T \in \mathbb{T}\).
The context-sensitive measure \(\vert \bullet \vert\)
denotes not only
the Lebesgue measure of Lebesgue sets in \(\R^2\),
but also
the modulus of real numbers,
the cardinality of finite sets,
and
the Euclidian norm of vectors in \(\R^2\).

\section{LSFEM for the Poisson model problem}
\label{sec:LSFEM}

For a right-hand side $f \in L^2(\Omega)$
on the polygonal Lipschitz domain \(\Omega \subset \R^2\),
the first-order system formulation of the Poisson model problem
seeks $(p,u) \in H(\ddiv,\Omega) \times H^1_0(\Omega)$ with
\begin{equation}
    \label{eq:PMP}
    f + \ddiv p = 0
    \quad\text{and}\quad
    p - \nabla u = 0
    \quad\text{in } \Omega.
\end{equation}
This paper employs standard notation for Sobolev and Lebesgue spaces
\(H^1(\Omega)\), \(H(\ddiv, \Omega)\), and \(L^2(\Omega)\).
Appropriate subscripts designate their usual norms
\(\Vert \bullet \Vert_{H^1(\Omega)}\),
\(\Vert \bullet \Vert_{H(\ddiv, \Omega)}\),
and \(\Vert \bullet \Vert_{L^2(\Omega)}\).

Let \(\T \in \mathbb{T}\) denote
a regular triangulation of \(\Omega\)
into closed triangles.
The lowest-order Raviart-Thomas function space
$RT_0(\T)\subset H(\ddiv,\Omega)$ and
the conforming piecewise polynomials of first order
$S_0^1(\T)\subset H_0^1(\Omega)$
allow for a unique discrete minimiser
$(p_\LS, u_\LS) \in RT_0(\T) \times S^1_0(\T)$
of the least-squares functional
\[
    LS(f;q_\LS,v_\LS)
    \coloneqq
    \Vert f + \ddiv q_\LS \Vert_{L^2(\Omega)}^2
    +
    \Vert q_\LS - \nabla v_\LS \Vert_{L^2(\Omega)}^2
\]
over all $(q_\LS,v_\LS) \in RT_0(\T) \times S^1_0(\T)$.
The fundamental equivalence of the homogeneous least-squares functional
\cite[Lem.~4.3]{MR0461948}
\begin{equation}
    \label{eq:equivalence}
    LS(0; q, v)
    \approx
    \Vert q \Vert_{H(\ddiv, \Omega)}^2
    +
    \Vert \nabla v \Vert_{L^2(\Omega)}^2
    \quad
    \text{for all }
    (q, v) \in H(\ddiv, \Omega) \times H^1_0(\Omega)
\end{equation}
provides well-posedness of the LSFEM.
This ensures convergence towards
the solution $(p,u)$ of \eqref{eq:PMP}
in the case of (quasi-)uniform mesh-refinement
\[
    \Vert p - p_\LS \Vert_{H(\ddiv,\Omega)}^2
    +
    \Vert \nabla(u - u_\LS) \Vert_{L^2(\Omega)}^2
    \lesssim
    \inf_{\substack{q_\LS \in RT_0(\T) \\ v_\LS \in S^1_0(\T)}}
    \hspace{-3pt}
    \big(
        \Vert p - q_\LS \Vert_{H(\ddiv,\Omega)}^2
        +
        \Vert \nabla(u - v_\LS) \Vert_{L^2(\Omega)}^2
    \big).
\]

Let \(P_k(\T)\) denote the space of piecewise polynomials
with respect to the triangulation \(\T\) and
$\Pi: L^2(\Omega) \to P_0(\T)$
the $L^2$-orthogonal projection onto \(P_0(\T)\).
The piecewise constant approximation
of some \(f \in L^2(\Omega)\) leads to
the data oscillations
\begin{equation}
    \label{eq:osc}
    \osc^2(f, \T)
    \coloneqq
    \sum_{T \in \T}
    h_T^2\,
    \Vert (1 - \Pi) f \Vert_{L^2(T)}^2.
\end{equation}

\section{Three adaptive algorithms}
\label{sec:algorithms}

The following subsections introduce each of the adaptive LSFEMs,
including the employed error estimators and the adaptive algorithm.
They recall the theoretical convergence results.

\subsection{Natural adaptive LSFEM}
\label{subsec:NALSFEM}

The contributions to
the built-in a~posteriori error estimator
\begin{equation}
    \label{eq:estimator_natural}
    \eta_\NAT^2(\T, T)
    \coloneqq
    \Vert f + \ddiv p_\LS \Vert_{L^2(T)}^2
    +
    \Vert p_\LS - \nabla u_\LS \Vert_{L^2(T)}^2
\end{equation}
sum up to the least-squares functional
\[
    \eta_\NAT^2(\T)
    \coloneqq
    \sum_{T \in \T}
    \eta_\NAT^2(\T, T)
    =
    LS(f; p_\LS, u_\LS).
\]
The fundamental equivalence~\eqref{eq:equivalence} ensures
that this estimator is reliable and efficient
even in the case of the inexact solution of the discrete problem
\[
    LS(f; q_\LS, v_\LS)
    \approx
    \Vert p - q_\LS \Vert_{H(\ddiv, \Omega)}^2
    +
    \Vert \nabla(u - v_\LS) \Vert_{L^2(\Omega)}^2.
\]
Moreover, the built-in error estimator is
even asymptotically exact
with respect to
the norm on \(H(\ddiv, \Omega) \times H^1_0(\Omega)\)
abbreviated by
\[
    \vvvert (q, v) \vvvert^2
    \coloneqq
    \Vert q \Vert_{H(\ddiv, \Omega)}^2
    +
    \Vert \nabla v \Vert_{L^2(\Omega)}^2
    \quad\text{for }
    (q, v) \in
    H(\ddiv, \Omega) \times H^1_0(\Omega).
\]
\begin{theorem}[{\cite[Thm.~3.1]{MR3820383}}]
    \label{thm:asymptotic}
    For all \(\varepsilon > 0\),
    there exists \(\delta > 0\)
    such that every \(\T \in \mathbb{T}\)
    with \(\max_{T \in \T} \operatorname{diam}(T) \leq \delta\)
    satisfies
    \[
        (1 - \varepsilon)
        \vvvert (p - p_\LS, u - u_\LS) \vvvert^2
        \leq
        LS(f; p_\LS, u_\LS)
        \leq
        (1 + \varepsilon)
        \vvvert (p - p_\LS, u - u_\LS) \vvvert^2.
    \]
\end{theorem}

Theorem~\ref{thm:asymptotic}
applies to standard conforming discretisations
of any order and various applications \cite{MR3820383}.
However,
an earlier asymptotic exactness result in
\cite{MR3372049,MR3452849}
relies on an unbalanced discretisation
of the two variables \(p_\LS\) and \(u_\LS\).

\begin{algorithm}[t]
    \caption{NALSFEM (natural adaptive LSFEM)}
    \label{alg:NALSFEM}
    \begin{algorithmic}
        \Require
        regular triangulation $\T_0$ and bulk parameter $0 < \theta \leq 1$.
        \For{$\ell = 0,1,2,\dots$}
        \State \textbf{Solve} LSFEM with respect to triangulation
        $\T_\ell$ for solution $(p_\ell, u_\ell)$.
        \State \textbf{Compute} $\eta_\NAT(\T_\ell, T)$
        from~\eqref{eq:estimator_natural}
        for all \(T \in \T_\ell\).
        \State \textbf{Mark}  minimal subset
        $\mathcal{M}_\ell \subseteq \T_\ell$
        by the D\"{o}rfler criterion~\eqref{eq:bulk_criterion}
        for \(\eta \equiv \eta_\NAT\).
        \State \textbf{Refine} $\T_\ell$ to $\T_{\ell + 1}$
        by NVB
        such that $\mathcal{M}_\ell\subseteq \T_\ell \setminus \T_{\ell+1}$.
        \EndFor
        \Ensure
        sequence of triangulations $\T_\ell$ with $(p_\ell, u_\ell)_{\ell}$
        and $\eta_\NAT(\T_\ell)$ for $\ell \in \mathbb{N}_0$.
    \end{algorithmic}
\end{algorithm}

The collective marking for the natural estimator \(\eta_\NAT\)
from~\eqref{eq:estimator_natural}
results in the adaptive Algorithm~\ref{alg:NALSFEM} (NALSFEM).
Independently of the choice
of the bulk parameter \(\theta\),
NALSFEM creates a convergent sequence
of discrete solutions \((p_\ell, u_\ell)_{\ell}\)
for \(\ell \in \N_0\).%
\begin{theorem}[{\cite[Thm.~2]{MR4138307},
    \cite[Thm.~3.3]{MR4216839}}]
    \label{thm:plain_convergence}
    For all \(0 < \theta \leq 1\),
    the output \((p_\ell, u_\ell)_\ell\) of NALSFEM satisfies
    \[
        \Vert p - p_\ell \Vert_{H(\ddiv, \Omega)}^2
        +
        \Vert \nabla (u - u_\ell) \Vert_{L^2(\Omega)}^2
        \to
        0
        \quad\text{as } \ell \to \infty.
    \]
\end{theorem}
The proofs in \cite{MR4138307,MR4216839}
employ the plain convergence framework
from \cite{MR2832786}
under mild assumptions on
the partial differential equation,
the marking strategy, and the mesh refinement.
It applies to higher-order discretisations as
well as to
more general marking criteria,
e.g., the maximum marking strategy or
the equilibrium marking strategy \cite[Sect.~2.6]{MR4138307}.

If the NVB in the step \textbf{Refine} ensures
the bisection of each edge of
the marked triangles in \(\mathcal{M}_\ell\),
then NALSFEM
converges Q-linearly in the following sense.
\begin{theorem}[{\cite[Thm.~4.1]{MR3671598}}]
    \label{thm:convergence_natural}
    Assume that the initial triangulation is sufficiently fine
    in that \(f = \Pi_{L+1} f\) is resolved exactly
    on the level \(\T_{L+1}\).
    There exist
    a minimal bulk parameter \(0 < \Theta_0 < 1\),
    a reduction factor \(0 < \varrho < 1\),
    and a constant \(0 < \Lambda < \infty\)
    such that, for all \(\Theta_0 \leq \theta \leq 1\),
    the modified estimator
    \[
        \widehat\eta_\NAT^2(\T_\ell)
        \coloneqq
        LS(f; p_\ell, u_\ell)
        +
        \Lambda\,
        \Vert
            (1 - \Pi_\ell) p_\ell
        \Vert_{L^2(\Omega)}^2
    \]
    with the output \((p_\ell, u_\ell)_\ell\) of NALSFEM satisfies
    \[
        \widehat\eta_\NAT^2(\T_{\ell+1})
        \leq
        \varrho\,
        \widehat\eta_\NAT^2(\T_\ell)
        \quad\text{for all } \ell = L, L+1, \ldots
    \]
\end{theorem}

The key difficulty in the proof of convergence with
rates as in Theorem~\ref{thm:convergence_natural}
consists of the reduction of the natural estimator
\(\eta_\NAT\) on refined triangles.
Within the frameworks \cite{MR3170325,MR3719030},
this relates to axiom~(A2)
for \(0 < \varrho < 1\) and \(0 < \Lambda\) such that
\begin{equation}
    \label{eq:reduction}
    \eta_\NAT(\T_{\ell+1}, \T_{\ell+1} \setminus \T_\ell)
    \leq
    \rho\,
    \eta_\NAT(\T_\ell, \T_\ell \setminus \T_{\ell+1})
    +
    \Lambda\,
    (LS(0; p_{\ell+1} - p_\ell, u_{\ell+1} - u_\ell))^{1/2}.
\end{equation}
The lack of prefactors in terms of the mesh-size
prevent the usual arguments for the proof of
\eqref{eq:reduction} for \(\eta_\NAT\),
cf.\ \cite{MR2324418,MR2421046,MR3170325}.
The earlier contributions \cite{MR1615154} and \cite{MR2650218}
to the convergence analysis of adaptive LSFEMs
prove the strict reduction,
for \(0 < \varrho < 1\),
\[
    \eta_\NAT^2(\T_{\ell+1})
    \leq
    \varrho\,
    \eta_\NAT^2(\T_\ell)
\]
in each refinement step under the explicit assumption
of a reduction property as \eqref{eq:reduction}
(called \emph{local saturation} in \cite{MR2650218}).
Note that both works  \cite{MR1615154,MR2650218}
employ a nonstandard marking routine and
include severe restrictions on the refinement region
(resp. on the shape of the domain \(\Omega\)).

It turns out that the linear convergence
(for small bulk parameter \(\theta\))
already implies the optimal convergence rate.

\begin{theorem}[{\cite[Prop.~15]{MR4138307}}]
    \label{thm:linear_implies_rate}
    There exists a maximal bulk parameter
    \(0 < \theta_0 < 1\) such that,
    for every \(0 < \theta \leq \theta_0\),
    the following implication holds.
    If the output \((p_\ell, u_\ell)_\ell\)
    of NALSFEM satisfies linear convergence
    with reduction factor \(0 < \varrho < 1\),
    for all \(\ell, m \in \N_0\),
    \begin{equation}
        \label{eq:linear_convergence_natural}
        \eta_\NAT(\T_{\ell + m})
        \lesssim
        \varrho^{m}\,
        \eta_\NAT(\T_\ell),
    \end{equation}
    then \((p_\ell, u_\ell)_\ell\) even converges
    with the optimal rate, i.e.,
    \[
        \sup_{\ell \in \N_0}
        (1 + \vert \T_\ell \vert - \vert \T_0 \vert)^s
        \eta_\NAT(\T_\ell)
        \approx
        \sup_{N \in \N_0}
        (1 + N)^s
        \min_{\T \in \mathbb{T}(N)}
        \eta_\NAT(\T).
    \]
\end{theorem}

This result solely provides a \emph{sufficient} condition
for optimal convergence rates.
However, the linear convergence~\eqref{eq:linear_convergence_natural}
in the case of a small bulk parameter
\(0 < \theta < \theta_0 \ll 1\)
remains an open question.
In particular,
the assumptions of a sufficiently large
\(0 \ll \Theta_0 \leq \theta\)
in Theorem~\ref{thm:convergence_natural}
and of a sufficiently small \(\theta \leq \theta_0\)
in Theorem~\ref{thm:linear_implies_rate}
appear incompatible.

\subsection{Alternative adaptive least-squares FEM with collective marking}
\label{subsec:CALSFEM}

In order to enable the reduction property
of the form~\eqref{eq:reduction}, the
convergence analysis with rates for least-squares FEMs
in \cite{MR3296614,BringmannDissertation,MR4011536,MR4271577}
introduces alternative explicit a~posteriori error estimators
in terms of the constitutive residual
\begin{equation}
    \label{eq:estimator_separate}
    \begin{split}
        \eta_\SEP^2(\T, T)
        &\coloneqq
        h_T^2\,
        \Vert \ddiv(p_\LS - \nabla u_\LS) \Vert_{L^2(T)}^2
        + h_T^2\,
        \Vert \ccurl(p_\LS - \nabla u_\LS) \Vert_{L^2(T)}^2\\
        &\phantom{{}\coloneqq{}}
        + h_T\,
        \sum_{E \in \E(T) \setminus \E(\partial\Omega)}
        \Vert [p_\LS - \nabla u_\LS]_E \cdot \nu_E \Vert_{L^2(E)}^2\\
        &\phantom{{}\coloneqq{}}
        + h_T\,
        \sum_{E \in \E(T)}
        \Vert [p_\LS - \nabla u_\LS]_E \cdot \tau_E \Vert_{L^2(E)}^2
    \end{split}
\end{equation}
The second term
$ \Vert \ccurl(p_\LS - \nabla u_\LS) \Vert_{L^2(T)}^2$
vanishes in the lowest-order case with
$(p_\LS, u_\LS) \in RT_0(\T) \times S_0^1(\T)$.
The discretisation of eigenvalue problems
in \cite{MR4410744},
based on first-order system least-squares formulations,
loses the built-in error estimation
property~\eqref{eq:equivalence}.
As a remedy,
an alternative error estimator similar to \(\eta_\SEP\)
enables a~posteriori error estimates in \cite[Sect.~5]{MR4410744}.

If the error in the flux variable is solely measured
in the \(L^2\) norm (and not the full \(H(\ddiv)\) norm)
the data oscillation term \eqref{eq:osc}
has to be included in the alternative error estimator
\cite{MR4011536,MR4271577}
\begin{equation}
    \label{eq:estimator_collective}
    \eta_\COL^2(\T)
    \coloneqq
    \sum_{T \in \T} \eta_\COL^2(\T, T)
    \quad\text{with}\quad
    \eta_\COL^2(\T, T)
    \coloneqq
    \eta_\SEP^2(\T, T)
    +
    h_T^2 \Vert (1-\Pi)f \Vert_{L^2(T)}^2.
\end{equation}
This provides a reliable and efficient error estimator
in the corresponding reduced norm
\cite[Eqn.~(5)]{MR4011536}
\begin{align*}
    LS(\Pi f; p_\LS, u_\LS)
    &\lesssim
    \Vert p - p_\LS \Vert_{L^2(\Omega)}^2
    +
    \Vert \nabla(u - u_\LS) \Vert_{L^2(\Omega)}^2
    \lesssim
    \eta_\COL^2(\T)
    \\
    &\lesssim
    LS(\Pi f; p_\LS, u_\LS)
    +
    \osc^2(f, \T).
\end{align*}
\begin{algorithm}[t]
    \caption{CALSFEM (collective marking adaptive LSFEM)}
    \label{alg:CALSFEM}
    \begin{algorithmic}
        \Require
        regular triangulation $\T_0$ and bulk parameter $0 < \theta \leq 1$.
        \For{$\ell = 0,1,2,\dots$}
            \State \textbf{Solve} LSFEM with respect to triangulation
            $\T_\ell$ for solution $(p_\ell, u_\ell)$.
            \State \textbf{Compute} $\eta_\COL(\T_\ell, T)$
            from~\eqref{eq:estimator_collective}
            for all \(T \in \T_\ell\).
            \State \textbf{Mark}  minimal subset
            $\mathcal{M}_\ell \subseteq \T_\ell$
            by the D\"{o}rfler criterion~\eqref{eq:bulk_criterion}
            for \(\eta \equiv \eta_\COL\).
            \State \textbf{Refine} $\T_\ell$ to $\T_{\ell + 1}$
            by NVB such that
            $\mathcal{M}_\ell \subseteq \T_\ell \setminus \T_{\ell+1}$.
        \EndFor
        \Ensure
        sequence of triangulations $\T_\ell$ with $(p_\ell, u_\ell)_{\ell}$
        and $\eta_\COL(\T_\ell)$ for $\ell \in \mathbb{N}_0$.
    \end{algorithmic}
\end{algorithm}

Replacing the built-in error estimator
\(\eta_\NAT\) in Algorithm~\ref{alg:NALSFEM} by
\(\eta_\COL\) leads to an alternative
adaptive Algorithm~\ref{alg:CALSFEM}
with collective marking (CALSFEM).
The estimator \(\eta_\COL\) guarantees
optimal convergence rates of CALSFEM with respect to
the reduced norm.
\begin{theorem}[{\cite[Sect.~2.5]{MR4011536}}]
    \label{thm:convergence_collective}
    For all \(0 < \theta \leq 1\),
    there exists \(0 < \varrho < 1\) such that
    the output \((p_\ell, u_\ell)_\ell\)
    of CALSFEM converges R-linearly,
    for all \(\ell, m \in \N_0\),
    \[
        \eta_\COL(\T_{\ell + m})
        \lesssim
        \rho^m\,
        \eta_\COL(\T_\ell).
    \]
    Moreover,
    there exists a maximal bulk parameter
    \(0 < \theta_0 < 1\)
    such that, for every
    \(0 < \theta \leq \theta_0\),
    the sequence \((p_\ell, u_\ell)_\ell\)
    converges with the optimal rate,
    i.e., for every \(0 < s < 1\),
    \[
        \sup_{\ell \in \N_0}
        (1 + \vert \T_\ell \vert - \vert \T_0 \vert)^s
        \eta_\COL(\T_\ell)
        \approx
        \sup_{N \in \N_0}
        (1 + N)^s
        \min_{\T \in \mathbb{T}(N)}
        \eta_\COL(\T).
    \]
\end{theorem}

Theorem~\ref{thm:convergence_collective}
generalises to higher-order discretisations
in three spatial dimensions
\cite[Sect.~2.8]{MR4271577}.

\subsection{Alternative adaptive least-squares FEM with separate marking}
\label{subsec:SALSFEM}

The optimal convergence rate of an adaptive algorithm
in the full \(H(\ddiv)\) norm for the flux variable
requires the reduction of the data approximation error
\begin{equation}
    \label{eq:data_error}
    \mu^2(\T)
    \coloneqq
    \sum_{T \in \T} \mu^2(T)
    \quad\text{with}\quad
    \mu^2(T)
    \coloneqq
    \Vert (1 - \Pi) f \Vert_{L^2(\Omega)}^2.
\end{equation}
The sum of this data error with the residual error estimator
from~\eqref{eq:estimator_separate} provides
a reliable and efficient error estimator
\cite[Thm.~3.1]{MR3296614}
\[
    \mu^2(\T) + \eta_\SEP^2(\T)
    \approx
    LS(f; p_\LS, u_\LS).
\]

Since the data error term \(\mu^2(T)\)
lacks any prefactor in terms of the mesh-size,
its strict reduction in the sense of the axioms of adaptivity
remains unclear.
In order to achieve optimal convergence rates,
Algorithm~{\ref{alg:SALSFEM}} (SALSFEM) employs
a separate marking strategy
\cite{MR2772091,MR3719030}.
\begin{algorithm}[t]
    \caption{SALSFEM (separate marking adaptive LSFEM)}
    \label{alg:SALSFEM}
    \begin{algorithmic}
        \Require
        regular triangulation $\T_0$,
        bulk parameter $0 < \theta \leq 1$,
        reduction parameter $0 < \rho < 1$,
        and separation parameter $0 < \kappa$.
        \For{$\ell = 0,1,2,\dots$}
            \State \textbf{Solve} LSFEM with respect to triangulation
            $\T_\ell$ for solution $(p_\ell, u_\ell)$.
            \State \textbf{Compute} $\eta_\SEP(\T_\ell, T)$
            from~\eqref{eq:estimator_separate}
            for all \(T \in \T_\ell\).
            \If{\textbf{Case A} $\mu^2(\T) \leq \kappa\,\eta_\SEP^2(\T)$}
                \State \textbf{Mark} minimal subset
                $\mathcal{M}_\ell \subseteq \T_\ell$
                by D\"{o}rfler criterion~\eqref{eq:bulk_criterion}
                for \(\eta \equiv \eta_\SEP\).
                \State \textbf{Refine} $\T_\ell$ to $\T_{\ell + 1}$
                by NVB such that
                $\mathcal{M}_\ell \subseteq \T_\ell \setminus \T_{\ell+1}$.
                \Else{ (\textbf{Case B}
                    $\kappa\,\eta_\SEP^2(\T) < \mu^2(\T)$)}
                \State \textbf{Compute} a refinement
                $\T_{\ell+1}$ of $\T_\ell$ of
                (almost) minimal
                cardinality
                with $\mu(\T_{\ell+1}) \leq \rho\, \mu(\T_\ell)$.
            \EndIf
        \EndFor
        \Ensure
        sequence of triangulations $\T_\ell$
        with $(p_\ell,u_\ell)_{\ell}$, $\eta_\SEP(\T_\ell)$,
        and $\mu(\T_\ell)$ for $\ell \in \mathbb{N}_0$.
    \end{algorithmic}
\end{algorithm}
If the residual error estimator \(\eta_\SEP^2(\T)\) dominates
the data error \(\mu^2(\T)\),
the former is refined by the standard D\"orfler marking and NVB.
Otherwise, the latter is reduced by
a suitable data approximation algorithm.
The data approximation in Case~B of SALSFEM employs the
approximation algorithm (AA) from~\cite{MR3377430}.
It consists of a slight modification of
the Thresholding Second Algorithm (TSA) from~\cite{MR2050076}
and utilises binary bins
to guarantee linear computational complexity
\cite[Rem.~5.3]{MR2050076}.
The algorithm considers
the refinement indicator \(\widetilde\mu(T_j)\)
for the two children \(T_1\) and \(T_2\)
of a bisected parent triangle \(T\)
defined, for \(j = 1,2\), by
\begin{equation}
    \label{eq:modified_data_error_indicator}
    \widetilde\mu(T_j)
    \coloneqq
    (\mu(T_1) + \mu(T_2))\, \widetilde\mu(T) / (\mu(T) + \widetilde\mu(T))
\end{equation}
with \(\widetilde\mu(T) \coloneqq \mu(T)\) for all initial triangles
\(T \in \T_0\).
The TSA is followed by a completion step
in order to ensure the output triangulation to be shape-regular.
\begin{algorithm}[t]
    \caption{Approximation Algorithm (AA)}
    \label{alg:approx}
    \begin{algorithmic}
        \Require
        initial regular triangulation \(\T_0\),
        error tolerance \(\textrm{Tol} > 0\)
        \State \textbf{Compute} \(\mu(T) = \widetilde\mu(T)\)
        for all \(T \in \T_0\) and set \(\widehat\T
        \coloneqq \T_0\).
        \While{\(\mu(\widehat\T) > \textrm{Tol}\)}
            \State {%
                \textbf{Select}
                the minimal \(k \in \mathbb{Z}\) such that
                \(
                    \widetilde\mu(T)
                    <
                    2^{k+1}
                \)
                for all \(T \in \widehat\T\).
                \State \textbf{Mark} the set
                \(
                    \mathcal{M}
                    \coloneqq
                    \{
                        T \in \widehat\T \;:\;
                        2^k \leq
                        \widetilde\mu(T)
                        < 2^{k+1}
                    \}
                \).
            }
            \State \textbf{Bisect} all triangles in \(\mathcal{M}\)
            to obtain a new \(\widehat\T\).
            \State {%
                \textbf{Compute} \(\mu(T)\) and
                \(\widetilde\mu(T)\)
                from~\eqref{eq:modified_data_error_indicator}
                for all newly created \(T \in \widehat\T\)
            }
        \EndWhile
        \State Apply \textbf{completion} on \(\widehat\T\)
        to obtain a regular refinement
        \(\T_\textup{Tol}\) of \(\T_0\).
        \Ensure \(\T_\textup{Tol}\)
    \end{algorithmic}
\end{algorithm}
The resulting Algorithm~\ref{alg:approx} (AA)
is instance optimal \cite{MR2050076,MR2050077}.
Then Algorithm~\ref{alg:SALSFEM} SALSFEM
converges with the optimal rate.
\begin{theorem}[{\cite[Thm.~6.1]{MR3296614}}]
    \label{thm:convergence_separate}
    For all \(0 < \theta \leq 1\),
    \(0 < \kappa\), and \(0 < \rho < 1\),
    there exists \(0 < \varrho < 1\) such that
    the output \((p_\ell, u_\ell)_\ell\)
    of SALSFEM converges R-linearly,
    for all \(\ell, m \in \N_0\),
    \[
        \mu^2(\T_{\ell + m})
        +
        \eta_\SEP^2(\T_{\ell + m})
        \lesssim
        \varrho^m\,
        \big(
            \mu^2(\T_\ell)
            +
            \eta_\SEP^2(\T_\ell)
        \big).
    \]
    Moreover,
    there exists a maximal bulk parameter
    \(0 < \theta_0 < 1\)
    and a maximal separation parameter
    \(0 < \kappa_0\)
    such that for all
    \(0 < \theta \leq \theta_0\),
    \(0 < \kappa \leq \kappa_0\),
    and \(0 < \rho < 1\),
    the sequence \((p_\ell, u_\ell)_\ell\)
    converges with the optimal rate,
    i.e., for all \(0 < s < 1\),
    \[
        \sup_{\ell \in \N_0}
        (1 + \vert \T_\ell \vert - \vert \T_0 \vert)^s
        (\mu(\T_\ell) + \eta_\SEP(\T_\ell))
        \approx
        \sup_{N \in \N_0}
        (1+N)^s
        \min_{\T \in \mathbb{T}(N)}
        (\mu(\T) + \eta_\SEP(\T)).
    \]
\end{theorem}

For the generalisation to higher-order polynomial
degrees and inhomogeneous mixed boundary conditions
in three spatial dimensions,
the reader is referred to \cite{MR4557622}.

\section{Numerical experiments}
\label{sec:experiments}
This section presents and compares the numerical results of NALSFEM, CALSFEM,
and SALSFEM for three benchmark examples of the Poisson model problem
and one of an elliptic problem with piecewise constant scalar
diffusion constant.
A primary focus consists of investigating
the data approximation in
Subsection~\ref{subsec:LshapeMicrostructure}
below.
The lowest-order discretisation
prevents any additional quadrature error
for this benchmark problem.

\subsection{Implementation and time measurement}
\label{subsec:implementation}
The empirical investigation was carried out using
the author's \textsc{Matlab} software package octAFEM
\cite{octAFEM}.
All experiments in this paper are reproducible
with the compute capsule on the Code Ocean platform.
The octAFEM package bases
on the in-house \textsc{Matlab} software
package~\cite{AFEMpackage}.
It was developed and tested under
\textsc{Matlab} version 9.14.0.2206163 (R2023a),
but should be executable in older versions as well.
Moreover, the code is completely compatible
with the open-source software Octave
(tested with version 8.1.0).
The realisation differs from~\cite{BringmannDissertation}
because the object-oriented implementation therein
employs a \texttt{Simplex} class
for the representation of every simplex separately
resulting in a huge computational overhead.
Instead, the \texttt{ApproxTriangulation}
in the implementation at hand
includes an array of indices
containing the complete history of simplices.
The data approximation Algorithm~\ref{alg:approx} (AA)
ensures linear complexity using binary bins
as described in \cite[Rem.~5.3]{MR2050076}.
It is incorporated into the \texttt{Triangulation} class
from~\cite{BringmannDissertation}.
This allows the separate marking strategy
in one triangulation object
containing the complete refinement history
of all simplices and thereby
avoiding the computation of the overlay
of \(\T_\textup{Tol}\) and \(\T_\ell\).
The data error \(\mu(T)\)
of each simplex \(T \in \T\) is
stored in the array of \texttt{ApproxTriangulation}
as well.
It is computed when creating the simplex \(T\).
This causes some general overhead to the refinement process
but may lead to some reduction of the runtime of AA
because it can reuse information already created during
a previous step of the NVB in a Case~A of the separate marking algorithm.

The transformation formula allows to reduce
the integral over any triangle to the
reference triangle
\(T_\textup{ref} \coloneqq \conv\{0, (1,0), (0,1)\}\).
The transformation
\(\Phi: [0,1]^2 \to T_\textup{ref}\),
\(y \mapsto (y_1, (1-y_1)y_2)^\top\)
from the unit square to the reference triangle
shows
\[
    \int_{T_\textup{ref}}
    f \,\mathrm{d}x
    =
    \int_0^1 \int_0^1
    (1- y_1) \, f(y_1, (1 - y_1)y_2)
    \,\mathrm{d}y_2 \,\mathrm{d}y_1.
\]
The first integral with respect to \(y_2\)
is approximated by the Gauss--Legendre quadrature.
The second integral with respect to \(y_1\)
employs the Gauss--Jacobi quadrature on
the interval \([0, 1]\) with
weight function \(w(\xi) = (1 - \xi)\).
Both one-dimensional quadrature nodes and weights
are computed using the Golub--Welsch algorithm
\cite{MR0245201}
with recursion coefficients from
\cite{MR3509205}.
The resulting conical product rules
with \(k^2\) function evaluations, \(k \in \N\),
are exact for the integration of
polynomials up to partial degree
\(2k - 1\).
For the evaluation of bilinear forms
or integration of polynomial input data,
the number of quadrature points is chosen
such that the quadrature is exact.

The experiments investigating the performance
of the algorithms in terms of the runtime are
carried out
with \textsc{Matlab} version 9.9.0.1467703 (R2020b)
on a compute server
using 16 out of 128 Intel(R) Xeon(R) E7-8867 CPUs of 2.50GHz
and 2\,TiB RAM.
The code employs parallel computing for local quantities
such as local stiffness matrices and
the integration of the right-hand side.
Since the \textsc{Matlab} command \texttt{cputime}
adds up the time for all parallel threads,
the documentation recommends the measurement
of real time.
Additionally this exemplifies the practical performance
as experienced by the user.
To this end, time is measured on carefully selected parts of the program
to distinguish the performance for the solution, estimation,
and refinement.
This allows to neglect possible overhead due to printing information
to the command line or saving the results to disk.
The time is measured in ten separate runs
and averaged for improved realiability.
The graphs below also indicate the maximal and the minimal measured
time by vertical error bars to visualise possible inaccuracies
of the measurement.
The only significant differences occur for the very first iterations
of each adaptive computation.

\subsection{L-shaped domain}
\label{subsec:Lshape}
\begin{figure}[p]
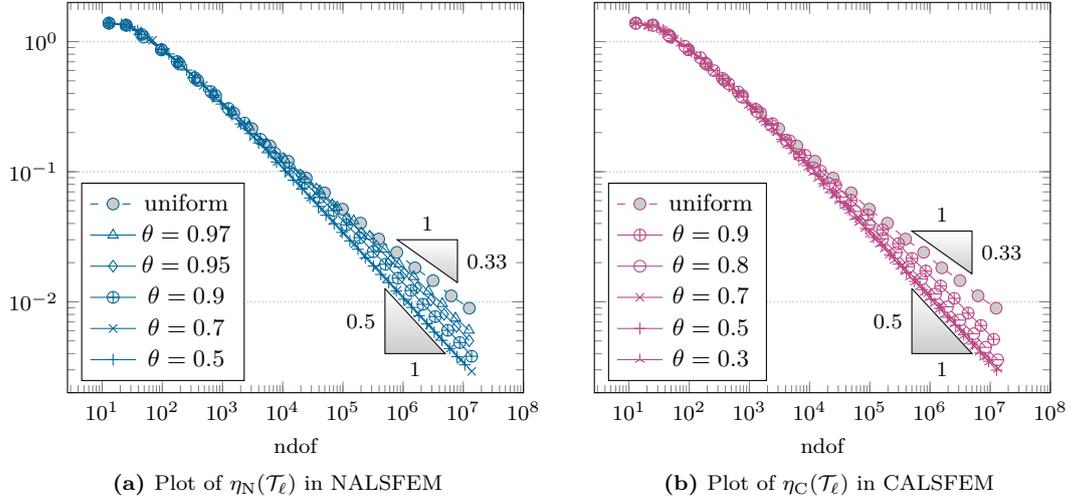

    \subfloat[Plot of \(\eta_\NAT(\T_\ell)\) in NALSFEM]{%
        \label{subfig:Lshape_convergence_natural}
        \input{figures/Poisson_Lshape_LS_adap_nat_theta.tex}
    }
    \subfloat[Plot of \(\eta_\COL(\T_\ell)\) in CALSFEM]{%
        \label{subfig:Lshape_convergence_collective}
        \input{figures/Poisson_Lshape_LS_adap_col_theta.tex}
    }
    \caption{%
        Comparison of various choices
        for the bulk parameter \(0 < \theta \leq 1\)
        in the adaptive mesh-refinement strategies
        (uniform refinement for \(\theta = 1\))
        for the benchmark problem on the L-shaped domain
        from Subsection~\ref{subsec:Lshape}.
    }
    \label{fig:Lshape_convergence}
\end{figure}

The Poisson model problem on the L-shaped domain
$\Omega = (-1,1)^2 \setminus [0,1)^2$
with constant right-hand side $f\equiv1$
is a standard benchmark for adaptive mesh-refinement.
The reentrant corner leads to reduced elliptic regularity
of the unknown exact solution
\(u \in H^{1+s-\varepsilon}(\Omega)\) with \(s = 2/3\)
for all \(\varepsilon > 0\).
This is why uniform refinement exhibits a suboptimal convergence rate
of \(1/3\) with respect to the number of degrees of freedom (ndof)
for the natural estimator \(\eta_\NAT\)
and the alternative estimator \(\eta_\COL\)
in Figure~\ref{fig:Lshape_convergence}.

The convergence result in Theorem~\ref{thm:convergence_collective}
asserts optimal rates for CALSFEM for sufficiently small
bulk parameters \(\theta < \theta_0\).
The upper bound
\(\theta_0 = (1 + C_\textup{stab}^2 C_\textup{drel})^{-1}\)
from~\cite[Prop.~4.2~(ii)]{MR3170325} includes the generic constants
of the stability and discrete reliability axiom.
These constants are bounded in \cite[Sect.~6]{MR3824773}
in the case of the Courant FEM
on a mesh with right-iscosceles triangles
for the Poisson model problem by
\begin{equation}
    \label{eq:constants_bounds}
    C_\textup{stab}^2 \leq 40.36
    \quad\text{and}\quad
    C_\textup{drel} \leq 9\,201.
\end{equation}
This leads to the small theoretical lower bound of
\(\theta_0 \geq 2.6 \times 10^{-6}\).
Nevertheless, Figure~\ref{subfig:Lshape_convergence_collective}
shows the optimal convergence rate already for moderate
bulk parameters \(\theta \leq 0.8\) in practice.

\begin{figure}[p]
    \centering
    \subfloat[NALSFEM (3\,946 triangles)]{
        \begin{tikzpicture}
    \begin{axis}[%
        axis equal image,%
        width=5.4cm,%
        xmin=-1.15, xmax=1.15,%
        ymin=-1.15, ymax=1.15,%
        font=\footnotesize%
    ]

        \addplot graphics [xmin=-1, xmax=1, ymin=-1, ymax=1]
        {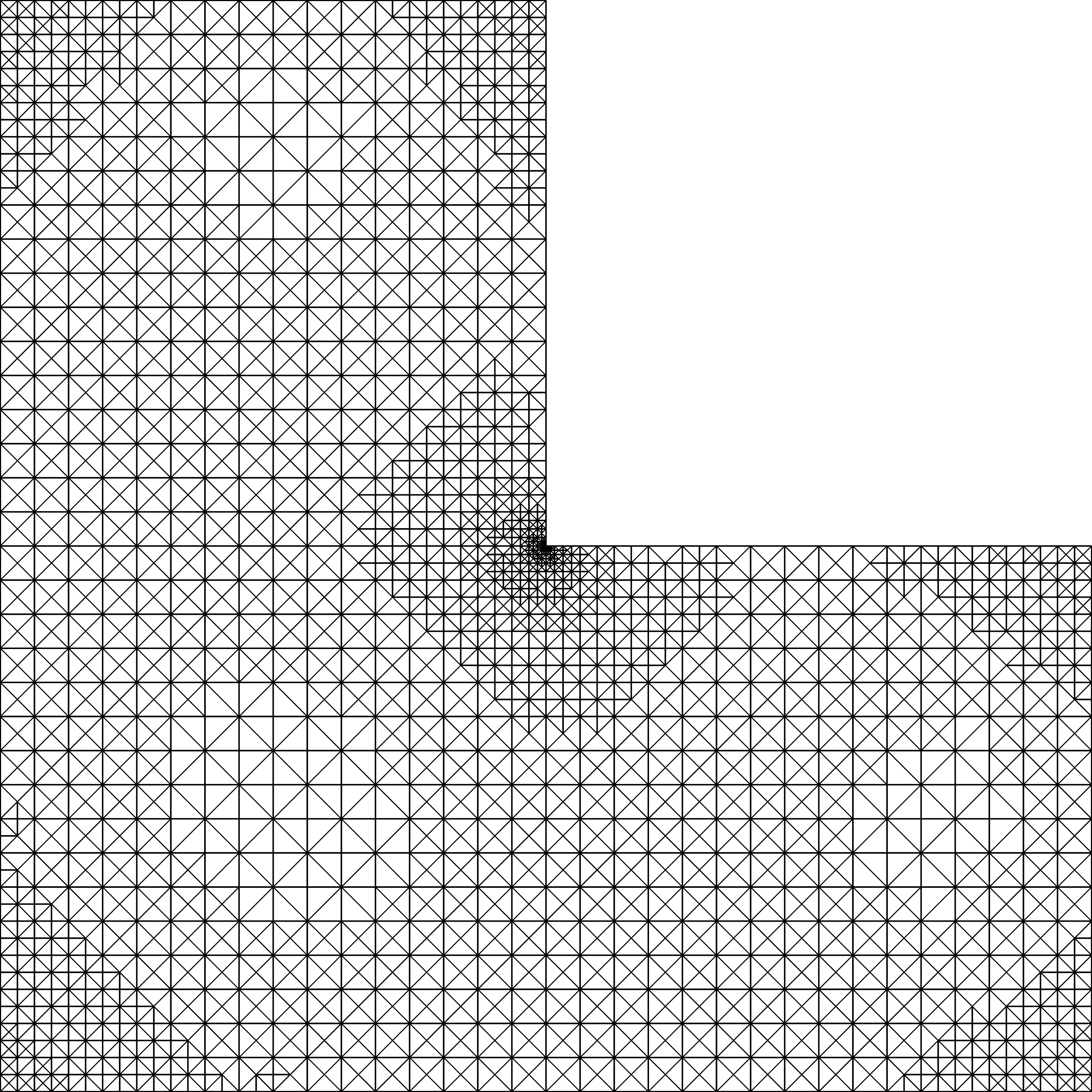};
    \end{axis}
\end{tikzpicture}
    }
    \subfloat[CALSFEM (3\,314 triangles)]{
        \begin{tikzpicture}
    \begin{axis}[%
        axis equal image,%
        width=5.4cm,%
        xmin=-1.15, xmax=1.15,%
        ymin=-1.15, ymax=1.15,%
        font=\footnotesize%
    ]

        \addplot graphics [xmin=-1, xmax=1, ymin=-1, ymax=1]
        {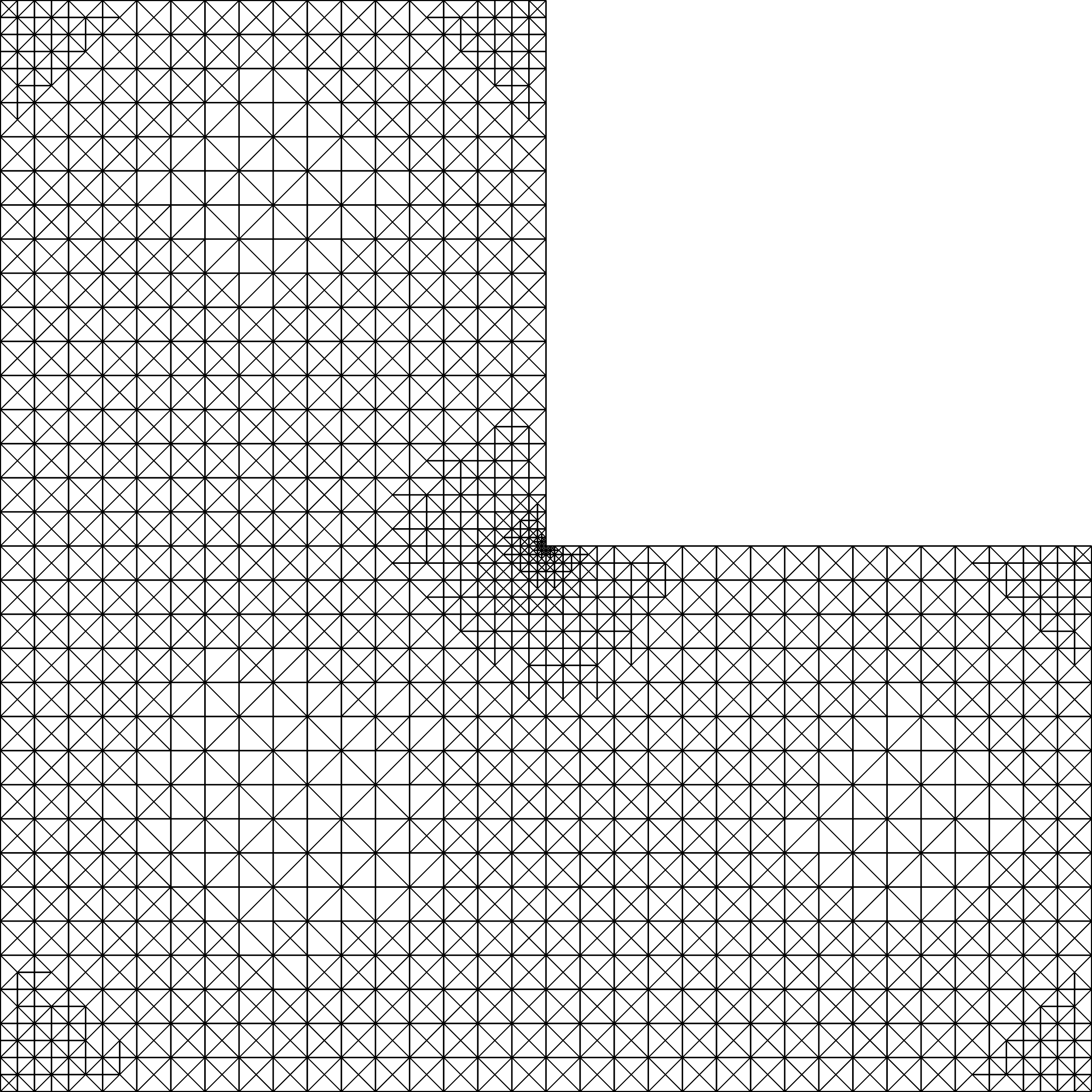};
    \end{axis}
\end{tikzpicture}
    }
    \hspace{4.5em}~

    \subfloat[NALSFEM (5\,314\,052 triangles)]{
        \begin{tikzpicture}
    \begin{axis}[%
        axis equal image,%
        width=5.4cm,%
        xmin=-1.15, xmax=1.15,%
        ymin=-1.15, ymax=1.15,%
        font=\footnotesize%
    ]

        \addplot graphics [xmin=-1, xmax=1, ymin=-1, ymax=1]
        {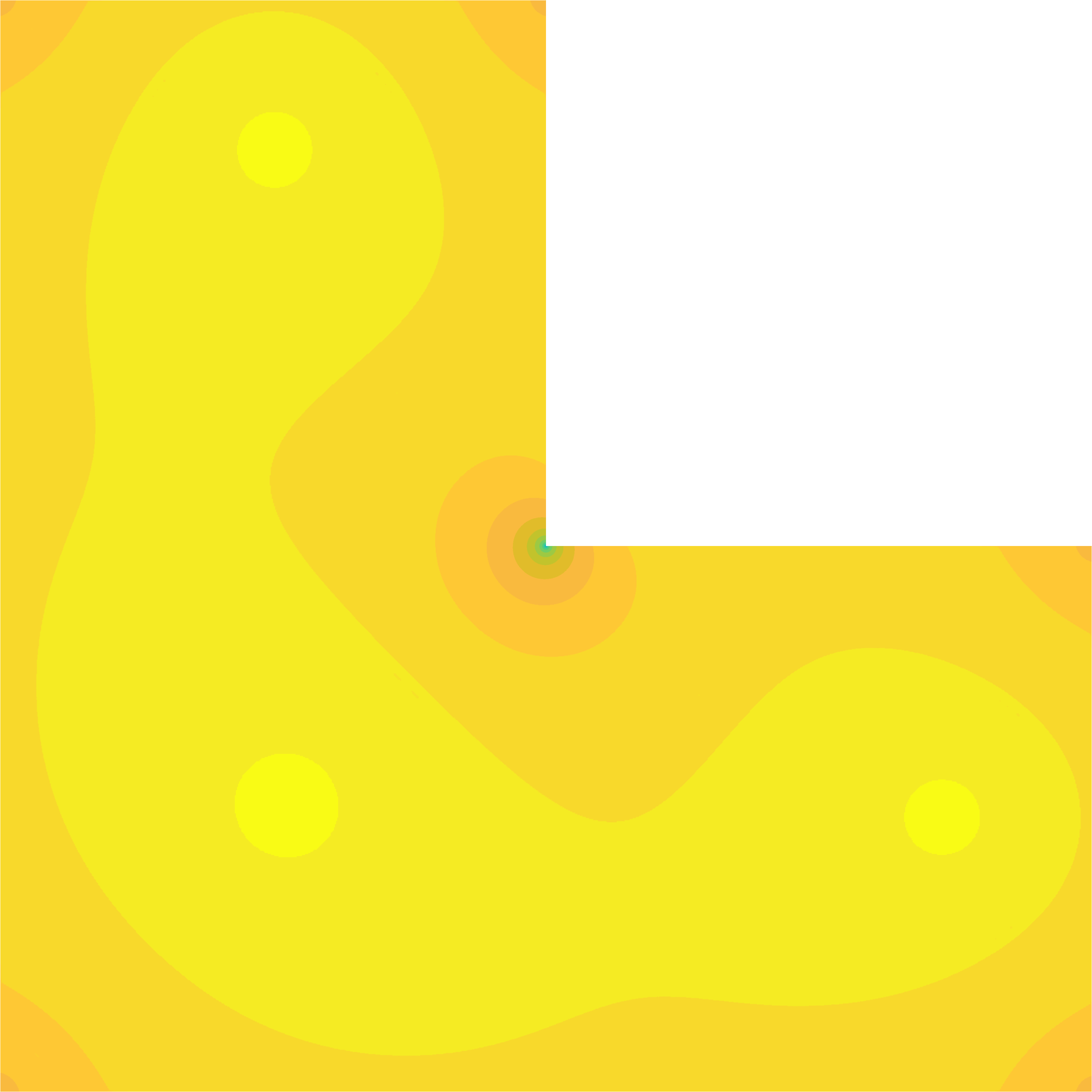};

    \end{axis}
\end{tikzpicture}
    }
    \subfloat[CALSFEM (6\,405\,512 triangles)]{
        \begin{tikzpicture}
    \begin{axis}[%
        axis equal image,%
        width=5.4cm,%
        xmin=-1.15, xmax=1.15,%
        ymin=-1.15, ymax=1.15,%
        font=\footnotesize,%
        point meta min=-6.32162991,%
        point meta max=-2.85978496,%
        colorbar,%
        colorbar style={%
            title={\(h_\ell\)},%
            font=\footnotesize,%
            width=2.5mm,%
            title style={yshift=-2mm},%
            yticklabel={$10^{\pgfmathprintnumber{\tick}}$},%
        },%
    ]

        \addplot graphics [xmin=-1, xmax=1, ymin=-1, ymax=1]
        {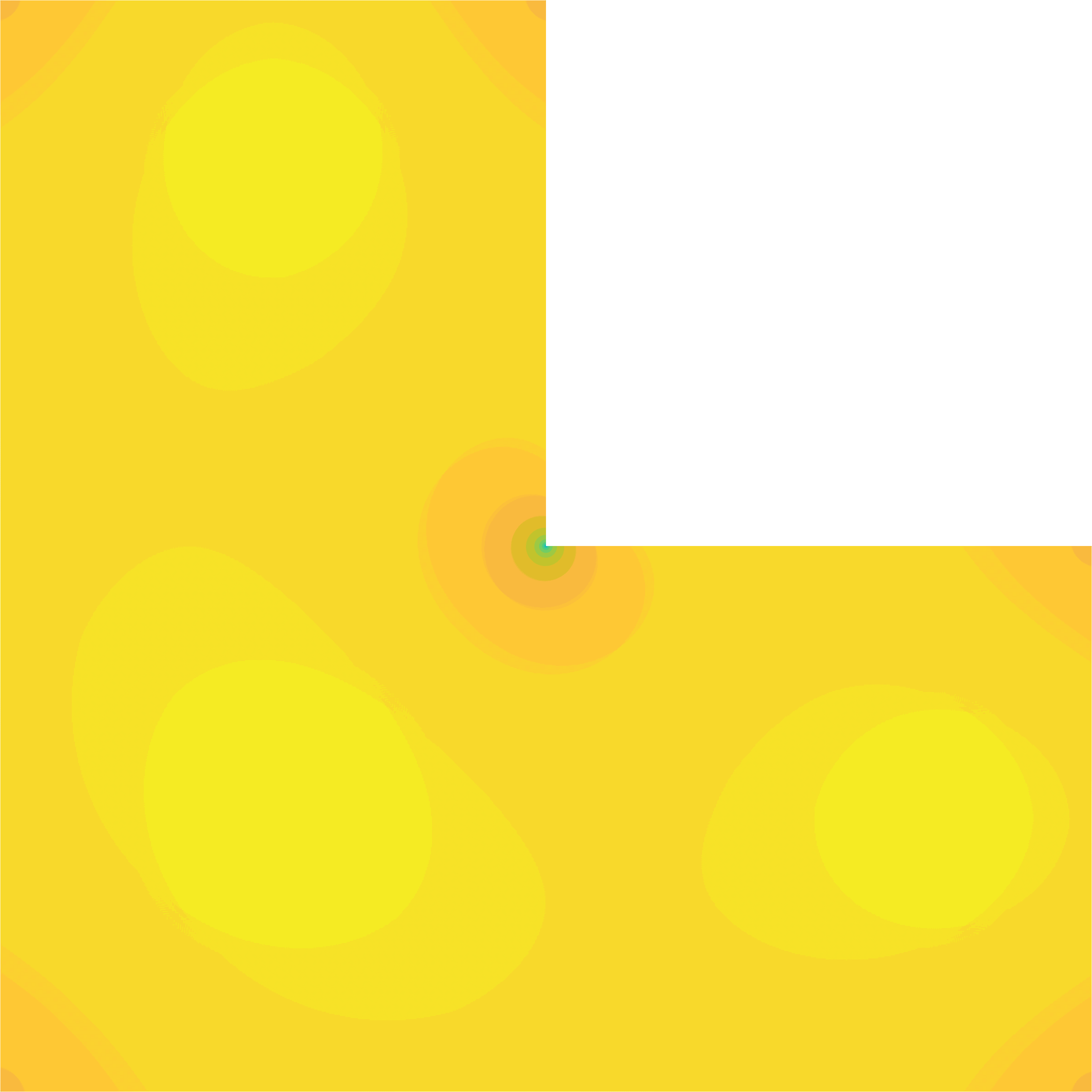};

    \end{axis}
\end{tikzpicture}
    }
    \caption{%
        Adaptively refined meshes from two different adaptive algorithms
        with bulk parameter \(\theta = 0.9\)
        for the benchmark problem from Subsection~\ref{subsec:Lshape}.
        The Subfigures~(c) and~(d) display
        the mesh-size \(h_\ell \vert_T \equiv \vert T \vert^{1/2}\)
        for the triangles \(T \in \T_\ell\) in a very fine triangulation.
        The same color scale enables the comparison
        of the mesh-size in the two algorithms.
    }
    \label{fig:Lshape_meshes}
\end{figure}

The algorithm NALSFEM converges with the optimal rate
for even larger bulk parameters \(\theta \leq 0.9\).
The alternative estimator \(\eta_\COL\)
focuses on the constitutive residual
while the natural estimator \(\eta_\NAT\)
includes the equilibrium residual as well.
This may explain the better performance of
the natural refinement strategy for large bulk parameters.
This difference is small and
the coarse adaptively generated meshes
look essentially identical
for both refinement strategies
as displayed in Figure~\ref{fig:Lshape_meshes}.
A closer investigation of the fine triangulations
with one million triangles and more exhibit
an increased adaptive refinement towards the reentrant corner
while at the same time allowing coarser triangles
in the remaining parts of the domain
for the NALSFEM compared to the CALSFEM

Since the data \(f\) is resolved exactly
on every triangulation \(\T \in \mathbb{T}\),
the data error and oscillation terms vanish
\(\mu(\T) = \osc(f, \T) = 0\).
Hence, the case~B in the separate marking does never hold and
the SALSFEM algorithm provides exactly the same results
as CALSFEM.

\subsection{L-shaped domain with microstructure}
\label{subsec:LshapeMicrostructure}

\begin{figure}[p]
    \centering
    \subfloat[Definition of right-hand side \(f_\epsilon\)]{
        \begin{tikzpicture}
    \begin{axis}[%
        axis equal image,%
        width = 0.36\textwidth,%
        ymin  = -1.3,%
        xmin  = -1.3,%
        font  = \footnotesize%
    ]

        \addplot[%
            mark      = none,%
            thick,%
            line cap  = round,%
            line join = round,%
        ] coordinates {
            (0, 0)
            (1, 0)
            (1, -1)
            (-1, -1)
            (-1, 1)
            (0, 1)
            (0, 0)
        };

        \fill[HUblue!66!white]
            (axis cs: -0.6, 0.6) -- (axis cs: -0.6, 0.4) --
            (axis cs: -0.4, 0.4) -- (axis cs: -0.4, 0.6) -- cycle;

        \draw[densely dotted,HUblue]
            (axis cs: -1, 0.6) -- (axis cs: -0.4, 0.6)
            (axis cs: -1, 0.4) -- (axis cs: -0.4, 0.4)
            (axis cs: -0.6, -1) -- (axis cs: -0.6, 0.6)
            (axis cs: -0.4, -1) -- (axis cs: -0.4, 0.6);

        \draw [decorate, decoration={brace, amplitude=2pt}, semithick, HUblue]
            (axis cs: -1, 0.4) -- (axis cs: -1, 0.6);
        \draw [decorate, decoration={brace, amplitude=2pt}, semithick, HUblue]
            (axis cs: -0.4, -1) -- (axis cs: -0.6, -1);

        \path (axis cs: -1, 0.5) node[left] {\color{HUblue}\(2\epsilon\)};
        \path (axis cs: -0.5, -1) node[below] {\color{HUblue}\(2\epsilon\)};

        \draw
            (axis cs: 0.3, 0.7) -- (axis cs: 0.3, 0.85) --
            (axis cs: 0.45, 0.85) -- (axis cs: 0.45, 0.7)
            node[midway,right] {\(f_\epsilon \equiv 0\)} -- cycle;
        \draw[fill=HUblue!66!white]
            (axis cs: 0.3, 0.35) -- (axis cs: 0.3, 0.5) --
            (axis cs: 0.45, 0.5) -- (axis cs: 0.45, 0.35)
            node[midway,right] {\(f_\epsilon \equiv 1\)} -- cycle;
    \end{axis}
\end{tikzpicture}
    }
    \hfil
    \subfloat[Discrete solution \(u_\ell\) for \(\epsilon = 2^{-5}\)]{
        \begin{tikzpicture}
    \pgfplotsset{/pgf/number format/fixed}

    \begin{axis}[%
        width=0.46\textwidth,%
        xmin=-1.1, xmax=1.1,%
        ymin=-1.1, ymax=1.1,%
        zmin=-0.0005, zmax=0.0025,%
        font=\footnotesize%
    ]
        \addplot3 graphics [%
            points={%
                (-1,1,0) => (0,360-269)
                (-1,-1,0) => (156,0)
                (-0.5,0.5,0.00205) => (90,360-111)
                (1,0,0) => (282,360-245)
            }
            ]
            {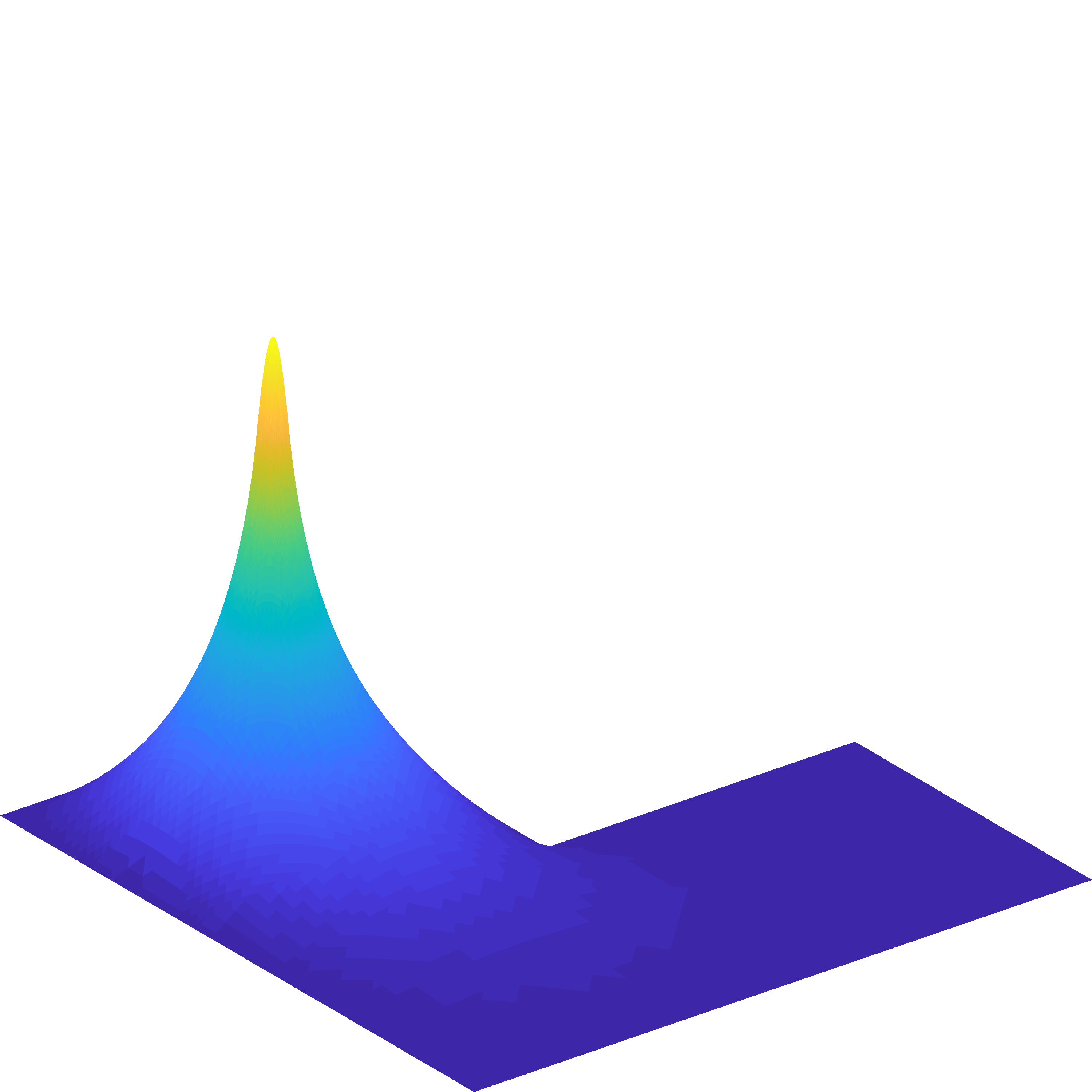};
    \end{axis}
\end{tikzpicture}
    }
    \caption{%
        Illustration of the microstructure in the benchmark problem
        from Subsection~\ref{subsec:LshapeMicrostructure}.
    }
    \label{fig:microstructure}
\end{figure}

The second benchmark considers the L-shaped domain
$\Omega = (-1,1)^2 \setminus [0,1)^2$
from the Section~\ref{subsec:Lshape}
with a right-hand side \(f_\epsilon \in L^2(\Omega)\)
for some parameter \(0 < \epsilon < 1/2\),
given in \cite[Sect.~3.4]{MR3377430} by
\begin{align*}
    f_{\epsilon}(x)
    \coloneqq
    \begin{cases}
        1, &\text{if }
        \vert x_1 + \frac12 \vert \leq \epsilon
        \text{ and } \vert x_2 - \frac12 \vert \leq \epsilon,\\
        0, &\text{otherwise}.
    \end{cases}
\end{align*}
Figure~\ref{fig:microstructure} illustrates the definition
of \(f_\epsilon\) and shows an example solution
for \(\epsilon = 2^{-5}\).

Due to the small support of the right-hand side \(f_\epsilon\),
the quadrature described in Section~\ref{subsec:implementation}
may be inaccurate, in particular for coarse triangulations.
This is why the integration of the right-hand side
for this benchmark is computed directly as the area
of the convex intersection polygon of the support
\(
    \operatorname{supp}(f_\epsilon)
    =
    (-1/2 - \epsilon, -1/2 + \epsilon)
    \times (1/2 - \epsilon, 1/2 + \epsilon)
\)
and any triangle \(T \in \T\).
First, the vertices of the intersection polygon
\(x_1, \dots, x_J\)
are determined by the Sutherland-Hodgman algorithm
\cite{SH74}.
Second, the area of the intersection polygon
is computed by the formula
\[
    \vert \operatorname{supp}(f_\epsilon) \cap T \vert
    =
    \frac12
    \Big\vert
    \sum_{j=1}^J
    (x_{j,1} x_{j+1,2} - x_{j,2} x_{j+1,1})
    \Big\vert
    \quad\text{with } x_{J+1} \equiv x_1.
\]
This procedure allows for the exact computation
of the piecewise constant approximation \(\Pi f_\epsilon\)
and the data error \(\mu(\T)\).

If the microstructure can be resolved exactly
for \(\epsilon = 2^{-m}\) with \(m \in \N\),
all three algorithms reach the point of exact data resolution
and converge with the best possible rate from then on
as displayed in Figure~\ref{subfig:microstructure_convergence_exact}
for \(\epsilon = 2^{-5}\).
Otherwise the data approximation plays a crucial role
throughout the whole computation as for \(\epsilon = 3^{-3}\)
in Figure~\ref{subfig:microstructure_convergence}.
For the algorithms NALSFEM and SALSFEM
the least-squares functional converges
with the optimal rate of \(0.5\).
The indication of the cases at the top of the plot
shows that the data approximation dominates
in the first eight iterations of SALSFEM.
The convergence behaviour of NALSFEM
turns out to be very close to the separate marking algorithm
but with significantly more intermediate solution steps.
The alternative estimator \(\eta_\textup{C}\) in CALSFEM
converges with the optimal rate as well,
as asserted by Theorem~\ref{thm:convergence_collective}.
However, it does not allow to control the data approximation
as part of the divergence contribution to the error of the flux variable.
This results in a suboptimal rate of \(0.25\)
at the beginning of the computation when the data oscillation
is presumably large enough.
Once the dominance of the data approximation ends at about
\(5 \times 10^4\) degrees of freedom,
the CALSFEM algorithm is not able to considerably
reduce the least-squares estimator any more.
\begin{figure}[p]
    \centering
    \subfloat[\(\epsilon = 2^{-5}\)]{%
        \label{subfig:microstructure_convergence_exact}
        \input{figures/Poisson_LshapeMicrostructure_LS_adap30_strategies.tex}
    }
    \subfloat[\(\epsilon = 3^{-3}\)]{%
        \label{subfig:microstructure_convergence}
        \input{figures/Poisson_LshapeMicrostructure_inexact_LS_adap30_strategies.tex}
    }

    \subfloat[Legend of
        Figures~\ref{subfig:microstructure_convergence_exact}%
        --\ref{subfig:microstructure_convergence}
        and~\ref{fig:performance}]{
        \label{subfig:legend_comparison}
        \begin{tikzpicture}
    \pgfplotsset{%
        eta/.style={%
            mark=*,%
            every mark/.append style={%
                solid,%
                scale=1.1,%
            }%
        },%
        res/.style={%
            mark=square*,
            every mark/.append style={%
                solid,%
                scale=1.1,%
            }%
        },%
        resL2/.style={%
            mark=halfsquare left*,%
            every mark/.append style={%
                solid,%
                scale=1.1,%
            }%
        },%
        resDiv/.style={%
            mark=halfsquare right*,%
            every mark/.append style={%
                solid,%
                scale=1.1,%
            }%
        },%
        mu/.style={%
            mark=diamond*,%
            every mark/.append style={%
                solid,%
                scale=1.1,%
            }%
        },%
        cases/.style={%
            mark=*,%
            only marks,%
            HUgreen,%
            every mark/.append style={%
                solid,%
                scale=0.5,%
            }%
        },%
        casesLetters/.style={%
            mark=,%
            only marks,%
            HUgreen,%
            nodes near coords,%
            point meta=explicit symbolic,%
            every node near coord/.style={%
                font=\footnotesize%
            }%
        },%
        uniform/.style={%
            dashed,%
            black,%
            every mark/.append style={fill=black!20!white}%
        },%
        natural/.style={%
            dashed,%
            HUblue,%
            every mark/.append style={%
                HUblue,%
                fill=white,%
            }%
        },%
        collective/.style={%
            dotted,%
            HUred,%
            every mark/.append style={%
                HUred,%
                fill=HUred!33!white,%
            }%
        },%
        separate/.style={%
            solid,%
            HUgreen,%
            every mark/.append style={
                HUgreen,%
                fill=HUgreen!66!white,%
            }%
        },%
    }

    \matrix [
        matrix of nodes,
        anchor = south,
        font = \scriptsize,
        column 1/.style={anchor=base east},
    ] at (0,0) {
                      & NALSFEM & CALSFEM & SALSFEM \\
        \hline \\
        \(LS(f;p_\ell,u_\ell)^{1/2}\)
        & \ref{leg:nat:eta}
        & \ref{leg:col:res}
        & \ref{leg:sep:res} \\
        Alternative estimator
        & \(\eta_\NAT\)
        & \(\eta_\COL\)
        & \((\eta_\SEP^2 + \mu^2)^{1/2}\) \\
        & \ref{leg:nat:eta}
        & \ref{leg:col:eta}
        & \ref{leg:sep:eta} \\
        cases A or B & & & \ref{leg:sep:cases} \\
    };
\end{tikzpicture}
    }
    \caption{%
        Comparison of the three adaptive refinement strategies
        with parameters \(\theta = 0.3\), \(\kappa = 1\), and \(\rho = 0.8\)
        for the solution of the
        benchmark problem
        from Subsection~\ref{subsec:LshapeMicrostructure}.
    }
\end{figure}

\begin{figure}
    \centering
    \subfloat[NALSFEM (\(123\,909\) triangles)]{%
        \label{subfig:microstructure_mesh_natural}
        \begin{tikzpicture}
    \begin{axis}[%
        axis equal image,%
        width=5.4cm,%
        xmin=-1.15, xmax=1.15,%
        ymin=-1.15, ymax=1.15,%
        font=\footnotesize,%
    ]

        \addplot graphics [xmin=-1, xmax=1, ymin=-1, ymax=1]
        {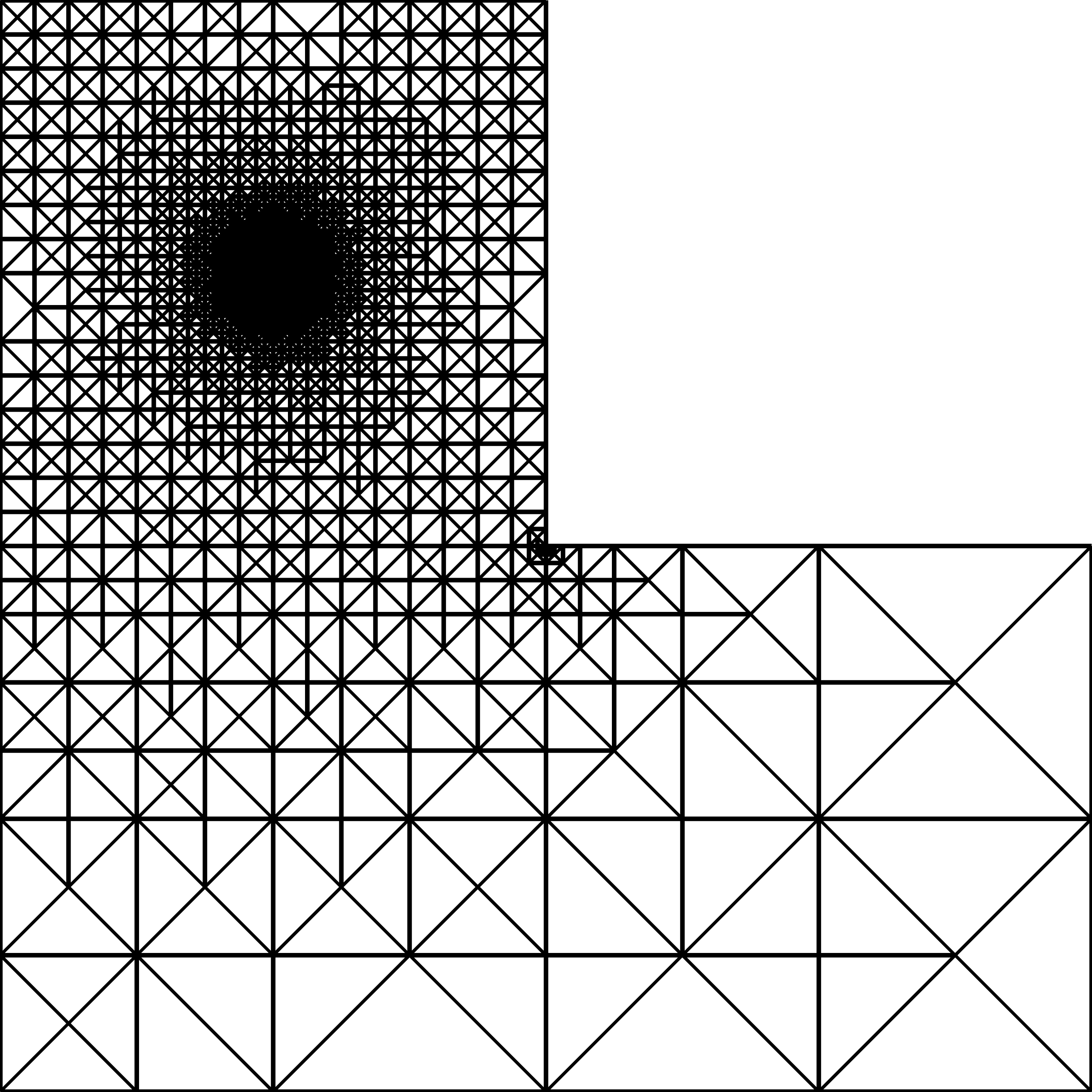};

    \end{axis}
\end{tikzpicture}
    }
    \hfil
    \subfloat[NALSFEM (\(5\,333\,087\) triangles)]{%
        \label{subfig:microstructure_meshsize_natural}
        \begin{tikzpicture}
    \begin{axis}[%
        axis equal image,%
        width=5.4cm,%
        xmin=-1.15, xmax=1.15,%
        ymin=-1.15, ymax=1.15,%
        font=\footnotesize,%
        point meta min=-6.62265990,%
        point meta max=-1.05360498,%
        colorbar,%
        colorbar style={%
            title={\(h_\ell\)},%
            font=\footnotesize,%
            width=2.5mm,%
            title style={yshift=-2mm},%
            yticklabel={$10^{\pgfmathprintnumber{\tick}}$},%
            ytick={-1,-2,-3,-4,-5,-6,-7},%
        },%
    ]

        \addplot graphics [xmin=-1, xmax=1, ymin=-1, ymax=1]
        {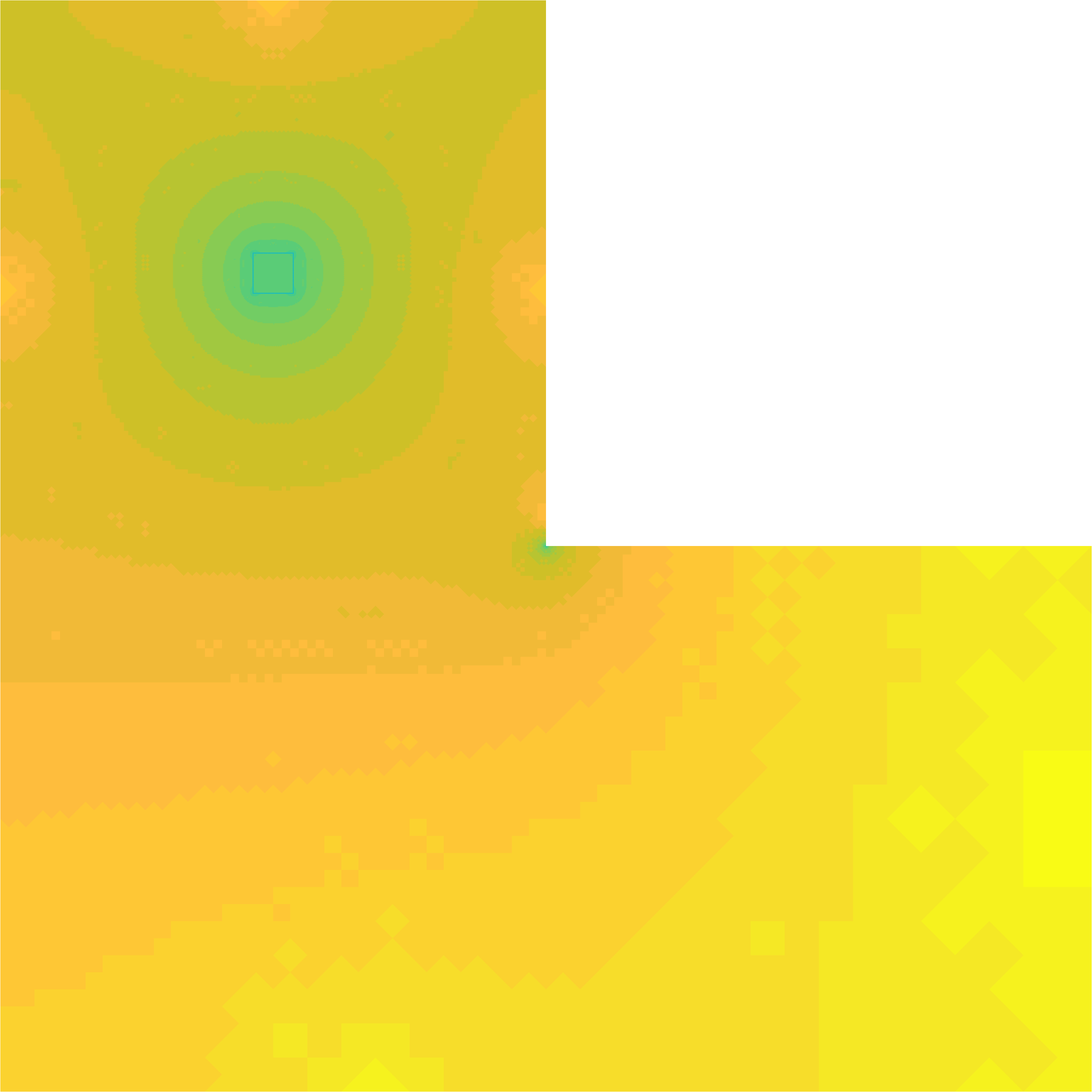};

    \end{axis}
\end{tikzpicture}
    }

    \subfloat[CALSFEM (\(8\,986\) triangles)]{%
        \label{subfig:microstructure_mesh_collective}
        \begin{tikzpicture}
    \begin{axis}[%
        axis equal image,%
        width=5.4cm,%
        xmin=-1.15, xmax=1.15,%
        ymin=-1.15, ymax=1.15,%
        font=\footnotesize,%
    ]

        \addplot graphics [xmin=-1, xmax=1, ymin=-1, ymax=1]
        {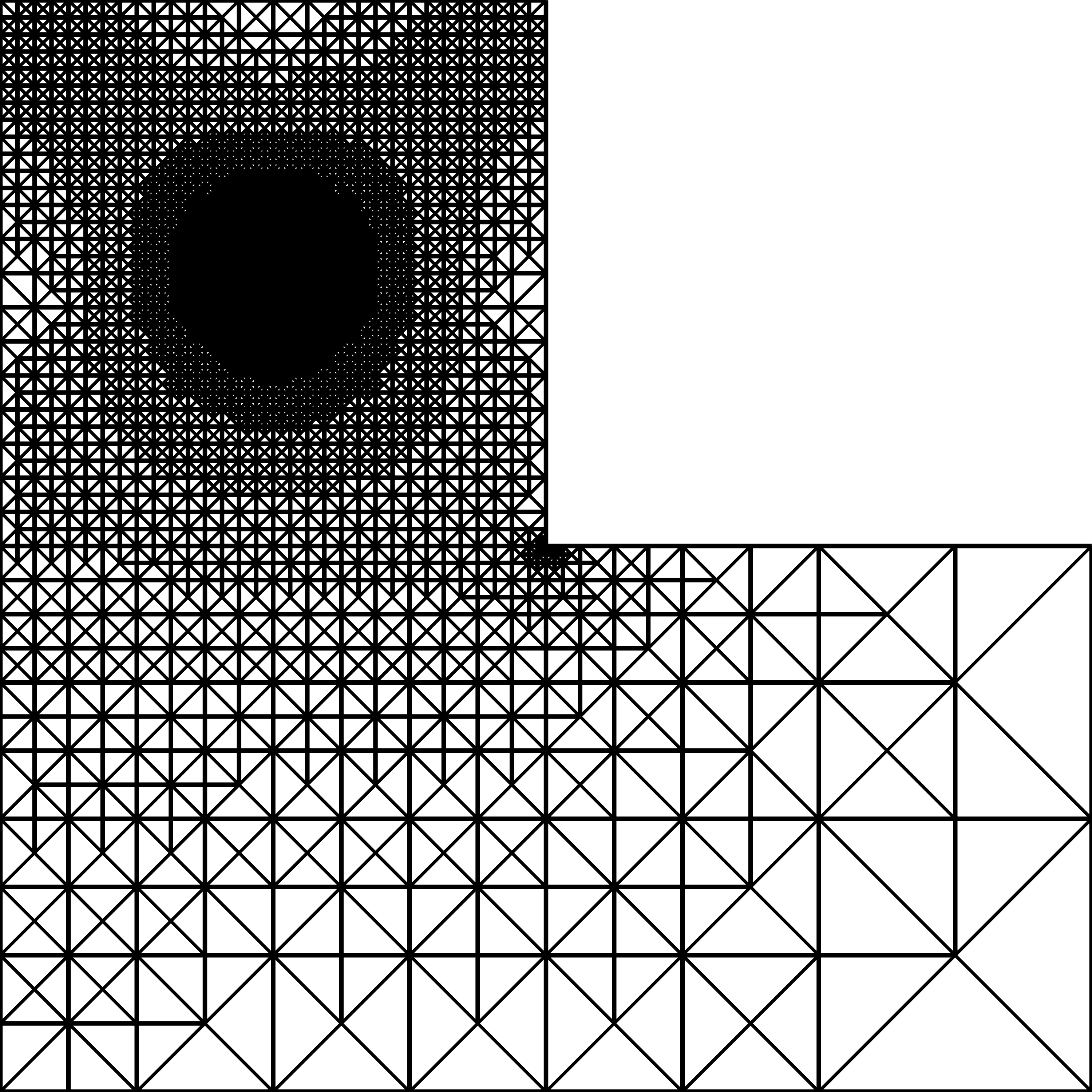};
    \end{axis}
\end{tikzpicture}
    }
    \hfil
    \subfloat[CALSFEM (\(5\,771\,531\) triangles)]{%
        \label{subfig:microstructure_meshsize_collective}
        \begin{tikzpicture}
    \begin{axis}[%
        axis equal image,%
        width=5.4cm,%
        xmin=-1.15, xmax=1.15,%
        ymin=-1.15, ymax=1.15,%
        font=\footnotesize,%
        point meta min=-6.62265990,%
        point meta max=-1.05360498,%
        colorbar,%
        colorbar style={%
            title={\(h_\ell\)},%
            font=\footnotesize,%
            width=2.5mm,%
            title style={yshift=-2mm},%
            yticklabel={$10^{\pgfmathprintnumber{\tick}}$},%
            ytick={-1,-2,-3,-4,-5,-6,-7},%
        },%
    ]

        \addplot graphics [xmin=-1, xmax=1, ymin=-1, ymax=1]
        {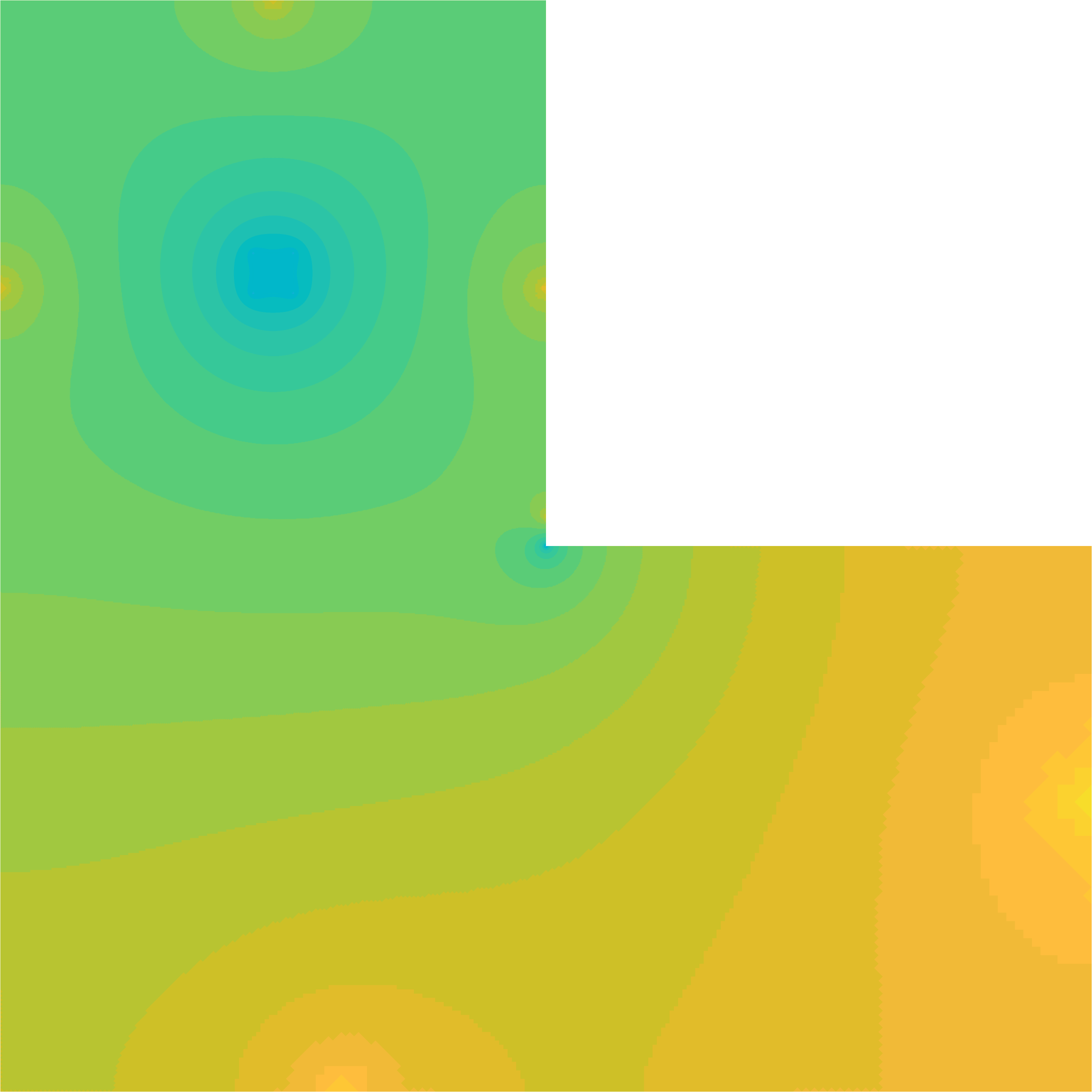};

    \end{axis}
\end{tikzpicture}
    }

    \subfloat[SALSFEM (\(132\,925\) triangles)]{%
        \label{subfig:microstructure_mesh_separate}
        \begin{tikzpicture}
    \begin{axis}[%
        axis equal image,%
        width=5.4cm,%
        xmin=-1.15, xmax=1.15,%
        ymin=-1.15, ymax=1.15,%
        font=\footnotesize,%
    ]

        \addplot graphics [xmin=-1, xmax=1, ymin=-1, ymax=1]
        {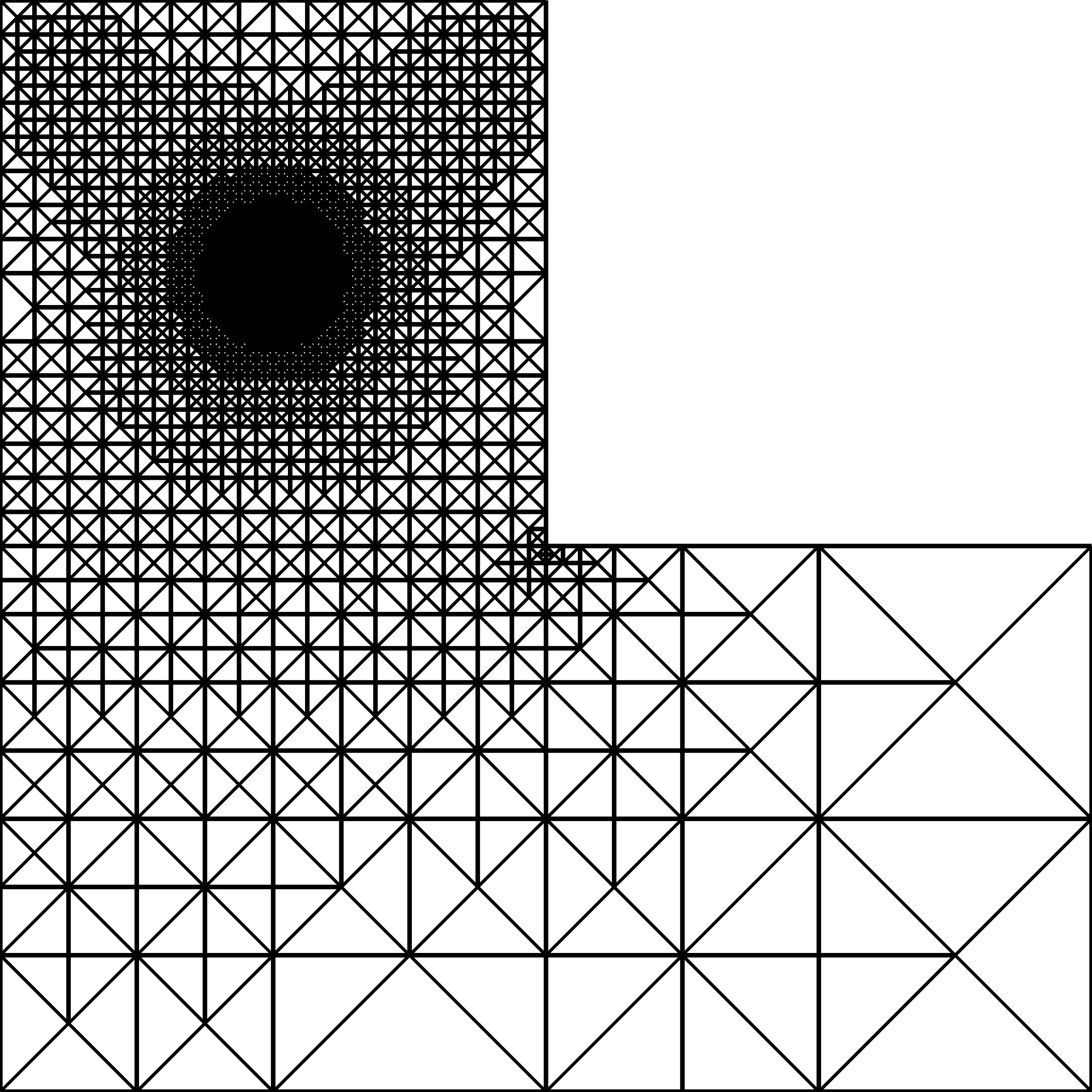};

    \end{axis}
\end{tikzpicture}
    }
    \hfil
    \subfloat[SALSFEM (\(3\,776\,010\) triangles)]{%
        \label{subfig:microstructure_meshsize_separate}
        \begin{tikzpicture}
    \begin{axis}[%
        axis equal image,%
        width=5.4cm,%
        xmin=-1.15, xmax=1.15,%
        ymin=-1.15, ymax=1.15,%
        font=\footnotesize,%
        point meta min=-6.62265990,%
        point meta max=-1.05360498,%
        colorbar,%
        colorbar style={%
            title={\(h_\ell\)},%
            font=\footnotesize,%
            width=2.5mm,%
            title style={yshift=-2mm},%
            yticklabel={$10^{\pgfmathprintnumber{\tick}}$},%
            ytick={-1,-2,-3,-4,-5,-6,-7},%
        },%
    ]

        \addplot graphics [xmin=-1, xmax=1, ymin=-1, ymax=1]
        {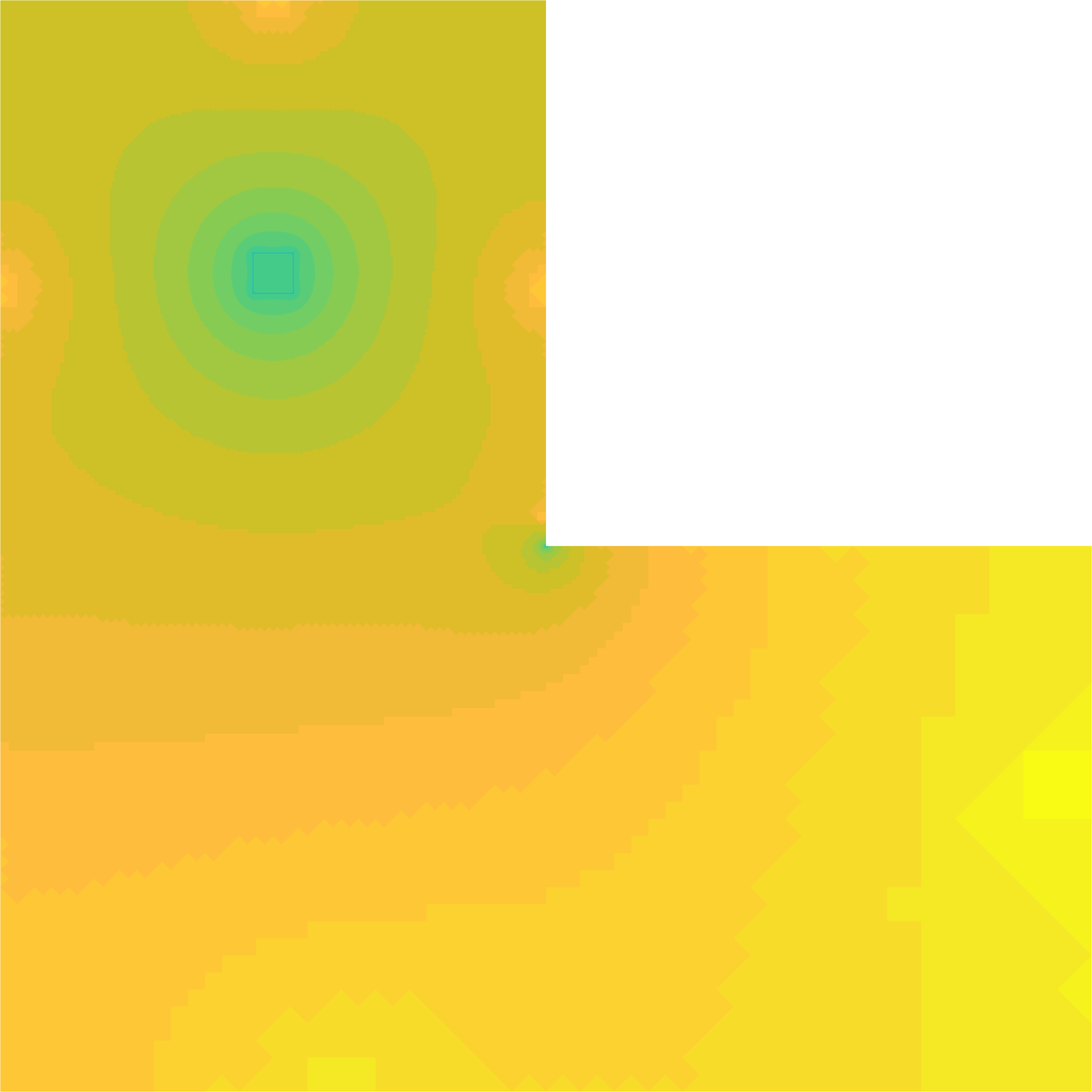};

    \end{axis}
\end{tikzpicture}
    }
    \caption{%
        Adaptively generated meshes by
        the three refinement strategies
        with parameters \(\theta = 0.3\),
        \(\kappa = 1\), and \(\rho = 0.8\)
        for the benchmark problem
        from Subsection~\ref{subsec:LshapeMicrostructure}
        with \(\epsilon = 3^{-3}\).
        The same color scale enables the comparison
        of the mesh-size \(h_\ell\vert_T \equiv \vert T \vert^{1/2}\)
        for \(T \in \T_\ell\)
        in the Subfigures~(b), (d), and~(f).
    }
    \label{fig:microstructure_meshes}
\end{figure}

The mesh plots in Figure~\ref{fig:microstructure_meshes}
demonstrate the different behaviour of the adaptive algorithms.
The Figures~\ref{subfig:microstructure_mesh_natural},
\ref{subfig:microstructure_mesh_collective},
and~\ref{subfig:microstructure_mesh_separate}
present the first level with an observable refinement towards
the reentrant corner.
While NALSFEM and SALSFEM focus on the microstructure up to
more than \(10^5\) triangles,
the CALSFEM already increases the refinement at the origin
at about \(7.4\times 10^3\) triangles.
The consideration of the very fine levels
of more than \(5.5 \times 10^5\) triangles
exhibits a highly adaptive refinement at the boundary of
the microstructure and the singularity at the reentrant corner.
However, the mesh of the CALSFEM appears more uniformly
with a larger minimal mesh-size
leading to the suboptimal convergence behaviour of the overall error
as seen in the convergence history plot in
Figure~\ref{subfig:microstructure_convergence}.

\begin{figure}
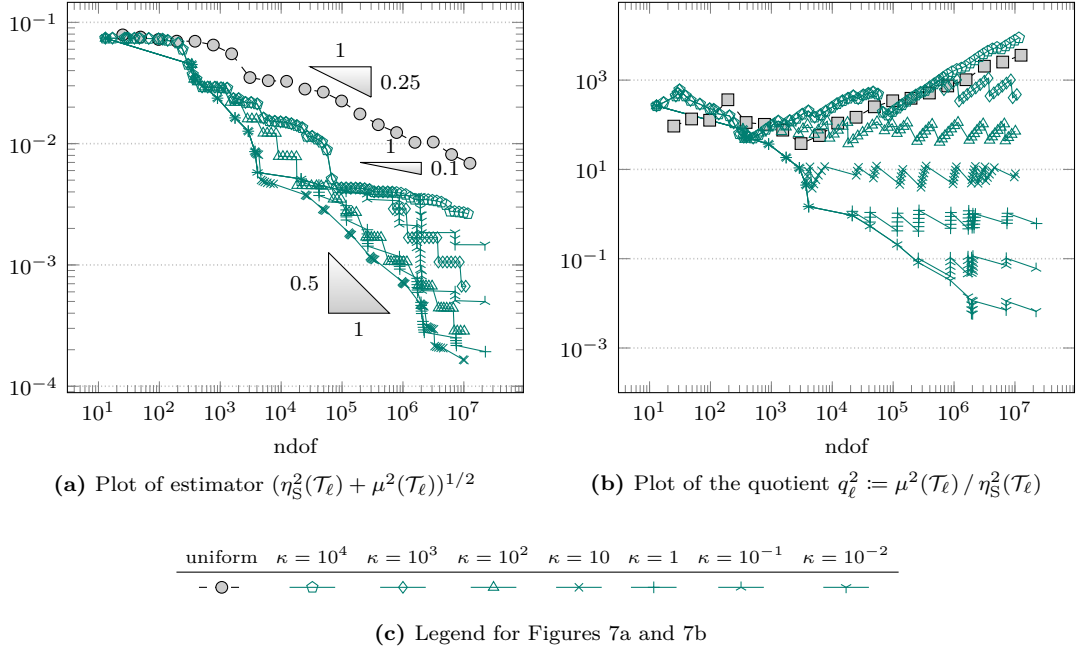

    \centering
    \subfloat[%
        Plot of estimator
        \((\eta_\SEP^2(\T_\ell) + \mu^2(\T_\ell))^{1/2}\)
    ]{%
        \label{subfig:microstructure_kappa}
        \input{figures/Poisson_LshapeMicrostructure_LS_adap30_sep_kappa.tex}
    }
    \subfloat[%
        Plot of the quotient
        \(q_\ell^2 \coloneqq
        \mu^2(\T_\ell)\, /\, \eta_{\textup{S}}^2(\T_\ell)\)
    ]{%
        \label{subfig:microstructure_kappa_quotient}
        \input{figures/Poisson_LshapeMicrostructure_LS_adap30_sep_kappa_quotient.tex}
    }

    \subfloat[Legend for
    Figures~\ref{subfig:microstructure_kappa}
    and~\ref{subfig:microstructure_kappa_quotient}]{
        \label{subfig:legend_kappa}
        \begin{tikzpicture}
    \pgfplotsset{%
        eta/.style={%
            mark=*,%
            every mark/.append style={%
                solid,%
                scale=1.1,%
            }%
        },%
        res/.style={%
            mark=square*,
            every mark/.append style={%
                solid,%
                scale=1.1,%
            }%
        },%
        resL2/.style={%
            mark=halfsquare left*,%
            every mark/.append style={%
                solid,%
                scale=1.1,%
            }%
        },%
        resDiv/.style={%
            mark=halfsquare right*,%
            every mark/.append style={%
                solid,%
                scale=1.1,%
            }%
        },%
        mu/.style={%
            mark=diamond*,%
            every mark/.append style={%
                solid,%
                scale=1.1,%
            }%
        },%
        cases/.style={%
            mark=,%
            only marks,%
            HUgreen,%
            nodes near coords,%
            point meta=explicit symbolic,%
            every node near coord/.style={%
                font=\footnotesize%
            }%
        },%
        uniform/.style={%
            dashed,%
            black,%
            every mark/.append style={fill=black!20!white}%
        },%
        natural/.style={%
            dashed,%
            HUblue,%
            every mark/.append style={%
                HUblue,%
                fill=white,%
            }%
        },%
        collective/.style={%
            dotted,%
            HUred,%
            every mark/.append style={%
                HUred,%
                fill=HUred!33!white,%
            }%
        },%
        separate/.style={%
            solid,%
            HUgreen,%
            every mark/.append style={
                HUgreen,%
                fill=HUgreen!66!white,%
            }%
        },%
    }

    \matrix [
        matrix of nodes,
        anchor = south,
        font = \scriptsize,
    ] at (0,0) {
        uniform &
        \(\kappa = 10^4\) &
        \(\kappa = 10^3\) &
        \(\kappa = 10^2\) &
        \(\kappa = 10\) &
        \(\kappa = 1\) &
        \(\kappa = 10^{-1}\) &
        \(\kappa = 10^{-2}\)
        \\
        \hline \\
        \ref{leg:sep:eta:unif} &
        \ref{leg:sep:eta:kappa10_4} &
        \ref{leg:sep:eta:kappa10_3} &
        \ref{leg:sep:eta:kappa10_2} &
        \ref{leg:sep:eta:kappa10} &
        \ref{leg:sep:eta:kappa1} &
        \ref{leg:sep:eta:kappa10_-1} &
        \ref{leg:sep:eta:kappa10_-2} \\
    };

\end{tikzpicture}
    }
    \caption{%
        Comparison of various choices
        for the separation parameter \(0 < \kappa\)
        in adaptive mesh-refinement with SALSFEM
        for the benchmark problem
        from Subsection~\ref{subsec:LshapeMicrostructure}
        with \(\epsilon = 3^{-3}\).
    }
\end{figure}

The convergence result in Theorem~\ref{thm:convergence_separate}
requires the separation parameter \(0 < \kappa < \kappa_0\)
to be sufficiently small.
The theoretical upper bound
\(
    \kappa_0
    =
    \min\{\widetilde\kappa, C_\textup{stab}^{-2} C_\textup{drel}^{-1}\}
\)
from \cite[Thm.~2.1]{MR3719030} incorporates two conditions.
If the data error \(\mu(\T)\) is monotonically decreasing
under mesh refinement (i.e., \(\Lambda_6 = 1\) in \cite{MR3719030}),
the estimator reduction in \cite[Thm.~4.1]{MR3719030}
and thus the plain convergence in \cite[Thm.~4.2]{MR3719030}
hold for arbitrary \(0 < \kappa < \infty\).
Hence, \(\widetilde\kappa = \infty\).
The proof of optimal convergence rates in \cite[Sect.~4.3]{MR3719030}
requires \(\kappa < C_\textup{stab}^{-2} C_\textup{drel}^{-1}\).
The estimates~\eqref{eq:constants_bounds}
for the Courant FEM with right-isosceles triangles lead to
\(\kappa_0 \geq 2.6 \times 10^{-6}\).
Despite this pessimistic theoretical bound,
the convergence rate of SALSFEM is optimal for the large range of
\(10^{-2} \leq \kappa \leq 10^2\) in practice
as displayed in Figure~\ref{subfig:microstructure_kappa}.
This suggests that the algorithm is fairly robust to the choice of
the parameter \(\kappa\).
Solely very large values exhibit suboptimal convergence rates.
For \(\kappa = 10^4\), every iteration carries out Case~A with
D\"orfler marking for the alternative estimator.
Hence, every larger value \(\kappa \geq 10^4\)
leads to exactly the same behaviour.

Figure~\ref{subfig:microstructure_kappa_quotient} displays
the quotient \(q_\ell^2 \coloneqq \mu^2(\T_\ell) / \eta_\textup{S}(\T_\ell)\)
used for the decision of the refinement strategy in the separate marking.
If this value is above the threshold \(\kappa\),
Case~B holds and the data approximation algorithm is carried out,
otherwise in Case~A, the D\"orfler marking for
the alternative estimator \(\eta_\textup{S}\) and NVB apply.
The reduction of \(\eta_\textup{S}\)
has rather no influence on the data error \(\mu\) and,
thus, leads to an increase of the quotient up to the threshold.
This reveals that, throughout the computation,
the SALSFEM algorithm ensures some balance
of error estimator and data error specified
by the parameter \(\kappa\).
For large quotients in the regime of the uniform refinement,
solely Case~A refinement is carried out
leading to the highly suboptimal convergence rate of about \(0.1\)
in Figure~\ref{subfig:microstructure_kappa}.
This suggests a choice of \(\kappa\) considerably smaller
than the values of the quotient \(q_\ell^2\) for uniform refinement.

As expected from the theoretical convergence result in
Theorem~\ref{thm:convergence_separate},
Figure~\ref{fig:microstructure_rho}
approves that
the choice of the parameter \(0 < \rho < 1\)
has no influence on the optimal convergence rate.
However, the reduction of the parameter \(\rho\)
decreases the number of solution steps significantly.
While \cite{Rabus2015} suggests a relatively small \(\rho\)
of about \(0.1\) for best overall performance,
a value close to one allows a more sensitive behaviour in the distinction
of the two refinement cases.
Accordingly, the choice of \(\rho =0.8\) in the remaining experiments
is preferable for an informative numerical comparison.
\begin{figure}
    \begin{minipage}{0.48\textwidth}
        \input{figures/Poisson_LshapeMicrostructure_LS_adap30_sep1_rho.tex}
        \captionof{figure}{%
            Comparison of estimator
            \((\eta_\SEP^2(\T_\ell) + \mu^2(\T_\ell))^{1/2}\)
            in SALSFEM for various parameters \(0 < \rho < 1\)
            with \(\theta = 0.3\) and \(\kappa = 1\)
            for the benchmark problem
            from Subsection~\ref{subsec:LshapeMicrostructure}
            with \(\epsilon = 3^{-3}\).
        }
        \label{fig:microstructure_rho}
    \end{minipage}
    \hfill
    \begin{minipage}{0.48\textwidth}
        \begin{tikzpicture}
    \colorlet{colMu}{TUyellow}
    \pgfplotsset{%
        mu/.style={colMu, mark=diamond*, every mark/.append style={solid,scale=1.1,fill=colMu!60!white}},%
    }
    %
    \begin{loglogaxis}[%
            width            = 0.78\textwidth,%
            xlabel           = ndof,%
            ylabel           = {data error \(\mu(\T_\ell)\)},%
            ymin             = 1e-6,%
            ymajorgrids      = true,%
            font             = \footnotesize,%
            legend style     = {
                legend columns = 1,
                legend pos     = south west,
                font           = \footnotesize
            }%
        ]
        \addlegendimage{mu};
        \addlegendentry{\(\epsilon = 3^{-1}\)};
        \addlegendimage{mu,mark=square*};
        \addlegendentry{\(\epsilon = 3^{-3}\)};
        \addlegendimage{mu,mark=triangle*};
        \addlegendentry{\(\epsilon = 3^{-5}\)};
        \addlegendimage{mu,mark=pentagon*};
        \addlegendentry{\(\epsilon = 3^{-7}\)};

        \addplot+ [mu, forget plot] table [x=ndof, y=mu]
        {
ndof	refinementTime	oscF	mu
6.00000000e+00	0.00000000e+00	3.51364184e-01	4.96903995e-01
2.30000000e+01	2.83220000e-02	9.21284664e-02	3.68513866e-01
8.20000000e+01	3.24350000e-02	3.47222222e-02	2.77777778e-01
1.58000000e+02	3.84870000e-02	2.14042153e-02	2.42161052e-01
3.50000000e+02	1.02782000e-01	9.74389944e-03	1.77864562e-01
8.62000000e+02	2.30114000e-01	3.63133692e-03	1.28049229e-01
1.79000000e+03	4.37112000e-01	1.31109440e-03	9.08103947e-02
3.79800000e+03	8.71205000e-01	4.58754734e-04	6.41186990e-02
7.67000000e+03	1.80385600e+00	1.63042742e-04	4.53719943e-02
1.55660000e+04	3.48001800e+00	5.74945236e-05	3.20710993e-02
3.12140000e+04	5.41905500e+00	2.03538805e-05	2.26818448e-02
6.26620000e+04	1.07517580e+01	7.19150098e-06	1.60370180e-02
1.25414000e+05	2.42071920e+01	2.54340744e-06	1.13404033e-02
2.51070000e+05	5.02234080e+01	8.99083985e-07	8.01869256e-03
5.02238000e+05	1.28791563e+02	3.17900062e-07	5.67013678e-03
1.00472600e+06	3.70074231e+02	1.12390072e-07	4.00936923e-03
        };

        \addplot+ [mu, mark=square*, forget plot] table [x=ndof, y=mu]
        {
ndof	refinementTime	oscF	mu
6.00000000e+00	0.00000000e+00	5.22343834e-02	7.38705735e-02
1.46000000e+02	5.22890000e-02	2.63692670e-03	5.96668400e-02
1.70000000e+02	1.13280000e-02	1.39448365e-03	4.46234767e-02
1.94000000e+02	1.20640000e-02	7.42912116e-04	3.36203645e-02
3.38000000e+02	5.05920000e-02	4.15662223e-04	2.86926852e-02
4.50000000e+02	4.44930000e-02	1.82763495e-04	2.33937273e-02
8.74000000e+02	1.48111000e-01	6.63692890e-05	1.58813773e-02
1.44200000e+03	1.94659000e-01	3.11756139e-05	1.21142537e-02
1.85800000e+03	1.49764000e-01	1.70170421e-05	7.76985478e-03
1.97800000e+03	5.87060000e-02	1.09178835e-05	6.29306003e-03
2.05000000e+03	4.74450000e-02	6.16418923e-06	4.46334977e-03
6.64200000e+03	1.53542000e+00	1.27556816e-06	3.59099793e-03
1.41940000e+04	2.52158900e+00	3.98678262e-07	2.89990678e-03
3.87780000e+04	8.33795000e+00	1.26723316e-07	2.07623482e-03
6.49380000e+04	9.35435400e+00	6.66954987e-08	1.49469830e-03
1.29274000e+05	2.47087370e+01	1.96402177e-08	1.05615232e-03
2.64674000e+05	5.48018290e+01	8.98578753e-09	7.48814729e-04
4.84874000e+05	1.05033498e+02	3.94089855e-09	6.55620397e-04
8.47018000e+05	2.63000463e+02	1.82540479e-09	4.37236453e-04
9.65506000e+05	2.11732426e+02	1.03159877e-09	2.82534891e-04
1.00966600e+06	2.15545647e+02	5.31832025e-10	1.97164810e-04
        };

        \addplot+ [mu, mark=triangle*, forget plot] table [x=ndof, y=mu]
        {
ndof	refinementTime	oscF	mu
6.00000000e+00	0.00000000e+00	5.81961178e-03	8.23017390e-03
2.90000000e+02	1.01676000e-01	3.86481185e-05	6.99605675e-03
3.14000000e+02	1.25640000e-02	2.15853518e-05	5.52585005e-03
3.38000000e+02	1.31000000e-02	6.78965269e-06	2.45811684e-03
5.38000000e+02	6.95860000e-02	2.10154569e-06	2.11147940e-03
7.14000000e+02	6.63940000e-02	8.84638666e-07	1.76513925e-03
1.06600000e+03	1.26971000e-01	3.56639248e-07	1.41111043e-03
1.81800000e+03	2.54997000e-01	2.49819874e-07	1.10603177e-03
3.06600000e+03	4.41614000e-01	8.08257191e-08	7.49379182e-04
7.01000000e+03	1.36316400e+00	4.03375644e-08	6.05340904e-04
1.25300000e+04	1.84009200e+00	1.09636596e-08	4.02554970e-04
1.75140000e+04	1.86172200e+00	5.61701537e-09	3.07721627e-04
2.49700000e+04	2.87256500e+00	2.62496874e-09	2.48165458e-04
6.13140000e+04	1.29557600e+01	1.24805900e-09	1.93682507e-04
9.97780000e+04	1.39569510e+01	4.82104285e-10	1.46766049e-04
1.68682000e+05	2.83988590e+01	2.39169706e-10	1.14389242e-04
2.35250000e+05	3.14686540e+01	1.05165955e-10	9.16811893e-05
3.77986000e+05	6.48830690e+01	7.46105540e-11	7.68565514e-05
4.16618000e+05	3.63987790e+01	5.57519496e-11	6.72667086e-05
5.65314000e+05	1.11949391e+02	2.35198129e-11	4.93246227e-05
1.25423400e+06	4.75430029e+02	1.41785255e-11	4.23508355e-05
        };

        \addplot+ [mu, mark=pentagon*, forget plot] table [x=ndof, y=mu]
        {
ndof	refinementTime	oscF	mu
6.00000000e+00	0.00000000e+00	6.46645163e-04	9.14494359e-04
4.34000000e+02	1.67954000e-01	5.57991374e-07	8.08057824e-04
4.58000000e+02	1.67950000e-02	3.34611364e-07	6.85284074e-04
4.82000000e+02	1.56780000e-02	1.10769855e-07	3.20823769e-04
5.54000000e+02	3.65350000e-02	3.23050854e-08	2.64643260e-04
7.30000000e+02	7.35450000e-02	1.29940950e-08	2.12895253e-04
1.65000000e+03	3.25704000e-01	3.54870043e-09	1.18343094e-04
1.77800000e+03	7.05010000e-02	2.42553297e-09	8.39407893e-05
1.82600000e+03	3.73000000e-02	1.59002450e-09	5.97205169e-05
2.11400000e+03	1.18422000e-01	7.85103659e-10	4.23599973e-05
8.71400000e+03	2.19958800e+00	1.36676609e-10	3.58126522e-05
2.16420000e+04	4.35137100e+00	9.18291773e-11	2.80219156e-05
6.14980000e+04	1.37193780e+01	9.90043878e-12	2.07627250e-05
9.18500000e+04	1.09778750e+01	5.91963024e-12	1.41177518e-05
1.35874000e+05	1.81464090e+01	3.26199088e-12	1.11908469e-05
2.58000000e+05	4.74181400e+01	1.57484561e-12	9.48223307e-06
3.38104000e+05	3.61980360e+01	8.56463321e-13	8.09036323e-06
5.28176000e+05	8.87873940e+01	4.96277338e-13	6.84238328e-06
9.46722000e+05	2.99138892e+02	3.11717261e-13	5.18360105e-06
1.68905800e+06	8.58389470e+02	7.82226825e-14	3.73456282e-06
        };
        \drawslopetriangle[ST1]{0.5}{1e4}{7e-3}
    \end{loglogaxis}
\end{tikzpicture}
        \captionof{figure}{%
            Investigation of AA algorithm with \(\rho = 0.9\)
            for the given data \(f_\epsilon\) with varying
            micro\-structure parameter \(0 < \epsilon < 1/2\) from
            Sub\-section~\ref{subsec:LshapeMicrostructure}.
            \\
        }%
        \label{fig:microstructure_approx_optimality}
    \end{minipage}
\end{figure}

\begin{figure}
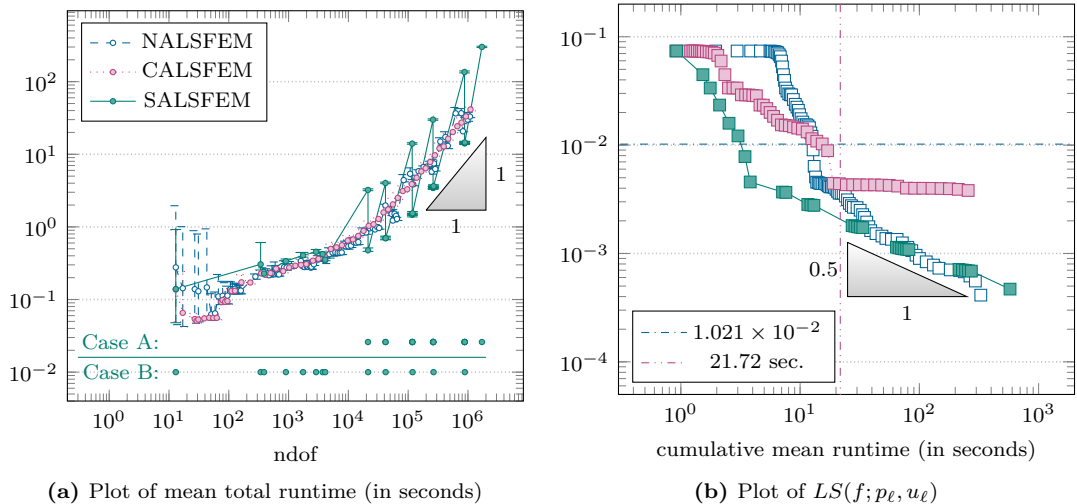

    \centering
    \subfloat[Plot of mean total runtime (in seconds)]{%
        \label{fig:microstructure_runtime}
        \input{figures/Poisson_LshapeMicrostructure_LS_adap30_runtime_inexact.tex}
    }
    \subfloat[Plot of \(LS(f; p_\ell, u_\ell)\)]{%
        \label{fig:microstructure_performance}
        \input{figures/Poisson_LshapeMicrostructure_LS_adap30_performance_inexact.tex}
    }
    \caption{%
        Comparison of the three adaptive strategies
        with parameters \(\theta = 0.3\), \(\kappa = 1\), and \(\rho = 0.8\)
        for the benchmark problem
        from Subsection~\ref{subsec:LshapeMicrostructure}
        with \(\epsilon = 3^{-3}\).
        Vertical error bars in Subfigure~(a)
        indicate the maximal and minimal measured time.
        Both figures employ the markers and
        line styles as introduced by the legend in
        Figure~\ref{subfig:legend_comparison}.
    }
    \label{fig:performance}
\end{figure}

In order to investigate the performance of the three algorithms,
Figure~\ref{fig:microstructure_runtime} displays
the mean total runtime in each iteration
from 10 independent runs of the adaptive loop.
All three refinement algorithms exhibit
almost linear complexity with respect to the number of
degrees of freedom. Note that the direct solution of the
algebraic linear system prevents linear complexity
of the overall implementation at hand.

The plot in Figure~\ref{fig:microstructure_performance}
displaying the estimator values versus the cumulative mean runtime
instead of the number of degrees of freedom
better represents the practical performance.
The adaptive algorithms ran up to \(10^5\) degrees of freedom.
In particular for the beginning of the computation
the SALSFEM is superior to the other refinement strategies.
This is because of the reduced number of solution steps
which may be further decreased by reducing the parameter \(\rho\).
Later NALSFEM and SALSFEM provide comparable results.

As a reference, Figure~\ref{fig:microstructure_performance}
displays the value
\(LS(f; p_\ell, u_\ell) \approx 1.02110264 \times 10^{-2}\)
resulting from a computation of a fine uniform mesh with
\(786\,432\) triangles
(\texttt{ndof}\({} = 1\,572\,865\)).
The solution and estimation took an average runtime of
about \(21.72\) seconds (without considering the time for the
generation of the fine mesh and the
computation of the alternative estimators).
The SALSFEM algorithm achieves the same accuracy already
after 3 seconds.
NALSFEM and even CALSFEM reach this
threshold after approximately 10 seconds
although the latter does not guarantee any
control of the data approximation error.
This is another striking evidence of the superiority of
adaptive mesh-refinement algorithms.

Finally,
Figure~\ref{fig:microstructure_approx_optimality}
confirms the quasi-optimality of the AA algorithm
\cite[axiom~(B1) in Sect.~2.4]{MR3719030}
with respect to the number of degrees of freedom.
\begin{figure}
    \centering
\end{figure}

\subsection{Waterfall benchmark}
\label{subsec:waterfall}

\begin{figure}
    \centering
    \subfloat[Solution \(u\)]{%
        \begin{tikzpicture}
    \pgfplotsset{/pgf/number format/fixed}

    \begin{axis}[%
        width=0.36\textwidth,%
        xmin=-0.1, xmax=1.1,%
        ymin=-0.1, ymax=1.1,%
        zmin=-0.0005, zmax=0.017,%
        font=\footnotesize,%
    ]
        \addplot3 graphics [%
            points={%
                (0,1,0) => (0,96-60.4)
                (0,0,0) => (61.6,0)
                (1,0,0) => (141.7,96-68.5)
                (0.5,0.5707,0.0158) => (66.5,96)
            }%
            ]
            {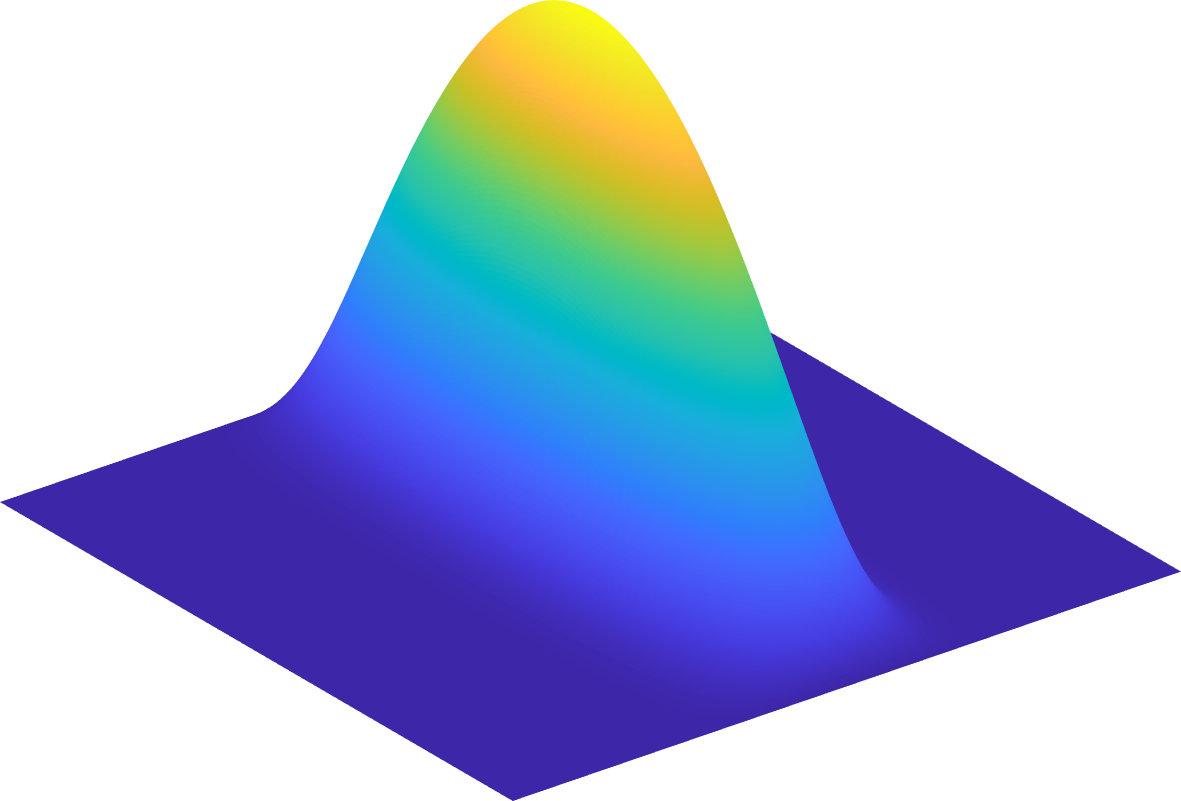};
    \end{axis}
\end{tikzpicture}
    }
    \hfil
    \subfloat[Right-hand side \(f\)]{%
        \label{subfig:waterfall_data}
        \begin{tikzpicture}
    \pgfplotsset{/pgf/number format/fixed}

    \begin{axis}[%
        width=0.36\textwidth,%
        xmin=-0.1, xmax=1.1,%
        ymin=-0.1, ymax=1.1,%
        zmin=-1.5, zmax=3.5,%
        font=\footnotesize,%
    ]
        \addplot3 graphics [%
            points={%
                (0,1,0) => (0,80.2-44.2)
                (0,0,0) => (61.6,0)
                (1,0,0) => (141.7,80.2-52.6)
                (0.5,0.587,3.38) => (64.4,80.2-1.1)
            }%
            ]
            {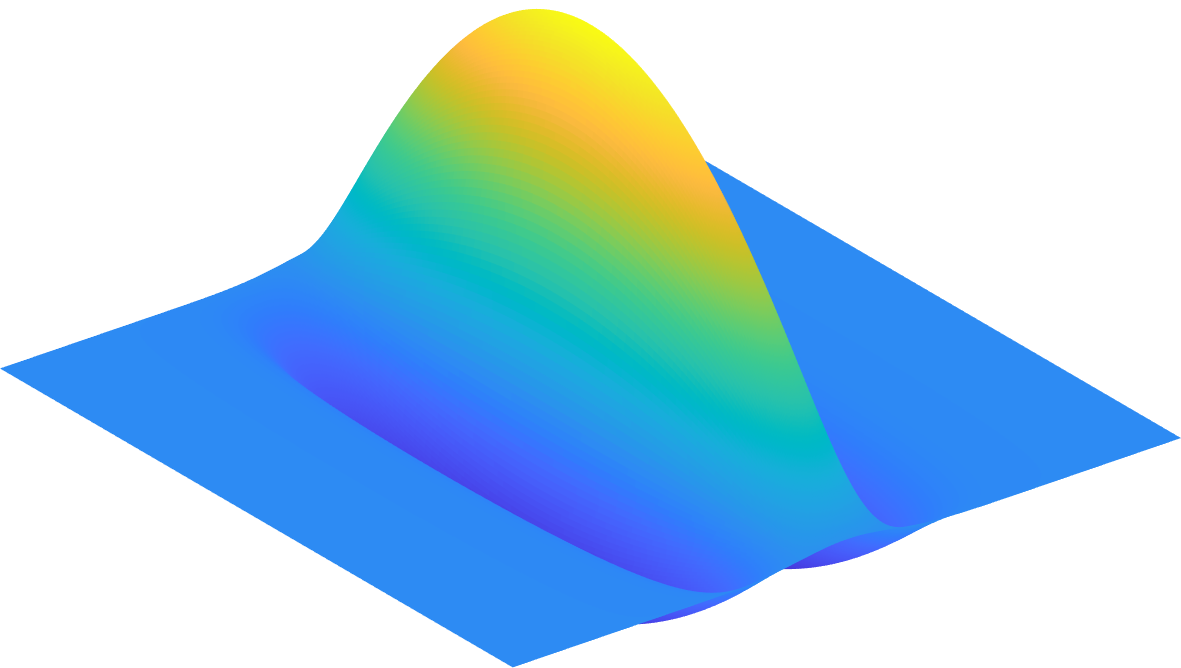};
    \end{axis}
\end{tikzpicture}
    }
    \caption{%
        Solution and right-hand side
        for the benchmark problem from Subsection~\ref{subsec:waterfall}.
    }
    \label{fig:waterfall_functions}
\end{figure}

This benchmark considers the exact solution
\(u \in H^1_0(\Omega)\) on the unit square
\(\Omega \coloneqq (0,1)^2\)
given in \cite[Sect.~4.2]{MR3279489} by
\begin{align*}
    u(x)
    \coloneqq
    x_1(x_1-1)x_2(x_2-1)\,
    \exp\!\big(
        -100(x_1 - 1/2)^2 -(x_2 - 117)^2 / 10000
    \big).
\end{align*}
The right-hand side is determined by \(f \coloneqq -\Delta u\).
Both functions are displayed in Figure~\ref{fig:waterfall_functions}.

For this benchmark with a smooth solution,
all adaptive algorithms exhibit optimal convergence rates
with a relatively small pre-asymptotic range.
Exemplarily, Figure~\ref{subfig:waterfall_convergence_natural}
presents the convergence graphs for the CALSFEM.
It confirms the equivalence of the estimator
with the exact error terms.
It is remarkable that even the data error
\(\mu^2(\T_\ell) \leq \Vert f + \ddiv p_\ell \Vert_{L^2(\Omega)}^2\)
converges with the optimal rate,
although this is not guaranteed by the theoretical convergence result.
Figure~\ref{subfig:waterfall_efficiency} displays
the efficiency indices of all three mesh-refinement
schemes.
The results illustrate the exactness of
the built-in error estimator \(LS(f; p_\ell, u_\ell)^{1/2}\)
from Theorem~\ref{thm:asymptotic}
already on the coarsest triangulations.
This is because the term
\(
    \Vert f + \ddiv p_\ell \Vert_{L^2(\Omega)}
    =
    \Vert \ddiv(p - p_\ell) \Vert_{L^2(\Omega)}
\)
dominates from the very beginning
in Figure~\ref{subfig:waterfall_convergence_natural}
and belongs to
the built-in error estimator \(LS(f; p_\ell, u_\ell)^{1/2}\)
and the error
\(e_\ell \coloneqq \vvvert (p - p_\ell, u - u_\ell) \vvvert\)
as well.
The fact that this dominating term
is not controlled by by the alternative estimator \(\eta_\COL\)
explains why the latter attains low efficiency indices only.
\begin{figure}
    \centering
    \subfloat[Plot of estimators and errors in
    CALSFEM]{%
        \input{figures/Poisson_Waterfall_LS_adap30_col.tex}
        \label{subfig:waterfall_convergence_natural}
    }
    \subfloat[Plot of efficiency indices]{%
        \label{subfig:waterfall_efficiency}
        \input{figures/Poisson_Waterfall_LS_adap30_comparison_efficiency.tex}
    }

    \subfloat[Legend of
    Figure~\ref{subfig:waterfall_efficiency}
    and~\ref{subfig:interface_efficiency}]{%
        \label{subfig:legend_efficiency}
        \begin{tikzpicture}
    \pgfplotsset{%
        eta/.style={%
            mark=*,%
            every mark/.append style={%
                solid,%
                scale=1.1,%
            }%
        },%
        res/.style={%
            mark=square*,
            every mark/.append style={%
                solid,%
                scale=1.1,%
            }%
        },%
        resL2/.style={%
            mark=halfsquare left*,%
            every mark/.append style={%
                solid,%
                scale=1.1,%
            }%
        },%
        resDiv/.style={%
            mark=halfsquare right*,%
            every mark/.append style={%
                solid,%
                scale=1.1,%
            }%
        },%
        mu/.style={%
            mark=diamond*,%
            every mark/.append style={%
                solid,%
                scale=1.1,%
            }%
        },%
        errL2Sigma/.style={%
            mark=triangle*,%
            every mark/.append style={%
                solid,%
                scale=1.1,%
            }%
        },%
        errGradU/.style={%
            mark=pentagon*,%
            every mark/.append style={%
                solid,%
                scale=0.9,%
            }%
        },%
        errL2U/.style={%
            mark=diamond*,%
            every mark/.append style={%
                solid,%
                scale=1.1,%
            }%
        },%
        cases/.style={%
            mark=,%
            only marks,%
            HUgreen,%
            nodes near coords,%
            point meta=explicit symbolic,%
            every node near coord/.style={%
                font=\footnotesize%
            }%
        },%
        uniform/.style={%
            dashed,%
            black,%
            every mark/.append style={fill=black!20!white}%
        },%
        natural/.style={%
            dashed,%
            HUblue,%
            every mark/.append style={%
                HUblue,%
                fill=white,%
            }%
        },%
        collective/.style={%
            dotted,%
            HUred,%
            every mark/.append style={%
                HUred,%
                fill=HUred!33!white,%
            }%
        },%
        separate/.style={%
            solid,%
            HUgreen,%
            every mark/.append style={
                HUgreen,%
                fill=HUgreen!66!white,%
            }%
        },%
    }

    %
    %
    \pgfplotstableset{%
        create on use/effEta/.style={%
            create col/expr={%
                \thisrow{eta} / sqrt(\thisrow{errSigma}^2 +
                \thisrow{errGradU}^2)
                }
        },
        create on use/effEtaMu/.style={%
            create col/expr={%
                sqrt(\thisrow{eta}^2 + \thisrow{mu}^2)
                / sqrt(\thisrow{errSigma}^2 + \thisrow{errGradU}^2)
                }
        },
        create on use/effRes/.style={%
            create col/expr={%
                \thisrow{res} / sqrt(\thisrow{errSigma}^2 +
                \thisrow{errGradU}^2)
            }
        }
    }

    \matrix [
        matrix of nodes,
        anchor = south,
        font = \scriptsize,
        column 1/.style={anchor=base east},
    ] at (0,0) {
        & NALSFEM
        & CALSFEM
        & SALSFEM \\
        \hline \\
        \(LS(f;p_\ell,u_\ell)^{1/2} / e_\ell\)
        & \ref{leg:eff:nat:eta}
        & \ref{leg:eff:col:res}
        & \ref{leg:eff:sep:res} \\
        Alternative estimator
        & \(\eta_\NAT / e_\ell\)
        & \(\eta_\COL / e_\ell\)
        & \((\eta_\SEP^2 + \mu^2)^{1/2} / e_\ell\) \\
        & \ref{leg:eff:nat:eta}
        & \ref{leg:eff:col:eta}
        & \ref{leg:eff:sep:eta} \\
    };
\end{tikzpicture}
    }

    \caption{%
        Convergence history plot and
        plot of the efficiency indices with respect to the error
        \(e_\ell \coloneqq \vvvert (p - p_\ell, u - u_\ell) \vvvert\)
        for the waterfall benchmark prob\-lem
        from Subsection~\ref{subsec:waterfall}.
        The parameters for the adaptive mesh-refinement
        strategies read \(\theta = 0.3\),
        \(\kappa = 1\), and \(\rho = 0.8\).
    }%
    \label{fig:waterfall_efficiency_comparison}
\end{figure}

\begin{figure}
    \centering
    \subfloat[NALSFEM (\(503\,428\) triangles)]{
        \label{subfig:waterfall_mesh_natural}
        \begin{tikzpicture}
    \begin{axis}[%
        axis equal image,%
        width=6cm,%
        xmin=-0.1, xmax=1.1,%
        ymin=-0.1, ymax=1.1,%
        font=\footnotesize,%
    ]

        \addplot graphics [xmin=0, xmax=1, ymin=0, ymax=1]
        {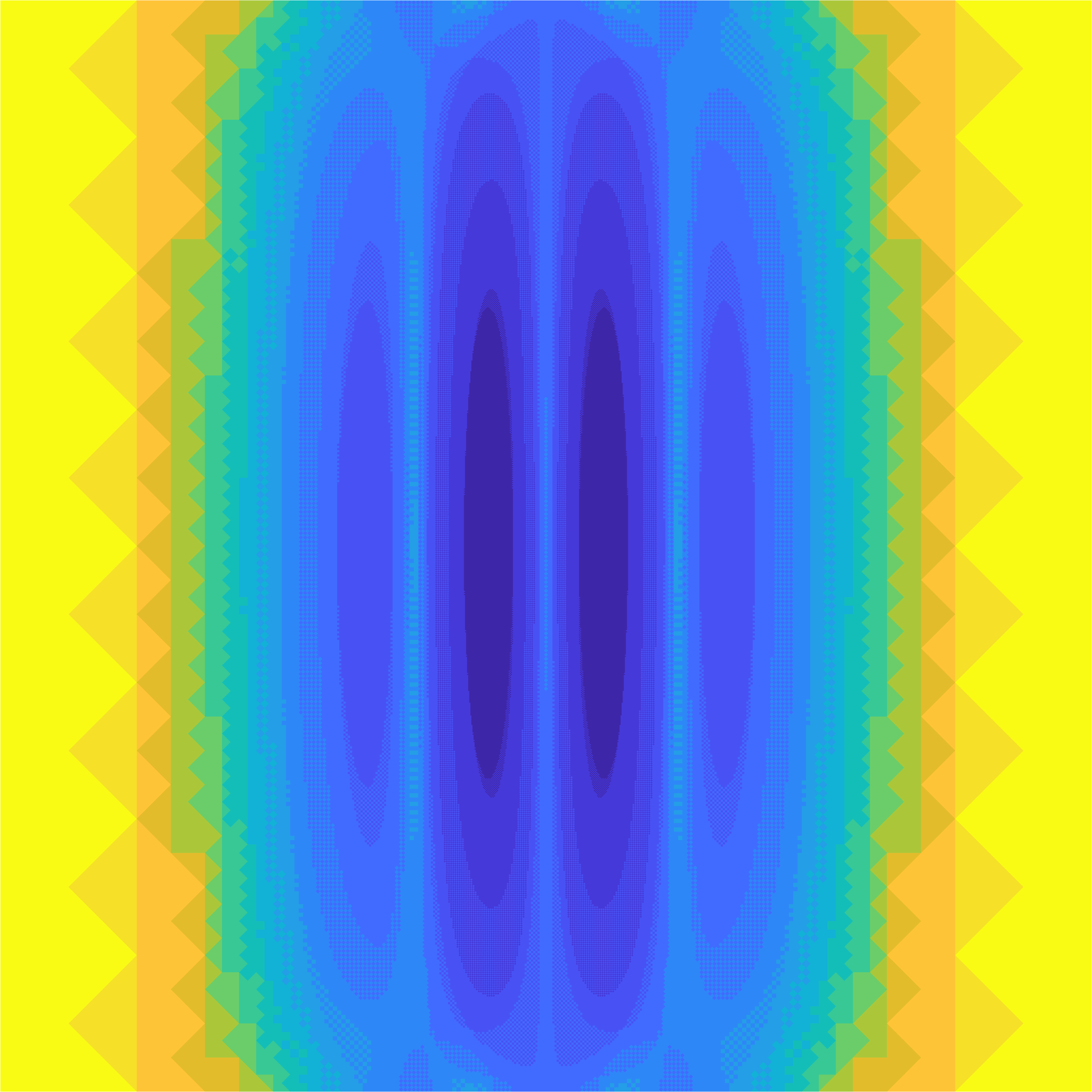};
    \end{axis}
\end{tikzpicture}
    }
    \hfil
    \subfloat[CALSFEM (\(531\,612\) triangles)]{
        \label{subfig:waterfall_mesh_collective}
        \begin{tikzpicture}
    \begin{axis}[%
        axis equal image,%
        width=6cm,%
        xmin=-0.1, xmax=1.1,%
        ymin=-0.1, ymax=1.1,%
        font=\footnotesize,%
        point meta min=-3.31132995,%
        point meta max=-1.20411998,%
        colorbar,%
        colorbar style={%
            title={\(h_\ell\)},%
            font=\footnotesize,%
            width=2.5mm,%
            title style={yshift=-2mm},%
            yticklabel={$10^{\pgfmathprintnumber{\tick}}$},%
        },%
    ]

        \addplot graphics [xmin=0, xmax=1, ymin=0, ymax=1]
        {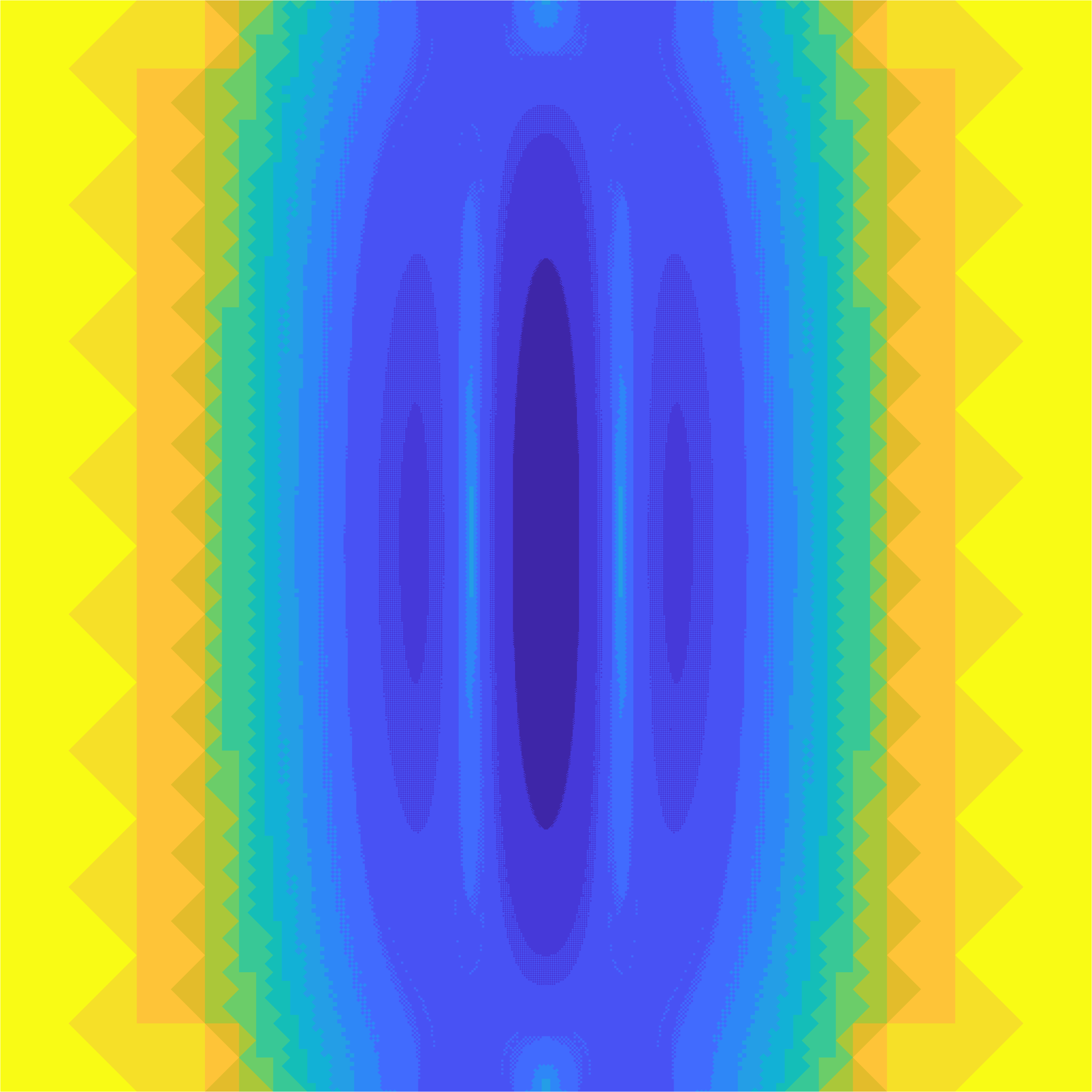};
    \end{axis}
\end{tikzpicture}
    }
    \medskip

    \subfloat[SALSFEM (\(610\,636\) triangles)]{
        \label{subfig:waterfall_mesh_separate}
        \begin{tikzpicture}
    \begin{axis}[%
        axis equal image,%
        width=6cm,%
        xmin=-0.1, xmax=1.1,%
        ymin=-0.1, ymax=1.1,%
        font=\footnotesize,%
    ]

        \addplot graphics [xmin=0, xmax=1, ymin=0, ymax=1]
        {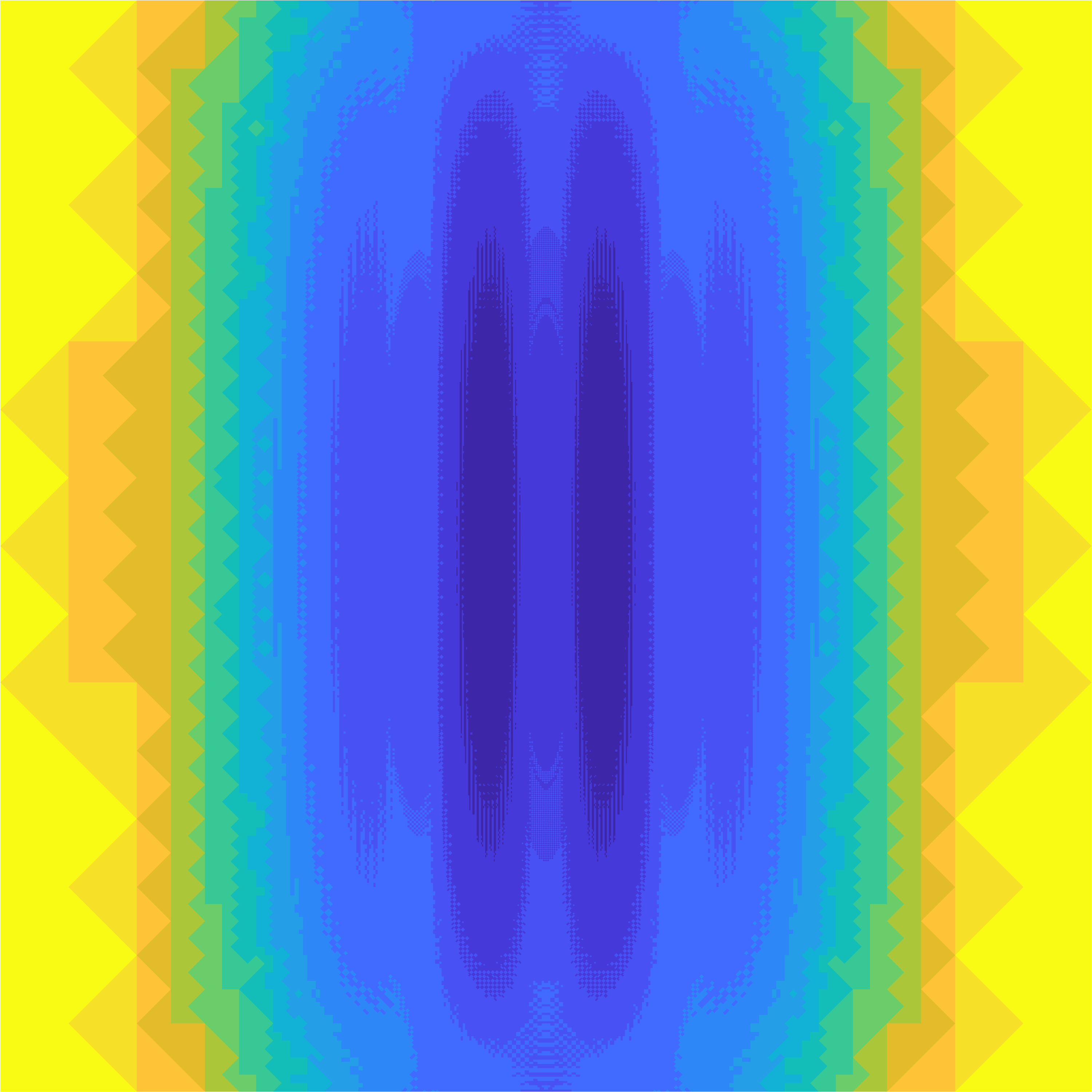};
    \end{axis}
\end{tikzpicture}
    }
    \hfil
    \subfloat[AA (\(565\,078\) triangles)]{
        \label{subfig:waterfall_mesh_approx}
        \begin{tikzpicture}
    \begin{axis}[%
        axis equal image,%
        width=6cm,%
        xmin=-0.1, xmax=1.1,%
        ymin=-0.1, ymax=1.1,%
        font=\footnotesize,%
        point meta min=-3.31132995,%
        point meta max=-1.20411998,%
        colorbar,%
        colorbar style={%
            title={\(h_\ell\)},%
            font=\footnotesize,%
            width=2.5mm,%
            title style={yshift=-2mm},%
            yticklabel={$10^{\pgfmathprintnumber{\tick}}$},%
        },%
    ]

        \addplot graphics [xmin=0, xmax=1, ymin=0, ymax=1]
        {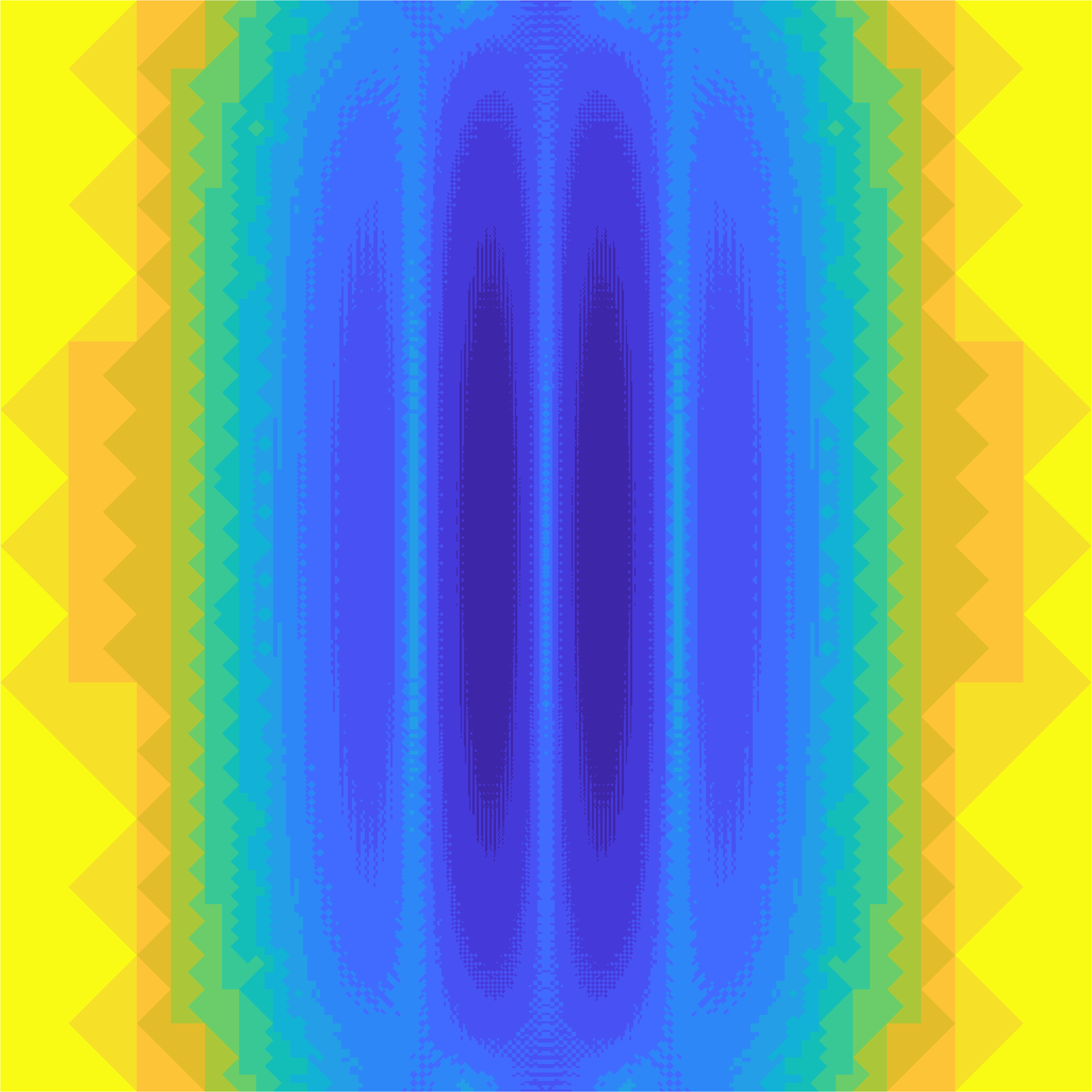};
    \end{axis}
\end{tikzpicture}
    }

    \caption{%
        Plots of adaptively refined meshes for the three refinement
        strategies with parameters \(\theta = 0.3\),
        \(\kappa = 10\), and \(\rho = 0.8\)
        for the benchmark problem from Subsection~\ref{subsec:waterfall}.
    }
    \label{fig:waterfall_meshes}
\end{figure}

The mesh plots in Figure~\ref{fig:waterfall_meshes}
illustrate the different behaviour of the adaptive algorithms.
The NALSFEM in Figure~\ref{subfig:waterfall_mesh_natural}
focusses on the regions with large gradients of the
right-hand side \(f\) (see Figure~\ref{subfig:waterfall_data})
in order to allow for a proper piecewise constant approximation.
The mesh is similar to the result from the data approximation by AA
in Figure~\ref{subfig:waterfall_mesh_approx}.
On the contrary,
the CALSFEM in Figure~\ref{subfig:waterfall_mesh_collective}
increases the refinement in regions with large absolute values
of \(f\).
The SALSFEM in Figure~\ref{subfig:waterfall_mesh_separate}
seemingly combines both aspects.

\subsection{Discontinuous coefficients}
\label{subsec:interface}

\begin{figure}
    \centering
    \subfloat[Solution \(u\)]{%
        \label{subfig:interface_solution}
        \begin{tikzpicture}
    \pgfplotsset{/pgf/number format/fixed}

    \begin{axis}[%
        width=0.36\textwidth,%
        xmin=-1.1, xmax=1.1,%
        ymin=-1.1, ymax=1.1,%
        zmin=-0.09, zmax=0.09,%
        font=\footnotesize,%
    ]
        \addplot3 graphics [%
            points={%
                (-1,-1,0.0812259) => (54.2,359.9-18.8)
                (1,-1,0) => (0.2,359.9-236.2)
                (-1,-0.00125156,0.078218) => (206.9,359.9-34.0)
                (1,1,-0.0812259) => (305.8,359.9-341.0)
            }%
            ]
            {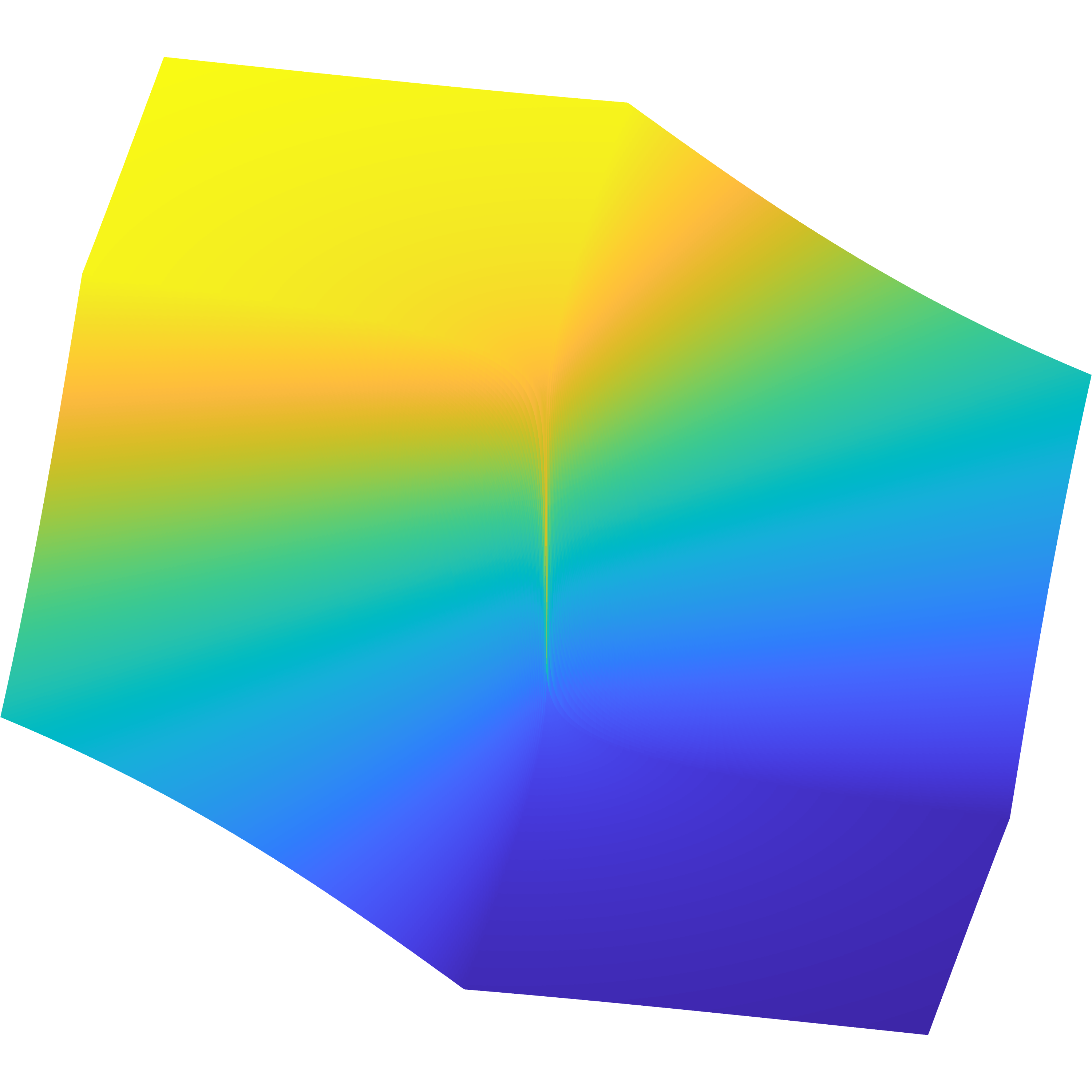};
    \end{axis}
\end{tikzpicture}
    }
    \hfil
    \subfloat[CALSFEM (\(5\,370\) triangles)]{
        \label{subfig:interface_adaptivity}
        \begin{tikzpicture}
    \begin{axis}[%
        axis equal image,%
        width=6cm,%
        xmin=-1.15, xmax=1.15,%
        ymin=-1.15, ymax=1.15,%
        font=\footnotesize,%
	point meta min=-6.02059991,%
        point meta max=-6.02059991e-01,%
        colorbar,%
        colorbar style={%
            title={\(h_\ell\)},%
            font=\footnotesize,%
            width=2.5mm,%
            title style={yshift=-2mm},%
            yticklabel={$10^{\pgfmathprintnumber{\tick}}$},%
            ytick={-1,-2,-3,-4,-5,-6},%
        },%
    ]

        \addplot graphics [xmin=-1, xmax=1, ymin=-1, ymax=1]
        {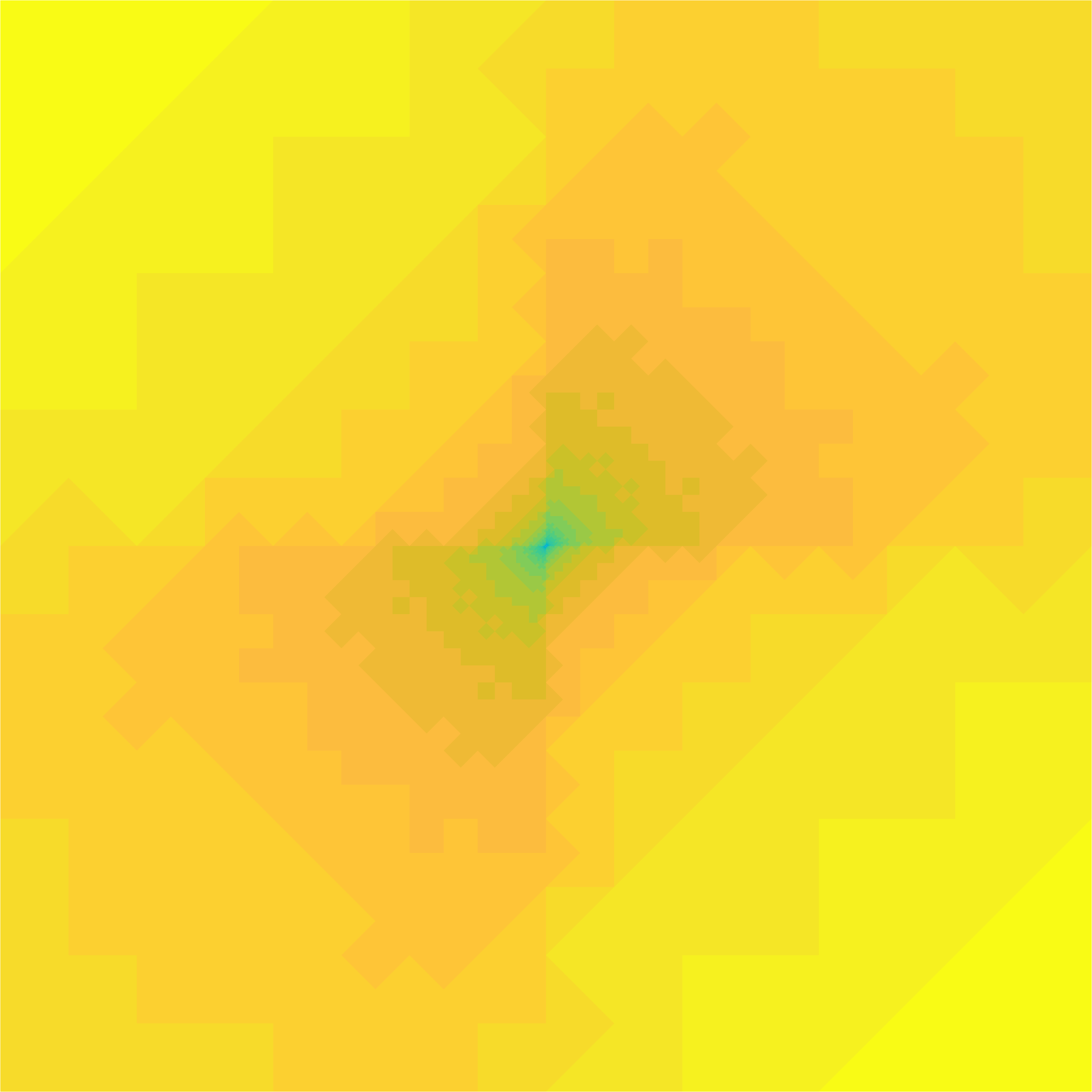};
    \end{axis}
\end{tikzpicture}
    }

    \caption{%
        Solution plot and
        plot of mesh-size
        \(h_\ell\vert_T \equiv \vert T \vert^{1/2}\)
        of adaptively refined mesh using CALSFEM
        with bulk parameter \(\theta = 0.7\)
        for the benchmark problem
        from Subsection~\ref{subsec:interface}.
    }
    \label{fig:interface:solution_mesh}
\end{figure}

The final benchmark considers the elliptic problem
with piecewise constant scalar diffusion coefficient
\(a \in L^\infty(\Omega)\), defined by
\[
    a(x) \coloneqq
    \begin{cases}
        a_1 & \text{if } 0 < x_1 x_2,
        \\
        a_2 & \text{if } x_1 x_2 < 0,
    \end{cases}
\]
and right-hand side \(f \equiv 0\)
on the square domain \(\Omega \coloneqq (-1, 1)^2\).
It seeks \((p, u) \in H(\ddiv, \Omega) \times H^1(\Omega)\) satisfying
\begin{equation}
    \label{eq:diffusion}
    f + \ddiv p = 0
    \quad\text{and}\quad
    a^{-1/2} p - a^{1/2} \nabla u = 0
    \quad\text{in } \Omega
    \quad\text{subject to}\quad
    u = u_{\textup{D}}
    \text{ on } \partial\Omega.
\end{equation}
The weighting of the second residual in~\eqref{eq:diffusion}
leads to the fundamental equivalence of the least-squares functional
\begin{equation}
    \label{eq:diffusion_functional}
    LS(f, a; p, u)
    \coloneqq
    \Vert f + \ddiv p \Vert_{L^2(\Omega)}^2
    +
    \Vert a^{-1/2} p - a^{1/2} \nabla u
    \Vert_{L^2(\Omega)}^2
\end{equation}
and the natural weighted \(H(\ddiv)\) and energy norm
\(
    \Vert \ddiv q \Vert_{L^2(\Omega)}^2
    +
    \Vert a^{-1/2} q \Vert_{L^2(\Omega)}^2
    +
    \Vert a^{1/2} \nabla u \Vert_{L^2(\Omega)}^2
\)
with equivalence constants solely depending
on the uniform lower bound of the diffusion coefficient.
Note that the inhomogeneous Dirichlet boundary conditions
\(u_{\textup{D}}\)
lead to an additional oscillation term in the estimators
\(\eta_\NAT\), \(\eta_\COL\), and \(\eta_\SEP\)
\cite{MR3715170,BringmannDissertation,MR4271577}
and the overall error
\begin{equation}
    \label{eq:diffusion_energy_norm}
    \begin{split}
        \vvvert (q, v) \vvvert_a^2
        &\coloneqq
        \Vert \ddiv q \Vert_{L^2(\Omega)}^2
        +
        \Vert a^{-1/2} q \Vert_{L^2(\Omega)}^2
        +
        \Vert a^{1/2} \nabla u \Vert_{L^2(\Omega)}^2
        \\
        &\phantom{{}\coloneqq{}}
        +
        \sum_{E \in \E(\partial\Omega)}
        \vert \omega_E \vert^{1/2}
        \Vert (1 - \Pi_{0,E}) \partial u_{\textup{D}} /\partial s
        \Vert_{L^2(E)}^2.
    \end{split}
\end{equation}
For some parameter \(0 < \gamma < 2\),
the exact weak solution to \eqref{eq:diffusion}
in polar coordinates from \cite{MR0393815}
reads \(u(r, \phi) \coloneqq r^\gamma \mu(\phi)\) and
\(p \coloneqq a \nabla u\) with
\[
    \mu(\phi)
    \coloneqq
    \begin{cases}
        \cos((\pi/2 - \sigma) \gamma)\,
        \cos((\phi - \pi/2 + \rho) \gamma)
        & \text{if } 0 \leq \phi < \pi/2,
        \\
        \cos(\rho \gamma)\,
        \cos((\phi - \pi + \sigma) \gamma)
        & \text{if } \pi/2 \leq \phi < \pi,
        \\
        \cos(\sigma \gamma)\,
        \cos((\phi - \pi - \rho) \gamma)
        & \text{if } \pi \leq \phi < 3\pi/2,
        \\
        \cos((\pi/2 - \rho) \gamma)\,
        \cos((\phi - 3\pi/2 - \sigma) \gamma)
        & \text{if } 3\pi/2 \leq \phi \leq 2\pi,
    \end{cases}
\]
and constants \(0 < \rho\) and \(\sigma < 0\).
The parameter \(\gamma\) determines the regularity of the
solution \(u \in H^{1+\gamma-\varepsilon}(\Omega)\)
for all \(0 < \varepsilon\).
The choice of \(\gamma = 0.1\) in \cite{MR1770058} leads to
the constants \(\rho = \pi/4\),
\(\sigma \approx -14.922\,565\,104\,551\,52\), the coefficients
\(a_1 \approx 161.447\,638\,797\,588\,1\), \(a_2 = 1\),
and the solution \(u\) displayed in
Figure~\ref{subfig:interface_solution}.
The nodal interpolation of the
exact solution \(u\) prescribes the inhomogeneous
boundary conditions \(u_{\textup{D}}\) in
the discrete minimization of~\eqref{eq:diffusion_functional}.

This benchmark problem models
intersecting interfaces with the difficulty of
a strong cross-point singularity at the origin.
Figure~\ref{subfig:interface_adaptivity} exhibits
the intense adaptive refinement of CALSFEM
towards the origin.
The heavy grading of the mesh leads to
ill-conditioned system matrices
already for a relatively small number of degrees of freedom,
e.g., from about \texttt{ndof} \({} = 4\,000\) for CALSFEM
with \(\theta = 0.7\).
For this reason, its unreliable results are omitted
in the figures.
Moreover, due to the lack of any data approximation error for
the right-hand side \(f\),
the results of CALSFEM and SALSFEM coincide.

\begin{figure}
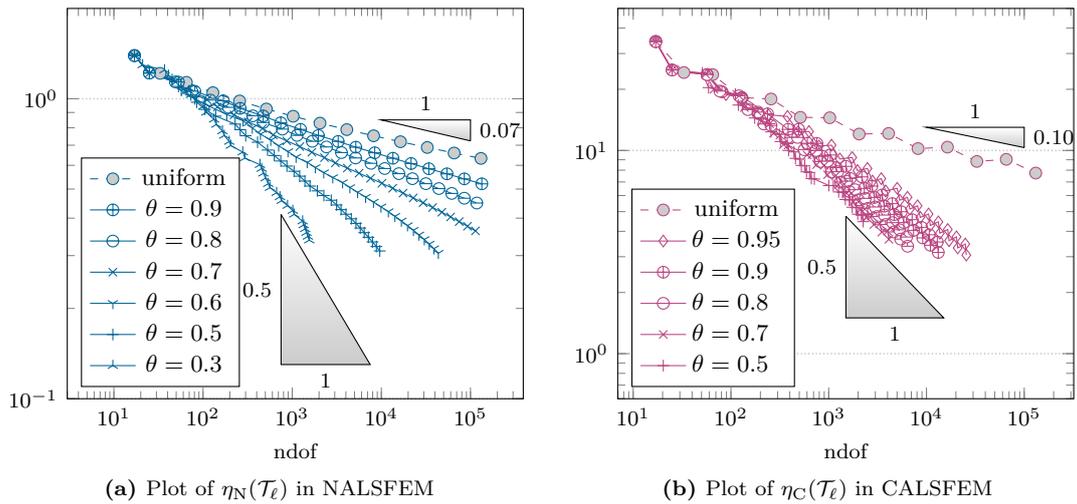

    \centering
    \subfloat[Plot of \(\eta_\NAT(\T_\ell)\) in NALSFEM]{%
        \input{figures/Diffusion_Kellogg_LS_adap_nat_theta.tex}
        \label{subfig:interface_theta_natural}
    }
    \hfil
    \subfloat[Plot of \(\eta_\COL(\T_\ell)\) in CALSFEM]{%
        \input{figures/Diffusion_Kellogg_LS_adap_col_theta.tex}
        \label{subfig:interface_theta_collective}
    }

    \caption{%
        Comparison of various choices for the bulk parameter
        \(0 < \theta \leq 1\) in the adaptive
        mesh-refinement strategies
        for the benchmark prob\-lem
        from Subsection~\ref{subsec:interface}.
    }%
    \label{fig:interface_theta}
\end{figure}

The plain convergence analysis for NALSFEM
from Theorem~\ref{thm:plain_convergence}
holds under general assumptions.
However, the convergence result in \cite{MR4138307}
requires nested discrete spaces
which is violated by the nodal interpolation
of the boundary data in the implementation at hand.
The analysis in \cite{MR4216839} covers inhomogeneous
boundary conditions if they are weakly enforced
by additional residuals in the least-squares functional.
Nevertheless, NALSFEM converges
for all choices of the bulk parameter
in Figure~\ref{subfig:interface_theta_natural}.
In contrast to Figure~\ref{subfig:Lshape_convergence_natural}
for the L-shaped domain benchmark,
the convergence of NALSFEM with the optimal rate
seems to require much smaller bulk parameters
\(\theta \leq 0.3\) in this benchmark problem.
In contrast to that,
CALSFEM appears to be much more robust with
respect to the choice of \(0 < \theta \leq 1\)
in Figure~\ref{subfig:interface_theta_collective}.

\begin{figure}
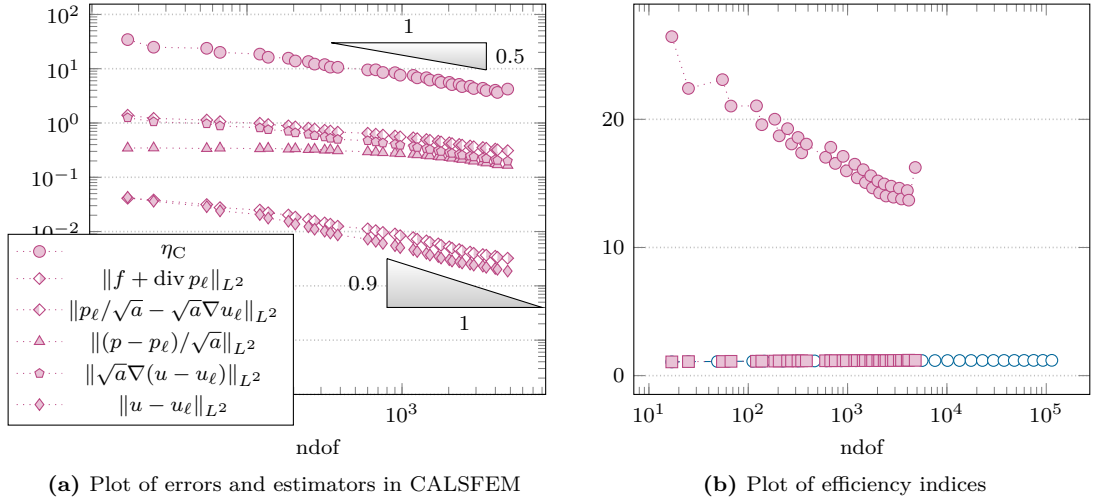

    \centering
    \subfloat[Plot of errors and estimators in CALSFEM]{%
        \input{figures/Diffusion_Kellogg_LS_adap70_col.tex}
        \label{subfig:interface_convergence_collective}
    }
    \hfil
    \subfloat[Plot of efficiency indices]{%
        \label{subfig:interface_efficiency}
        \input{figures/Diffusion_Kellogg_LS_adap70_comparison_efficiency.tex}
    }

    \caption{%
        Convergence history plot of CALSFEM
        with bulk parameter \(\theta = 0.7\)
        and efficiency indices with respect to the error
        \(e_\ell \coloneqq \vvvert (p - p_\ell, u - u_\ell) \vvvert_a\)
        from~\eqref{eq:diffusion_energy_norm}
        for the diffusion benchmark prob\-lem
        from Subsection~\ref{subsec:interface}.
        Figure~\ref{subfig:interface_efficiency}
        employs the graphs as introduced
        by the legend in
        Figure~\ref{subfig:legend_efficiency}.
    }%
    \label{fig:interface_efficiency_comparison}
\end{figure}

The case of piecewise constant diffusion coefficient
is included in the analysis of \cite{BringmannDissertation}.
Hence, Theorem~\ref{thm:convergence_separate}
for the optimal convergence rates of SALSFEM
(and so of CALSFEM) generalises
to the elliptic problem~\eqref{eq:diffusion}
as well.
Figure~\ref{subfig:interface_convergence_collective}
confirms the optimal convergence rates
even for a rather large bulk parameter
\(\theta = 0.7\).
The efficiency indices of the built-in error estimator
\(LS(f; p_\ell, u_\ell)^{1/2}\) range from \(1\) to \(1.25\)
in Figure~\ref{subfig:interface_efficiency}
providing further empirical evidence for its
accurate error estimation properties.
The slight increase of these indices
might result from the approximation of inhomogeneous
Dirichlet boundary conditions which are not covered
by \cite{MR3820383}.

\section{Conclusion and open questions}
\label{sec:conclusion}

The numerical experiments show that the adaptive algorithm
with separate marking is superior in particular
on moderate levels and for obtaining an overall high accuracy.
However, the realisation of the separate marking and
of the data approximation algorithm is more involved
and usually not included in standard FEM software
packages.
Since the natural mesh-refinement leads to comparable results
as the separate marking algorithm,
it is a good alternative.
The investigation of the efficiency indices
confirmed the exactness of the built-in error estimator
even on coarse meshes.

Choices of moderate bulk parameters of
\(0.3 \leq \theta \leq 0.5\)
provide optimal convergence rates while still ensuring
a tolerable number of solution steps.
The investigation of the separation parameter \(\kappa\) in
Subsection~\ref{subsec:LshapeMicrostructure}
suggests a choice of \(\kappa\) of one order of magnitude less
than the quotient
\(q_\ell^2 = \mu^2(\T_\ell) / \eta_\textup{S}^2(\T_\ell) \)
in the case of uniform refinement.
The evaluation of \(q_\ell\) for \(\ell = 0\)
or small levels \(\ell > 0\) allow for
a justified a~priori choice of \(\kappa\).
The convergence rate is robust with respect to the
parameter \(\rho\).
Smaller values of \(\rho\) significantly
reduce the number of solution steps
while larger values enable the adaptive algorithm to balance
error estimator reduction and data approximation
more accurately.

The adaptive LSFEM is well-established
and convincing in many applications.
However, important mathematical questions remain open.
The proof of Theorem~\ref{thm:convergence_natural} on
Q-linear convergence of the natural adaptive LSFEM
heavily relies on the lowest-order arguments
such that a straight-forward generalisation
to higher polynomial degrees seems inaccessible.
Additionally, the restriction to sufficiently large bulk
parameters appears artificial and the case of small
\(\theta\) is not covered yet.
Once this has been solved,
the linear convergence would imply optimal convergence rates
with respect to the number of degrees of freedom
by Theorem~\ref{thm:linear_implies_rate}.
Moreover, the study of
optimal convergence rates with respect to the computational costs
in the spirit of \cite{MR4280291}
represents an important task for future research.
While the collective marking algorithm fits into the
framework of \cite{MR4280291},
the application to an adaptive algorithm with separate
marking and data approximation requires a major
modification as for the axioms of adaptivity
in \cite{MR3719030}.

\section*{Acknowledgement}

It is my pleasure to acknowledge fruitful discussions
with Prof.\ Carsten Carstensen and Dr.\ Rui Ma.
This research has been supported by
the Austrian Science Fund (FWF) through the project
\emph{Computational nonlinear PDEs} (grant P33216).

%
%
\hbadness5000
\bibliography{Bibliography}

\begin{thebibliography}{10}

\bibitem{MR2765484}
{\sc J.~H. Adler, T.~A. Manteuffel, S.~F. McCormick, J.~W. Nolting, J.~W. Ruge,
  and L.~Tang}, {\em Efficiency based adaptive local refinement for first-order
  system least-squares formulations}, SIAM J. Sci. Comput., 33 (2011),
  pp.~1--24.

\bibitem{MR2551156}
{\sc F.~S. Attia, Z.~Cai, and G.~Starke}, {\em First-order system least squares
  for the {S}ignorini contact problem in linear elasticity}, SIAM J. Numer.
  Anal., 47 (2009), pp.~3027--3043.

\bibitem{MR1615154}
{\sc M.~Berndt, T.~A. Manteuffel, and S.~F. McCormick}, {\em Local error
  estimates and adaptive refinement for first-order system least squares
  ({FOSLS})}, Electron. Trans. Numer. Anal., 6 (1997), pp.~35--43.
\newblock Special issue on multilevel methods (Copper Mountain, CO, 1997).

\bibitem{MR3817763}
{\sc F.~Bertrand}, {\em First-order system least-squares for interface
  problems}, SIAM J. Numer. Anal., 56 (2018), pp.~1711--1730.

\bibitem{MR4410744}
{\sc F.~Bertrand and D.~Boffi}, {\em First order least-squares formulations for
  eigenvalue problems}, IMA J. Numer. Anal., 42 (2022), pp.~1339--1363.

\bibitem{MR2050077}
{\sc P.~Binev, W.~Dahmen, and R.~DeVore}, {\em Adaptive finite element methods
  with convergence rates}, Numer. Math., 97 (2004), pp.~219--268.

\bibitem{MR2050076}
{\sc P.~Binev and R.~DeVore}, {\em Fast computation in adaptive tree
  approximation}, Numer. Math., 97 (2004), pp.~193--217.

\bibitem{MR2896812}
{\sc M.~Brezina, J.~Garcia, T.~Manteuffel, S.~McCormick, J.~Ruge, and L.~Tang},
  {\em Parallel adaptive mesh refinement for first-order system least squares},
  Numer. Linear Algebra Appl., 19 (2012), pp.~343--366.

\bibitem{BringmannDissertation}
{\sc P.~Bringmann}, {\em Adaptive least-squares finite element method with
  optimal convergence rates}, PhD thesis,  (2021).
\newblock Humboldt-{U}niversit{\"a}t zu {B}erlin.

\bibitem{MR4557622}
\leavevmode\vrule height 2pt depth -1.6pt width 23pt, {\em How to prove optimal
  convergence rates for adaptive least-squares finite element methods}, J.
  Numer. Math., 31 (2023), pp.~43--58.

\bibitem{octAFEM}
\leavevmode\vrule height 2pt depth -1.6pt width 23pt, {\em {octAFEM}}, 2023.
\newblock \textsc{Matlab}/Octave software package, available on Code Ocean.
  DOI: \href{https://doi.org/10.24433/CO.6310426.v1}{10.24433/CO.6310426.v1}.

\bibitem{MR3599566}
{\sc P.~Bringmann and C.~Carstensen}, {\em An adaptive least-squares {FEM} for
  the {S}tokes equations with optimal convergence rates}, Numer. Math., 135
  (2017), pp.~459--492.

\bibitem{MR3715170}
\leavevmode\vrule height 2pt depth -1.6pt width 23pt, {\em {$h$}-adaptive
  least-squares finite element methods for the 2{D} {S}tokes equations of any
  order with optimal convergence rates}, Comput. Math. Appl., 74 (2017),
  pp.~1923--1939.

\bibitem{MR3757107}
{\sc P.~Bringmann, C.~Carstensen, and G.~Starke}, {\em An adaptive
  least-squares {FEM} for linear elasticity with optimal convergence rates},
  SIAM J. Numer. Anal., 56 (2018), pp.~428--447.

\bibitem{MR4433562}
{\sc P.~Bringmann, C.~Carstensen, and N.~T. Tran}, {\em Adaptive least-squares,
  discontinuous {P}etrov-{G}alerkin, and hybrid high-order methods}, in
  Non-standard discretisation methods in solid mechanics, vol.~98 of Lect.
  Notes Appl. Comput. Mech., Springer, Cham, 2022, pp.~107--147.

\bibitem{MR3372049}
{\sc Z.~Cai, V.~Carey, J.~Ku, and E.-J. Park}, {\em Asymptotically exact a
  posteriori error estimators for first-order div least-squares methods in
  local and global {$L_2$} norm}, Comput. Math. Appl., 70 (2015), pp.~648--659.

\bibitem{MR2100304}
{\sc Z.~Cai, J.~Korsawe, and G.~Starke}, {\em An adaptive least squares mixed
  finite element method for the stress-displacement formulation of linear
  elasticity}, Numer. Methods Partial Differential Equations, 21 (2005),
  pp.~132--148.

\bibitem{MR2084237}
{\sc Z.~Cai and G.~Starke}, {\em Least-squares methods for linear elasticity},
  SIAM J. Numer. Anal., 42 (2004), pp.~826--842.

\bibitem{10.1016/j.jnnfm.2009.02.004}
{\sc Z.~Cai and C.~Westphal}, {\em An adaptive mixed least-squares finite
  element method for viscoelastic fluids of oldroyd type}, Journal of
  Non-Newtonian Fluid Mechanics, 159 (2009), pp.~72--80.

\bibitem{MR4011536}
{\sc C.~Carstensen}, {\em Collective marking for adaptive least-squares finite
  element methods with optimal rates}, Math. Comp., 89 (2020), pp.~89--103.

\bibitem{MR3170325}
{\sc C.~Carstensen, M.~Feischl, M.~Page, and D.~Praetorius}, {\em Axioms of
  adaptivity}, Comput. Math. Appl., 67 (2014), pp.~1195--1253.

\bibitem{MR3279489}
{\sc C.~Carstensen, D.~Gallistl, F.~Hellwig, and L.~Weggler}, {\em Low-order
  d{PG}-{FEM} for an elliptic {PDE}}, Comput. Math. Appl., 68 (2014),
  pp.~1503--1512.

\bibitem{MR3824773}
{\sc C.~Carstensen and F.~Hellwig}, {\em Constants in discrete {P}oincar\'{e}
  and {F}riedrichs inequalities and discrete quasi-interpolation}, Comput.
  Methods Appl. Math., 18 (2018), pp.~433--450.

\bibitem{MR4271577}
{\sc C.~Carstensen and R.~Ma}, {\em Collective marking for arbitrary order
  adaptive least-squares finite element methods with optimal rates}, Comput.
  Math. Appl., 95 (2021), pp.~271--281.

\bibitem{AFEMpackage}
{\sc C.~Carstensen and {Numerical Analysis Group}}, {\em {AFEM}}.
\newblock Unpublished \textsc{Matlab} software package, 2009.

\bibitem{MR3296614}
{\sc C.~Carstensen and E.-J. Park}, {\em Convergence and optimality of adaptive
  least squares finite element methods}, SIAM J. Numer. Anal., 53 (2015),
  pp.~43--62.

\bibitem{MR3671598}
{\sc C.~Carstensen, E.-J. Park, and P.~Bringmann}, {\em Convergence of natural
  adaptive least squares finite element methods}, Numer. Math., 136 (2017),
  pp.~1097--1115.

\bibitem{MR2772091}
{\sc C.~Carstensen and H.~Rabus}, {\em An optimal adaptive mixed finite element
  method}, Math. Comp., 80 (2011), pp.~649--667.

\bibitem{MR3719030}
{\sc C.~Carstensen and H.~Rabus}, {\em Axioms of adaptivity with separate
  marking for data resolution}, SIAM J. Numer. Anal., 55 (2017),
  pp.~2644--2665.

\bibitem{MR3820383}
{\sc C.~Carstensen and J.~Storn}, {\em Asymptotic exactness of the
  least-squares finite element residual}, SIAM J. Numer. Anal., 56 (2018),
  pp.~2008--2028.

\bibitem{MR2421046}
{\sc J.~M. Cascon, C.~Kreuzer, R.~H. Nochetto, and K.~G. Siebert}, {\em
  Quasi-optimal convergence rate for an adaptive finite element method}, SIAM
  J. Numer. Anal., 46 (2008), pp.~2524--2550.

\bibitem{MR2870014}
{\sc J.~H. Chaudhry, S.~D. Bond, and L.~N. Olson}, {\em A weighted adaptive
  least-squares finite element method for the {P}oisson-{B}oltzmann equation},
  Appl. Math. Comput., 218 (2012), pp.~4892--4902.

\bibitem{DanischDissertation}
{\sc G.~Danisch}, {\em Gemischte {Finite} {Elemente} {Least-Squares} {Methoden}
  {f\"ur} die {Flachwassergleichung} mit kleiner {Viskosit\"at} (german)
  [{Mixed} least-squares finite element method for the shallow water equation
  with small viscosity]}, PhD thesis,  (2007).
\newblock Gottfried Wilhelm Leibniz {Universit\"at} Hannover.

\bibitem{MR1393904}
{\sc W.~D{\"o}rfler}, {\em A convergent adaptive algorithm for {P}oisson's
  equation}, SIAM J. Numer. Anal., 33 (1996), pp.~1106--1124.

\bibitem{MR1638080}
{\sc J.~M. Fiard, T.~A. Manteuffel, and S.~F. McCormick}, {\em First-order
  system least squares ({FOSLS}) for convection-diffusion problems: numerical
  results}, SIAM J. Sci. Comput., 19 (1998), pp.~1958--1979.

\bibitem{MR4050087}
{\sc T.~F\"{u}hrer}, {\em First-order least-squares method for the obstacle
  problem}, Numer. Math., 144 (2020), pp.~55--88.

\bibitem{MR4425909}
{\sc T.~F\"{u}hrer, N.~Heuer, and M.~Karkulik}, {\em M{INRES} for second-order
  {PDE}s with singular data}, SIAM J. Numer. Anal., 60 (2022), pp.~1111--1135.

\bibitem{MR4242919}
{\sc T.~F\"{u}hrer and M.~Karkulik}, {\em Space-time least-squares finite
  elements for parabolic equations}, Comput. Math. Appl., 92 (2021),
  pp.~27--36.

\bibitem{MR4138307}
{\sc T.~F\"{u}hrer and D.~Praetorius}, {\em A short note on plain convergence
  of adaptive least-squares finite element methods}, Comput. Math. Appl., 80
  (2020), pp.~1619--1632.

\bibitem{MR4280291}
{\sc G.~Gantner, A.~Haberl, D.~Praetorius, and S.~Schimanko}, {\em Rate
  optimality of adaptive finite element methods with respect to overall
  computational costs}, Math. Comp., 90 (2021), pp.~2011--2040.

\bibitem{MR4216839}
{\sc G.~Gantner and R.~Stevenson}, {\em Further results on a space-time {FOSLS}
  formulation of parabolic {PDE}s}, ESAIM Math. Model. Numer. Anal., 55 (2021),
  pp.~283--299.

\bibitem{MR3509205}
{\sc W.~Gautschi}, {\em Orthogonal polynomials in {MATLAB}}, vol.~26 of
  Software, Environments, and Tools, Society for Industrial and Applied
  Mathematics (SIAM), Philadelphia, PA, 2016.
\newblock Exercises and solutions.

\bibitem{MR0245201}
{\sc G.~H. Golub and J.~H. Welsch}, {\em Calculation of {G}auss quadrature
  rules}, Math. Comp. 23 (1969), 221-230; addendum, ibid., 23 (1969),
  pp.~A1--A10.

\bibitem{MR3103833}
{\sc H.~Gu and H.~Li}, {\em An adaptive least-squares mixed finite element
  method for nonlinear parabolic problems}, Comput. Math. Model., 20 (2009),
  pp.~192--206.

\bibitem{MR0461948}
{\sc D.~C. Jespersen}, {\em A least squares decomposition method for solving
  elliptic equations}, Math. Comp., 31 (1977), pp.~873--880.

\bibitem{10.1002/nme.1620240308}
{\sc B.-N. Jiang and G.~F. Carey}, {\em Adaptive refinement for least-squares
  finite elements with element-by-element conjugate gradient solution},
  International Journal for Numerical Methods in Engineering, 24 (1987),
  pp.~569--580.

\bibitem{MR3097045}
{\sc M.~Karkulik, D.~Pavlicek, and D.~Praetorius}, {\em On 2{D} newest vertex
  bisection: optimality of mesh-closure and {$H^1$}-stability of
  {$L_2$}-projection}, Constr. Approx., 38 (2013), pp.~213--234.

\bibitem{KayserHeroldDissertation}
{\sc O.~Kayser-Herold}, {\em Least-{Squares} {Methods} for the {Solution} of
  {Fluid-Structure} {Interaction} {Problems}}, PhD thesis,  (2006).
\newblock Technische Universit\"at Braunschweig.

\bibitem{MR0393815}
{\sc R.~B. Kellogg}, {\em On the {P}oisson equation with intersecting
  interfaces}, Applicable Anal., 4 (1974/75), pp.~101--129.

\bibitem{MR1329875}
{\sc I.~Kossaczk\'{y}}, {\em A recursive approach to local mesh refinement in
  two and three dimensions}, J. Comput. Appl. Math., 55 (1994), pp.~275--288.

\bibitem{MR3580409}
{\sc R.~Krause, B.~M\"{u}ller, and G.~Starke}, {\em An adaptive least-squares
  mixed finite element method for the {S}ignorini problem}, Numer. Methods
  Partial Differential Equations, 33 (2017), pp.~276--289.

\bibitem{MR3452849}
{\sc J.~Ku}, {\em Local error estimates for least-squares finite element
  methods for first-order system}, J. Comput. Appl. Math., 299 (2016),
  pp.~92--100.

\bibitem{MR2671052}
{\sc J.~Ku and E.-J. Park}, {\em A posteriori error estimators for the
  first-order least-squares finite element method}, J. Comput. Appl. Math., 235
  (2010), pp.~293--300.

\bibitem{MR1793582}
{\sc J.-L. Liu}, {\em Exact a posteriori error analysis of the least squares
  finite element method}, Appl. Math. Comput., 116 (2000), pp.~297--305.

\bibitem{MR4127415}
{\sc Q.~Liu and S.~Zhang}, {\em Adaptive flux-only least-squares finite element
  methods for linear transport equations}, J. Sci. Comput., 84 (2020),
  pp.~Paper No. 26, 22.

\bibitem{MR4087177}
{\sc Q.~Liu and S.~Zhang}, {\em Adaptive least-squares finite element methods
  for linear transport equations based on an {$\rm H(div)$} flux
  reformulation}, Comput. Methods Appl. Mech. Engrg., 366 (2020), pp.~113041,
  25.

\bibitem{MR1885710}
{\sc M.~Majidi and G.~Starke}, {\em Least-squares {G}alerkin methods for
  parabolic problems. {II}. {T}he fully discrete case and adaptive algorithms},
  SIAM J. Numer. Anal., 39 (2001/02), pp.~1648--1666.

\bibitem{MR2650218}
{\sc T.~Manteuffel, S.~McCormick, J.~Nolting, J.~Ruge, and G.~Sanders}, {\em
  Further results on error estimators for local refinement with first-order
  system least squares ({FOSLS})}, Numer. Linear Algebra Appl., 17 (2010),
  pp.~387--413.

\bibitem{MR1311687}
{\sc J.~M. Maubach}, {\em Local bisection refinement for {$n$}-simplicial grids
  generated by reflection}, SIAM J. Sci. Comput., 16 (1995), pp.~210--227.

\bibitem{MR1770058}
{\sc P.~Morin, R.~H. Nochetto, and K.~G. Siebert}, {\em Data oscillation and
  convergence of adaptive {FEM}}, SIAM J. Numer. Anal., 38 (2000),
  pp.~466--488.

\bibitem{MR3343602}
{\sc S.~M\"{u}nzenmaier}, {\em First-order system least squares for
  generalized-{N}ewtonian coupled {S}tokes-{D}arcy flow}, Numer. Methods
  Partial Differential Equations, 31 (2015), pp.~1150--1173.

\bibitem{MR2783231}
{\sc S.~M\"{u}nzenmaier and G.~Starke}, {\em First-order system least squares
  for coupled {S}tokes-{D}arcy flow}, SIAM J. Numer. Anal., 49 (2011),
  pp.~387--404.

\bibitem{OlsonDissertation}
{\sc L.~N. Olson}, {\em Multilevel {Least-Squares} {Finite} {Element} {Methods}
  for {Hyperbolic} {Partial} {Differential} {Equations}}, PhD thesis,  (2003).
\newblock University of Colorado.

\bibitem{MR4136545}
{\sc C.-M. Pfeiler and D.~Praetorius}, {\em D\"{o}rfler marking with minimal
  cardinality is a linear complexity problem}, Math. Comp., 89 (2020),
  pp.~2735--2752.

\bibitem{MR4173220}
{\sc W.~Qiu and S.~Zhang}, {\em Adaptive first-order system least-squares
  finite element methods for second-order elliptic equations in nondivergence
  form}, SIAM J. Numer. Anal., 58 (2020), pp.~3286--3308.

\bibitem{Rabus2015}
{\sc H.~Rabus}, {\em Quasi-optimal convergence of {AFEM} based on separate
  marking, {P}art {I}}, J. Numer. Math., 23 (2015), pp.~137--156.

\bibitem{MR3377430}
{\sc H.~Rabus}, {\em Quasi-optimal convergence of {AFEM} based on separate
  marking, {P}art {I} and {II}}, J. Numer. Math., 23 (2015), pp.~137--156.

\bibitem{MR2832786}
{\sc K.~G. Siebert}, {\em A convergence proof for adaptive finite elements
  without lower bound}, IMA J. Numer. Anal., 31 (2011), pp.~947--970.

\bibitem{MR2139398}
{\sc G.~Starke}, {\em A first-order system least squares finite element method
  for the shallow water equations}, SIAM J. Numer. Anal., 42 (2005),
  pp.~2387--2407.

\bibitem{MR2285860}
\leavevmode\vrule height 2pt depth -1.6pt width 23pt, {\em An adaptive
  least-squares mixed finite element method for elasto-plasticity}, SIAM J.
  Numer. Anal., 45 (2007), pp.~371--388.

\bibitem{MR2324418}
{\sc R.~Stevenson}, {\em Optimality of a standard adaptive finite element
  method}, Found. Comput. Math., 7 (2007), pp.~245--269.

\bibitem{MR2353951}
\leavevmode\vrule height 2pt depth -1.6pt width 23pt, {\em The completion of
  locally refined simplicial partitions created by bisection}, Math. Comp., 77
  (2008), pp.~227--241.

\bibitem{SH74}
{\sc I.~E. Sutherland and G.~W. Hodgman}, {\em Reentrant polygon clipping},
  Commun. ACM, 17 (1974), p.~32–42.

\bibitem{MR1475530}
{\sc C.~T. Traxler}, {\em An algorithm for adaptive mesh refinement in {$n$}
  dimensions}, Computing, 59 (1997), pp.~115--137.

\end{thebibliography}

\end{document}